# 3x + 1 Problem. Rev-16

**Vicente Padilla**
**April 2025**
vpadilla@aertecsolutions.com

## Abstract

The 3x+1 Problem, also known as Collatz Problem or Syracuse Problem, is stated as follows:

Let any positive integer $a_1$ be the first element of a sequence S

If $a_1$ is even, $a_1 = 2P$; then the next element of the sequence is

$$a_2 = \frac{a_1}{2}$$

If $a_1$ is odd, $a_1 = 2P + 1$; then the next element of the sequence is

$$a_2 = 3 \cdot a_1 + 1$$

If we repeat the same process indefinitely, the sequence eventually reaches "1".

$a_1$, the initial element of the sequence, is obviously a finite number.

$$a_1 < 2^{k_{00}}$$

In this paper we prove that the conjecture is true, and all sequences converge to 1. We first prove that all sequences can be broken up in many different cycles. Each cycle follows the same pattern:

1) Upward trajectory. Odd and even numbers alternate – upward and downward steps - until the cycle reaches an upper bound.
2) Downward trajectory. Two or more consecutive even numbers follow – downward steps only - until it reaches another odd number.

At this point, it's the beginning of the following cycle. All cycles begin and end with an odd number. It is important to note that the final number of a given cycle can be larger or smaller than the initial number of the cycle. Any sequence is evidently made of many consecutive cycles.

We also prove that, for any combination of cycles made of different upward and downward steps, it is always possible to build a Sequence-D made of said combination of cycles whose last cycle converges to 1.

In other words, Sequence-A and Sequence-D follow the same pattern of upward and downward steps except for the very end.

Finally, using Sequence-D, we prove that, if the initial Sequence-A unfolds indefinitely, the initial element of Sequence-A must be larger than $2^{k_{00}}$. This is a contradiction since the initial element is always a finite element. As a result, Sequence-A must necessarily converge to 1.

## 1. Definitions

**Sequence:** The collection of all the numbers, from the initial element of the sequence to the last one, which must obviously be 1.

**Step:** One single move forward in the sequence. The step can be upwards or downwards.

**Subscript:** Subscript "o" means ODD number.

$a_O^1$ - **Initial Element of the Sequence.** It is the initial element that generates the sequence. The initial element is finite.

$$a_O^1 < 2^{K_{00}}$$

**Cycle.** All sequences can be broken up in cycles. Each cycle follows the same pattern:

1) **Upward trajectory**. Odd and even numbers alternate until the cycle reaches an upper bound
2) **Downward trajectory**. Two or more consecutive even numbers follow until it reaches another odd number

The initial and the final elements of each cycle are ODD numbers. The final element of one cycle becomes obviously the initial element of the following cycle.

$a_O^s$ - **Initial Element of Cycle "s".** It is the initial element that generates the Cycle "s".

$a_F^s$ - **Final Element of Cycle "s".** It is the final element of the Cycle "s".

**Cycle Upper Bound.** The element where a cycle reaches its maximum height. Each cycle has its own upper bound. It is always an EVEN number.

**Closing Cycle.** It is the last and final cycle that makes the sequence converge to 1.

**Grade "n" of a Cycle.** It is the number of upward steps in a cycle.

**α**. It is the number of downward steps in a cycle. It includes the downward steps during both, the upward and the trajectories. Obviously, α > n

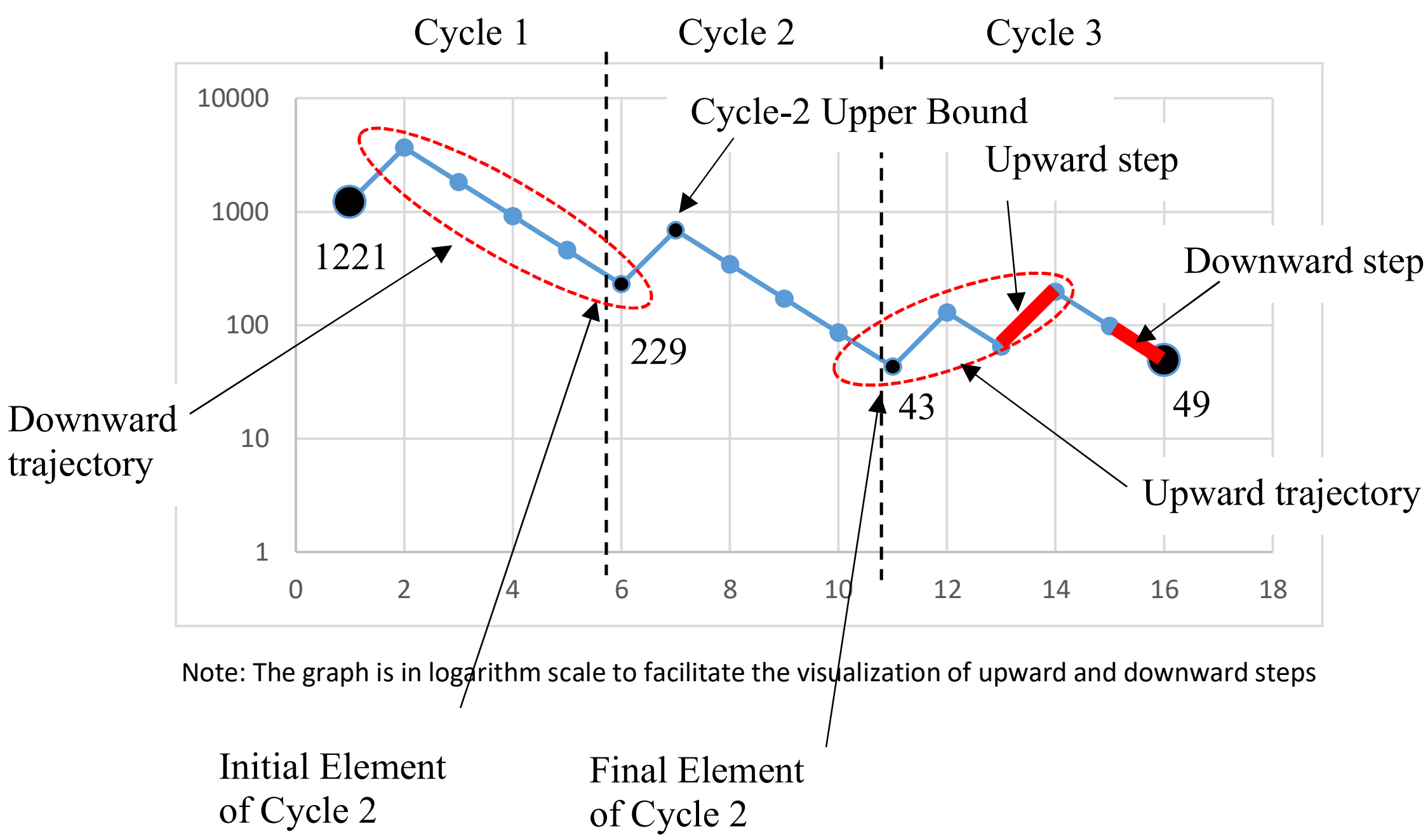

Note: The graph is in logarithm scale to facilitate the visualization of upward and downward steps

## 2. Introduction

This document is a revision of the previous document.

**3x + 1 Problem. Rev-14**
**October 2019**
**vpadilla@aertecsolutions.com**

There are significant changes to this document from the previous one. Nevertheless, the identities that constitute the foundation of the demonstration remain the same. The main change is that the final proof of the conjecture has been modified and simplified.

In **Appendix 1** we prove that any odd number can be expressed as:

$$(A1.1) \qquad a_O = 2^{\alpha} \cdot Q_O - \frac{2^{3^{n-1} \cdot j_o^{\beta}} + 1}{3^n} \cdot (3^n - 2^n)$$

$$n, \alpha \in \mathbb{N}$$

$$Q_O, j_o^{\beta} \in \{odd\}$$

In **Appendix 2** we define "n" as the number of upward steps until the first cycle of the sequence reaches an upper bound. "α" represents the number of downward steps until this first cycle reaches its end. As we shall see, in every cycle the number of upward steps is smaller than the number of downward steps. As a result,

$$\alpha > n$$

$$\forall \, a_O$$

In **Appendix 3** we prove that the sequence is made up of many cycles, and each cycle follows the same pattern. First there is an upward trajectory; odd and even numbers alternate until the cycle reaches an upper bound. Second, there is a downward trajectory; two or more consecutive even numbers follow until the sequence reaches another odd number.

At this point, it's the beginning of the following cycle. It is important to note that the final number of a given cycle can be larger or smaller than the initial number of the cycle. Any sequence is evidently made of many consecutive cycles.

Each cycle of a sequence starts and ends with an odd number. The expression for the initial odd number of the sequence is:

$$(A1.1) \qquad a_O = 2^{\alpha} \cdot Q_O - \frac{2^{3^{n-1} \cdot j_o^{\beta}} + 1}{3^n} \cdot (3^n - 2^n)$$

$$n, \alpha \in \mathbb{N}$$

$$Q_O, j_o^{\beta} \in \{odd\}$$

$$3^{n-1} \cdot j_o^{\beta} > \alpha > n$$

The final odd number at the end of the downward trajectory can be expressed as:

$$(A3.2) \qquad a_F = 3^n \cdot Q_O - 2^{3^{n-1} \cdot j_o^{\beta} - \alpha} \cdot (3^n - 2^n)$$

There is an example of a cycle in Appendix 3.

Since $a_F$ is an odd number, the sequence is followed by a new cycle. It climbs up again – upward trajectory – until it reaches another upper bound – even number - and then, two or more even numbers follow, and the sequence descends – downward trajectory - until it meets another odd number. And a new cycle begins.

In **Appendix 4** we prove that all fractions used in the identities of this document are always integer numbers. All numbers are integer numbers.

In **Appendix 5** we prove that all odd numbers can also be written using the following identity:

$$a_O = 3^m \cdot K_O - 2^n \cdot \frac{2^{3^{n-1}\cdot q_o}+1}{3^n}$$

$$\forall\, m, n \in \mathbb{N}$$

$$K_O, q_O \in \{odd\}$$

In **Appendix 6** we combine "i" cycles together. The initial element - odd number – of a sequence made of "i" cycles can be expressed as:

$$(A6.1) \quad a_O^1 = 2^{\sum_{t=1}^{t=i}\alpha_t} \cdot Q_{O.i} - \frac{2^{3^{\sum_{t=1}^{t=i} n_t - 1}\cdot j_{o.i}^{\delta}}+1}{3^{\sum_{t=1}^{t=i} n_t}} \cdot \sum_{s=1}^{i} 3^{\sum_{t=s+1}^{t=i} n_t} \cdot 2^{\sum_{t=1}^{t=s-1}\alpha_t} \cdot (3^{n_s} - 2^{n_s})$$

$$Q_{O.i}; j_{o.i}^{\delta} \in \{odd\}$$

where the above parameters $n_s, \alpha_s$ are the standard parameters of any Cycle-s as described in Appendix 3. $n_s$ represents the number of upward steps for Cycle-s and $\alpha_s$ represents the number of downward steps for said Cycle-s. The identity above also has the same restrictions as described in Appendix 1 and Appendix 2.

$$\alpha_s > n_s \qquad \forall s$$

$$(A6.2) \quad 3^{\sum_{t=1}^{t=i} n_s - 1} \cdot j_{o.i}^{\delta} > \sum_{s=1}^{s=i} \alpha_s$$

We also prove that the final element of this sequence after "$i$" cycles is:

$$(A6.3) \quad a_F^i = 3^{\sum_{t=1}^{t=i} n_t} \cdot Q_{O.i} - 2^{3^{\sum_{t=1}^{t=i} n_t - 1}\cdot j_{o.i}^{\delta} - \sum_{t=1}^{t=i}\alpha_t} \cdot \sum_{s=1}^{i} 3^{\sum_{t=s+1}^{t=i} n_t} \cdot 2^{\sum_{t=1}^{t=s-1}\alpha_t} \cdot (3^{n_s} - 2^{n_s})$$

In other words, $a_O^1$ generates a sequence, and after $\sum_{t=1}^{t=i} n_t$ upward steps and $\sum_{t=1}^{t=i} \alpha_t$downward steps, during "i" cycles, it reaches $a_F^i$. There is an example of a 3-cycle sequence in Appendix 6.

In **Appendix 7** we prove that, for any combination of cycles as shown in Appendix 6, it always exists a sequence that follows the exact combination of said cycles but whose last cycle converges to 1. In this document we always use the letter "d" or "D" to refer to an odd number $d_O^1$ that generates a Sequence-D that converges to 1 in its final Cycle-k.

The initial element of Sequence-D is:

$$(A7.6) \quad d_O^1 = 2^{\sum_{t=1}^{t=k-1} \alpha_t + 3^{n_k - 1} \cdot q_o + n_k} \cdot Q_{O.d} - \frac{2^{3^{\sum_{t=1}^{t=k} n_t - 1} \cdot j_{o.k}^{\delta}} + 1}{3^{\sum_{t=1}^{t=k} n_t}} \cdot \sum_{s=1}^{k} 3^{\sum_{t=s+1}^{t=k} n_t} \cdot 2^{\sum_{t=1}^{t=s-1} \alpha_t} \cdot (3^{n_s} - 2^{n_s})$$

And its final element after “k” cycles is obviously 1.

$$d_F^k = 1$$

$q_o$ and $Q_{O.d}$ are mutually dependent variables.

In **Appendix 8** we introduce the concept of parallel paths, where two sequences follow the exact pattern of upward and downward steps during a given number of downward steps. We prove that only one number smaller than $2^{\beta}$ can generate a sequence that follows its own pattern during $\beta$ downward steps.
Example:

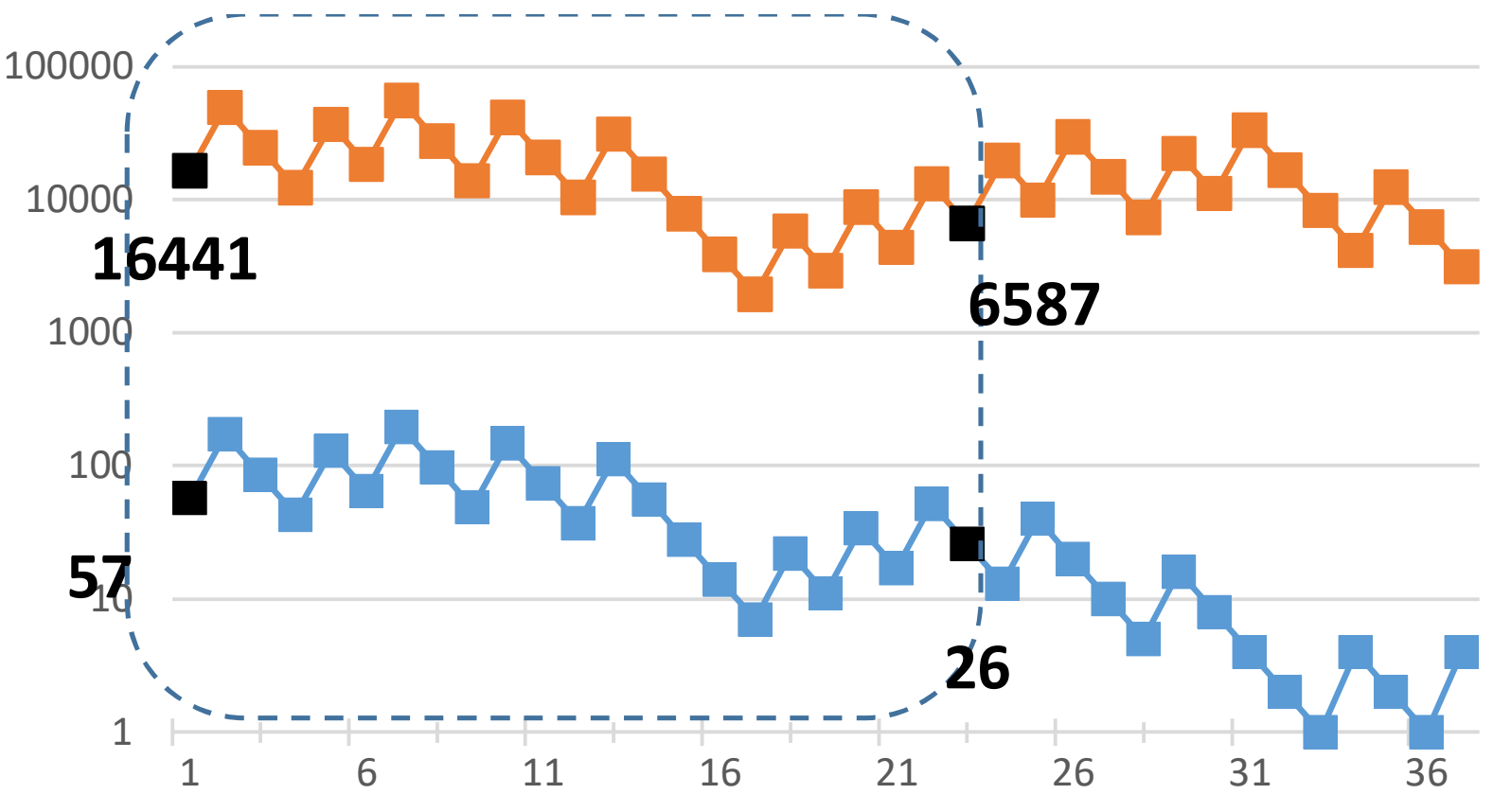

*Ordinate axis in logarithmic scale to facilitate visualization*

In **Appendix 9** we prove several relevant lemmas needed for the final proof.

**Single-cycle Sequence**

In Appendix 3, we also prove that any sequence generated by the following odd number:

$$(A3.7) \quad a_O = 2^n \cdot \frac{2^{3^{n-1} \cdot q_o} + 1}{3^n} - 1$$

$$q_o \in \{odd\}$$

converges to “1” in one single cycle. In other words:

$$a_F = 1$$

This single-cycle sequence follows the same pattern as the others. There is first an upward trajectory in which “n” upward and downward steps follow, until it reaches an upper bound – an even number - and at this stage, the sequence follows a downward trajectory made of "$3^{n-1} \cdot q_o + 1$" downward steps until it reaches an odd number, in this case 1.

There is an example of a single-cycle sequence in Appendix 3.

## 3. Sequence-A and Sequence-D

To prove the conjecture, two sequences will be built. Sequence-A is any sequence that we need to prove it converges to 1. Sequence-D is a sequence similar to Sequence-A but converges to 1. This is how we do it.

Sequence-A is generated by any odd number $a_O^1$. Obviously, this initial number $a_O^1$ is a finite number.

$$a_O^1 < 2^{k_{00}}$$

$$k_{00} \in \mathbb{N}$$

Please note that the double "0" in the subscript of k is used to denote that it is a fixed number.

If after "t" cycles, all final elements $a_F^i$ of each cycle were smaller than the initial element $a_O^1$ of the sequence,

$$a_O^1 > a_F^i$$

$$\forall i > t$$

Then, the sequence would converge, and the conjecture would be proven. Therefore, from now on, we assume that Sequence-A unfolds indefinitely and does not converge. In other words, after "t" cycles, all final elements $a_F^i$ of each cycle are equal to or larger than the initial element $a_O^1$ of the sequence,

$$a_O^1 \leq a_F^i$$

$$\forall i > t$$

As we have seen, we prove in Appendix 3 that the sequence can be broken down in many different cycles.

We then build another Sequence-D that follows a parallel path to the original one. The initial element of Sequence-D is $d_O^1$. A parallel path means that both sequencies follow the same pattern of upward and downward steps. Appendix 7 proves that this is always possible. Sequence-D is built in such a way that both sequences - A and D - follow parallel paths during k cycles until the number of downward steps is larger than $k_{00}$. Please be aware that since $a_O^1$ and $k_{00}$ are finite numbers

$$a_O^1 < 2^{k_{00}}$$

it is always possible to find a Cycle-k for which:

$$\sum_{s=1}^{s=k} \alpha_s > k_{00}$$

At that point, the sequences' paths diverge and follow different patterns. Sequence-A keeps unfolding and Sequence-D directly converges to 1. Once Sequence-D reaches 1, its following cycles are always the trivial cycle 4-2-1. In other words, after the Cycle-k, Sequence-D is only made of cycles with one step up and two steps down:

$$n_s = 1$$
$$\alpha_s = 2$$
$$\forall\, s \geq k+1$$

The initial element of Sequence-A - $a_O^1$ – must be smaller than the initial element of Sequence-D - $d_O^1$.

$$a_O^1 < d_O^1$$

This is the case since Appendix 8 proves that if

$$a_O^1 < 2^{k_{00}} < 2^{\sum_{s=1}^{s=k} \alpha_s}$$

$a_O^1$ is the only number smaller than $2^{K_{00}}$ whose sequence follows its own pattern during at least $k_{00}$ downward steps. Since Sequence-A and Sequence-D follow parallel paths (the same pattern of upward and downward steps) during the first k cycles, and as we just saw above

$$\sum_{s=1}^{s=k} \alpha_s > k_{00}$$

Also, since the initial assumption is that

$$a_O^1 < 2^{k_{00}}$$

then, we must conclude that $\mathrm{d}_\mathrm{O}^1$ must be larger than $2^{\mathrm{k}_{00}}$ and $2^{\sum_{s=1}^{s=k} \alpha_s}$, and as a result, larger than $a_O^1$.

$$(3.1) \quad d_O^1 > 2^{\sum_{s=1}^{s=k} \alpha_s} > 2^{k_{00}} > a_O^1$$

The graph below shows an example of two sequences that follow parallel paths during 2 cycles. Sequence-A keeps unfolding while Sequence-D converges to 1 in Cycle-2.

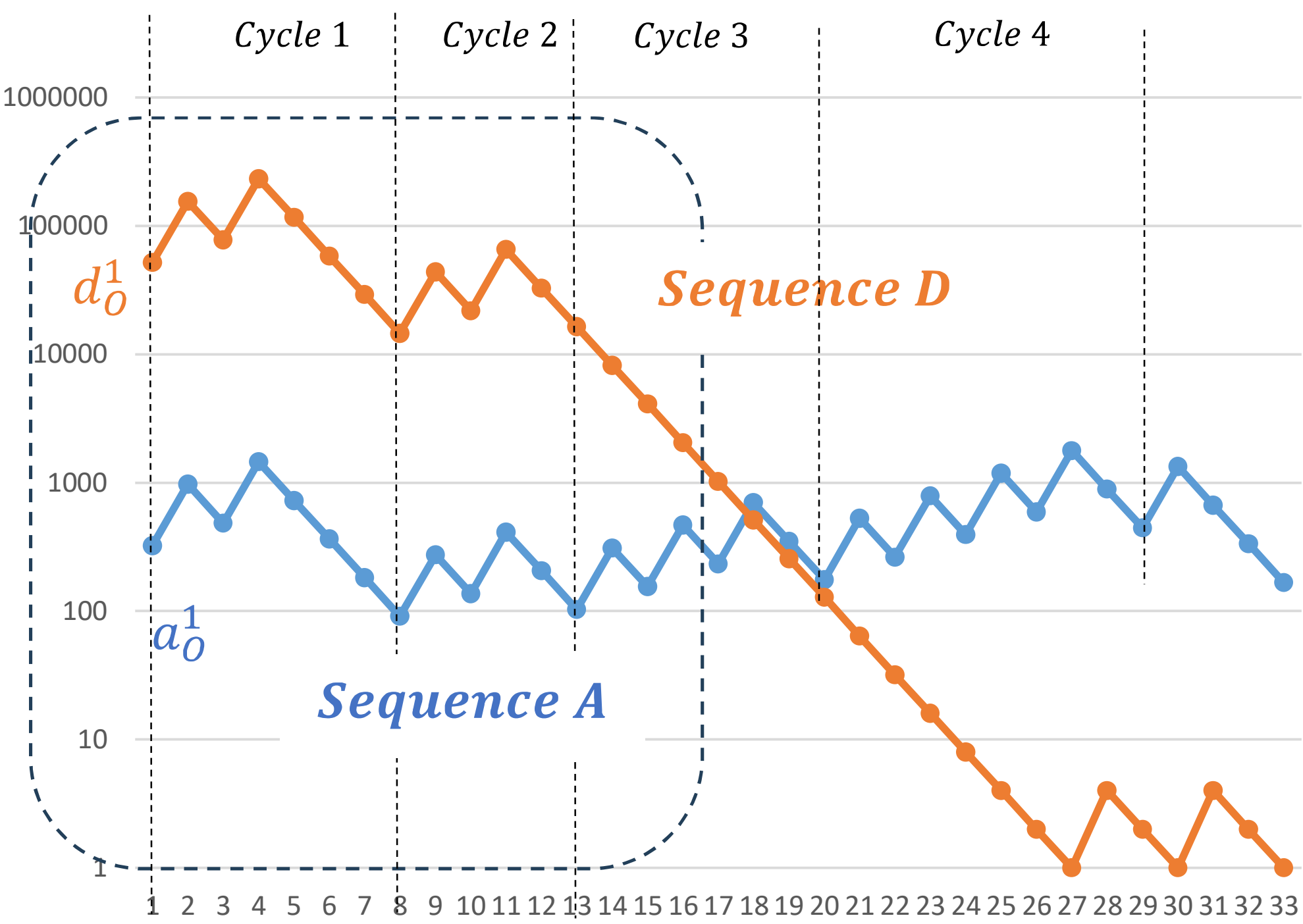

*Ordinate axis is in logarithmic scale so improve visualization.*

Please note that all parameters $\alpha_t$ and $n_t$ are identical except for the final number of downward steps $\alpha_k$ for Sequence-D's Cycle-k. The final number of downward steps for this final Cycle-k is given by Identity (A3.8) of Appendix 3.

$$(A3.8) \quad (\alpha_k)_{Sequence\ D} = 3^{n_k - 1} \cdot q_O + n_k$$
$$q_O \in \{odd\}$$

After Cycle-k, Sequence-A continues unfolding without converging.

Using Identity (A6.1) from Appendix 6, the initial element of Sequence-A can be expressed as:

$$(3.2) \quad a_O^1 = 2^{\sum_{t=1}^{t=i} \alpha_t} \cdot Q_{O.i}^A - \frac{2^{3^{\sum_{t=1}^{t=i} n_t - 1} \cdot j_{o.i}^{\delta}} + 1}{3^{\sum_{t=1}^{t=i} n_t}} \cdot \sum_{s=1}^{i} 3^{\sum_{t=s+1}^{t=i} n_t} \cdot 2^{\sum_{t=1}^{t=s-1} \alpha_t} \cdot (3^{n_s} - 2^{n_s})$$

We divide the summation above into three blocks. The First Block contains the first k cycles of Sequence-A until

$$\sum_{s=1}^{s=k} \alpha_s > k_{00}$$

The Second Block contains from Cycle-k+1 until the number of downward steps of Sequence-A matches the number of downward steps of Sequence-D at Cycle-k. Please, be reminded that Sequence-D follows a direct downward trajectory – only downward steps - until it converges to 1. The end of the Second Block in Sequence-A is reached at Cycle-j.

The Third Block of Sequence-A contains all cycles onward, in other words, from Cycle-j+1 to Cycle-i.

Therefore. Identity (3.2) becomes:

$$a_O^1 = 2^{\sum_{t=1}^{t=i}\alpha_t} \cdot Q_{O.i}^A - \frac{2^{3^{\sum_{t=1}^{t=i} n_t - 1} \cdot j_{o.i}^{\delta}} + 1}{3^{\sum_{t=1}^{t=i} n_t}} \cdot \sum_{s=1}^{k} 3^{\sum_{t=s+1}^{t=i} n_t} \cdot 2^{\sum_{t=1}^{t=s-1}\alpha_t} \cdot (3^{n_s} - 2^{n_s})$$

$$- \frac{2^{3^{\sum_{t=1}^{t=i} n_t - 1} \cdot j_{o.i}^{\delta}} + 1}{3^{\sum_{t=1}^{t=i} n_t}} \cdot \sum_{s=k+1}^{j} 3^{\sum_{t=s+1}^{t=i} n_t} \cdot 2^{\sum_{t=1}^{t=s-1}\alpha_t} \cdot (3^{n_s} - 2^{n_s})$$

$$- \frac{2^{3^{\sum_{t=1}^{t=i} n_t - 1} \cdot j_{o.i}^{\delta}} + 1}{3^{\sum_{t=1}^{t=i} n_t}} \cdot \sum_{s=j+1}^{i} 3^{\sum_{t=s+1}^{t=i} n_t} \cdot 2^{\sum_{t=1}^{t=s-1}\alpha_t} \cdot (3^{n_s} - 2^{n_s})$$

We extract the common factors $3^{\sum_{t=k+1}^{t=j} n_t}$ and $3^{\sum_{t=j+1}^{t=i} n_t}$ from the first summation, the common factors $2^{\sum_{t=1}^{t=k}\alpha_t}$ and $3^{\sum_{t=j+1}^{t=i} n_t}$ from the second summation and the common factors $2^{\sum_{t=1}^{t=k}\alpha_t}$ and $2^{\sum_{t=k+1}^{t=j}\alpha_t}$ from the third summation. We obtain:

$$(3.3) \qquad a_O^1 = 2^{\sum_{t=1}^{t=i}\alpha_t} \cdot Q_{O.i}^A +$$

$$- \frac{2^{3^{\sum_{t=1}^{t=i} n_t - 1} \cdot j_{o.i}^{\delta}} + 1}{3^{\sum_{t=1}^{t=i} n_t}} \cdot 3^{\sum_{t=k+1}^{t=j} n_t} \cdot 3^{\sum_{t=j+1}^{t=i} n_t} \cdot \sum_{s=1}^{k} 3^{\sum_{t=s+1}^{t=k} n_t} \cdot 2^{\sum_{t=1}^{t=s-1}\alpha_t} \cdot (3^{n_s} - 2^{n_s})$$

$$- \frac{2^{3^{\sum_{t=1}^{t=i} n_t - 1} \cdot j_{o.i}^{\delta}} + 1}{3^{\sum_{t=1}^{t=i} n_t}} \cdot 2^{\sum_{t=1}^{t=k}\alpha_t} \cdot 3^{\sum_{t=j+1}^{t=i} n_t} \sum_{s=k+1}^{j} 3^{\sum_{t=s+1}^{t=j} n_t} \cdot 2^{\sum_{t=k+1}^{t=s-1}\alpha_t} \cdot (3^{n_s} - 2^{n_s})$$

$$- \frac{2^{3^{\sum_{t=1}^{t=i} n_t - 1} \cdot j_{o.i}^{\delta}} + 1}{3^{\sum_{t=1}^{t=i} n_t}} \cdot 2^{\sum_{t=1}^{t=k}\alpha_t} \cdot 2^{\sum_{t=k+1}^{t=j}\alpha_t} \sum_{s=j+1}^{i} 3^{\sum_{t=s+1}^{t=i} n_t} \cdot 2^{\sum_{t=j+1}^{t=s-1}\alpha_t} \cdot (3^{n_s} - 2^{n_s})$$

We apply the same logic to the first element of Sequence-D. We start with Identity (A7.6) from Appendix 7, -, but instead of unfolding “i” cycles as in Sequence-A, we unfold Sequence-D “m” cycles. This is the case since, once Sequence-D reaches Cycle-k, Sequence-A and Sequence-D unfold at different speeds.

$$(3.4) \qquad d_O^1 = 2^{\sum_{t=1}^{t=k-1}\alpha_t + 3^{n_k - 1} \cdot q_o + n_k + \sum_{t=k+1}^{t=m}\alpha_t} \cdot Q_{O.m} +$$

$$- \frac{2^{3^{\sum_{t=1}^{t=m} n_t - 1} \cdot j_{o.m}^{\delta}} + 1}{3^{\sum_{t=1}^{t=m} n_t}} \cdot \sum_{s=1}^{m} 3^{\sum_{t=s+1}^{t=m} n_t} \cdot 2^{\sum_{t=1}^{t=s-1}\alpha_t} \cdot (3^{n_s} - 2^{n_s})$$

Please, be reminded, that the parameters $\alpha_t$ and $n_t$ for Sequence-D after Cycle-k – after converging to 1 - are no longer the same as in Sequence-A but the ones corresponding to the trivial cycle 4-2-1:

$$\alpha_t = 2$$
$$n_t = 1$$
$$\forall t \geq k+1$$

As a result, for Sequence-D:

$$\sum_{t=k+1}^{t=m} \alpha_t = \sum_{t=k+1}^{t=m} 2 = 2 \cdot (m-k)$$
$$\sum_{t=k+1}^{t=m} n_t = \sum_{t=k+1}^{t=m} 1 = m-k$$

Therefore, we can divide the summation of Sequence-D into two blocks. The First Block - as in Sequence-A – contains the first "k" cycles, nevertheless, Sequence-D converges to 1 in Cycle-k. The Second Block, once the sequence has already converged to 1, the trivial cycle 4-2-1 repeats indefinitely. As a result, after "m" cycles, the initial element of Sequence-D can be expressed as:

$$(3.5) \qquad d_O^1 = 2^{\sum_{t=1}^{t=k-1} \alpha_t + 3^{n_k - 1} \cdot q_o + n_k + 2 \cdot (m-k)} \cdot Q_{O.m}^D +$$
$$- \frac{2^{3^{\sum_{t=1}^{t=k} n_t + m - k - 1} \cdot j_{O.m}^{\delta}} + 1}{3^{\sum_{t=1}^{t=k} n_t + m - k}} \cdot 3^{(m-k)} \cdot \sum_{s=1}^{k} 3^{\sum_{t=s+1}^{t=k} n_t} \cdot 2^{\sum_{t=1}^{t=s-1} \alpha_t} \cdot (3^{n_s} - 2^{n_s})$$
$$- \frac{2^{3^{\sum_{t=1}^{t=k} n_t + m - k - 1} \cdot j_{O.m}^{\delta}} + 1}{3^{\sum_{t=1}^{t=k} n_t + m - k}} \cdot 2^{\sum_{t=1}^{t=k} \alpha_t} \cdot \sum_{s=k+1}^{m} 3^{\sum_{t=s+1}^{t=m} 1} \cdot 2^{\sum_{t=k+1}^{t=s-1} 2} \cdot (3^1 - 2^1)$$

In Appendix 9, Lemma A9.4 proves that the second summation of the identity above can be written in a simpler way:

$$(A9.12) \quad \sum_{s=k+1}^{m} 3^{\sum_{t=s+1}^{t=m} 1} \cdot 2^{\sum_{t=k+1}^{t=s-1} 2} \cdot (3^1 - 2^1) = 4^{(m-k)} - 3^{(m-k)}$$

We then apply Identity (A9.12) to Identity (3.5) and obtain

$$(3.6) \qquad d_O^1 = 2^{\sum_{t=1}^{t=k-1} \alpha_t + 3^{n_k - 1} \cdot q_o + n_k + 2 \cdot (m-k)} \cdot Q_{O.m}^D +$$
$$- \frac{2^{3^{\sum_{t=1}^{t=k} n_t + m - k - 1} \cdot j_{O.m}^{\delta}} + 1}{3^{\sum_{t=1}^{t=k} n_t + m - k}} \cdot 3^{(m-k)} \cdot \sum_{s=1}^{k} 3^{\sum_{t=s+1}^{t=k} n_t} \cdot 2^{\sum_{t=1}^{t=s-1} \alpha_t} \cdot (3^{n_s} - 2^{n_s})$$
$$- \frac{2^{3^{\sum_{t=1}^{t=k} n_t + m - k - 1} \cdot j_{O.m}^{\delta}} + 1}{3^{\sum_{t=1}^{t=k} n_t + m - k}} \cdot 2^{\sum_{t=1}^{t=k} \alpha_t} \cdot \left[4^{(m-k)} - 3^{(m-k)}\right]$$

In order to simplify the Identity (3.3) and Identity (3.6) above, we use the following **Nomenclature**:

## Nomenclature

-

$$\alpha_1^k = \sum_{t=1}^{t=k} \alpha_t \quad ; \quad \alpha_k^j = \sum_{t=k+1}^{t=j} \alpha_t \quad ; \quad \alpha_j^i = \sum_{t=j+1}^{t=i} \alpha_t$$

$$n_1^k = \sum_{t=1}^{t=k} n_t \quad ; \quad n_k^j = \sum_{t=k+1}^{t=j} n_t \quad ; \quad n_j^i = \sum_{t=j+1}^{t=i} n_t$$

$$\Delta = 3^{n_k - 1} \cdot q_O + n_k - \alpha_k$$

$$w_i = 3^{\sum_{t=1}^{t=i} n_t - 1} \cdot j_{o.i}^{\delta}$$

$$w_m = 3^{\left(\sum_{t=1}^{t=k} n_t + m - k - 1\right)} \cdot j_{o.m}^{\delta}$$

$$\sum K = \sum_{s=1}^{k} 3^{\sum_{t=s+1}^{t=k} n_t} \cdot 2^{\sum_{t=1}^{t=s-1} \alpha_t} \cdot \left(3^{n_s} - 2^{n_s}\right)$$

$$\sum J = \sum_{s=k+1}^{j} 3^{\sum_{t=s+1}^{t=j} n_t} \cdot 2^{\sum_{t=k+1}^{t=s-1} \alpha_t} \cdot \left(3^{n_s} - 2^{n_s}\right)$$

$$\sum I = \sum_{s=j+1}^{i} 3^{\sum_{t=s+1}^{t=i} n_t} \cdot 2^{\sum_{t=j+1}^{t=s-1} \alpha_t} \cdot \left(3^{n_s} - 2^{n_s}\right)$$

Therefore, using this Nomenclature, Identity (3.3) and Identity (3.6) can be expressed as:

$$(3.7) \quad d_O^1 = 2^{\alpha_1^k + \Delta + 2\cdot(m-k)} \cdot Q_{O.m}^D - \frac{2^{w_m} + 1}{3^{n_1^k + m - k}} \cdot \left[3^{(m-k)} \cdot \sum K + 2^{\alpha_1^k + \Delta} \cdot \left(4^{(m-k)} - 3^{(m-k)}\right)\right]$$

$$(3.8) \quad a_O^1 = 2^{\alpha_1^k + \alpha_k^j + \alpha_j^i} \cdot Q_{O.i}^A - \frac{2^{w_i} + 1}{3^{n_1^k + n_k^j + n_j^i}} \cdot \left(3^{n_k^j} \cdot 3^{n_j^i} \cdot \sum K + 2^{\alpha_1^k} \cdot 3^{n_j^i} \sum J + 2^{\alpha_1^k} \cdot 2^{\alpha_k^j} \sum I\right)$$

We need to prove that, as the number of cycles of Sequence-A unfold, at some point $a_O^1$ must be larger than $d_O^1$ to provoke a contradiction – see Inequality (3.1) -. That's how we prove the conjecture. In order to do it, we have to be able to compare the identities (3.7) and (3.8) above. The existing format makes it difficult to compare. We are going to make some modifications in both expressions in such a way that the identities above are going to resemble something like:

$$a_O^1 = 2^\alpha \cdot Q_O^A - \frac{2^w + 1}{3^n} \cdot \sum A$$

$$d_O^1 = 2^\alpha \cdot Q_O^D - \frac{2^w + 1}{3^n} \cdot \sum D$$

That way, both expressions are comparable. Please note that the number $\alpha$ representing the downward steps and the fraction

$$\frac{2^w + 1}{3^n}$$

are the same for both expressions.

**Stage 1**

In this first stage we are going to match the number of downward steps for both sequences. Please be reminded that Sequence-A and Sequence-D follow parallel paths – same pattern of upward and downward steps – until Cycle-k. Once they reach this cycle, they unfold independently following different paths. Therefore, each sequence can unfold faster or slower as needed.

Initially, both sequences unfold at the same rate until they both reach Cycle-k. As we have already seen, Cycle-k occurs when the number of downward steps is larger than $k_{00}$. The final element of Sequence-A after "k" cycles is the odd number $a_F^k$. On the other hand, Sequence-D follows a downward trajectory ( "$\Delta + \alpha_k - n_k$" additional downward steps ) until it reaches 1.

Next, Sequence-A continues unfolding until its number of downward steps matches the number of downward steps of Sequence-D. This is the end of the Second Block of Sequence-A. The total number of downward steps for Sequence-A after "k" cycles is:

$$\alpha_1^k + \alpha_k^j$$

And the total number of downward steps for Sequence-D is:

$$\alpha_1^k + \Delta$$

Therefore, if we match both figures

$$\alpha_1^k + \Delta = \alpha_1^k + \alpha_k^j$$

In other words:

$$\alpha_k^j = \Delta$$

Please note that, if the number of downward steps does not exactly match, since Sequence-D enters the trivial cycle 4-2-1, we can always extend Sequence-D as many cycles as necessary, even make the last cycle one downward step short, to make the number of downward steps match exactly.

At this point, both sequences, A and D, keep unfolding but they do it in such a way that both sequences have always the same number of downward steps. This is the Third Block for Sequence-A and the Second Block for Sequence-D.

Sequence-D, after converging to 1 in Cycle-k, is a repetition of the trivial cycle 4-2-1. The number of downward steps for Sequence-D is (once the sequence enters the trivial cycle repetition):

$$2 \cdot (m-k)$$

If we match the number of downward steps for both sequences after Cycle-j:

$$(3.9) \qquad \alpha_1^k + \Delta + 2 \cdot (m-k) = \alpha_1^k + \alpha_k^j + \alpha_j^i = \alpha$$

As we just mentioned, if the number of downward steps does not exactly match, we can always stop Sequence-D one downward step short to make it match exactly. In this case, Identity (3.9) would become

$$\alpha_1^k + \Delta + 2 \cdot (m-k) - 1 = \alpha_1^k + \alpha_k^j + \alpha_j^i = \alpha$$

To make the math less cumbersome, we assume that it matches. Therefore, we carry on using Identity (3.9) but the proof would not change if it was the other case. As a result, using Identity (3.9), the number m of cycles of Sequence-D that need to be unfolded to have the exact number of downward steps in both sequences is:

$$(3.10) \qquad m - k = \frac{\alpha_j^i}{2}$$

Therefore, after unfolding Sequence-A "i" cycles and Sequence-D the exact $\frac{\alpha_j^i}{2}$ cycles as shown in Identity (3.10), both sequences have unfolded for the same number of downward steps - $\alpha$ - and identities (3.7) and (3.8) can be expressed as:

$$(3.11) \quad a_O^1 = 2^\alpha \cdot Q_{O.i}^A - \frac{2^{w_i}+1}{3^{n_1^k+n_k^j+n_j^i}} \cdot \left(3^{n_k^j} \cdot 3^{n_j^i} \cdot \sum K + 2^{\alpha_1^k} \cdot 3^{n_j^i} \sum J + 2^{\alpha_1^k} \cdot 2^\Delta \sum I\right)$$

$$(3.12) \quad d_O^1 = 2^\alpha \cdot Q_{O.m}^D - \frac{2^{w_m}+1}{3^{n_1^k+\frac{\alpha_j^i}{2}}} \cdot \left[3^{\frac{\alpha_j^i}{2}} \cdot \sum K + 2^{\alpha_1^k+\Delta} \cdot \left(4^{\frac{\alpha_j^i}{2}} - 3^{\frac{\alpha_j^i}{2}}\right)\right]$$

Please, be reminded that, even though both sequences have the same number of downward steps, the number of upward steps is different. In the final block, the number of upward steps in Sequence-A is $n_j^i$ and the number of upward steps in Sequence-D is $\frac{\alpha_j^i}{2}$.

**Stage 2**

In this second stage we are going to match the power of 3 in the denominator of the fraction

$$\frac{2^{w_j}+1}{3^n}$$

If we could multiply the denominator and numerator of the fraction of Identity (3.12) by the factor

$$3^{n_k^j+n_j^i-\frac{\alpha_j^i}{2}}$$

Identity (3.12) would not change but then Identity (3.11) and Identity (3.12) would have the same fraction

$$\frac{2^{w_j}+1}{3^n}=\frac{2^{w_j}+1}{3^{n_1^k+n_j^i+n_j^i}}$$

Before we can do that, we first have to prove that

$$3^{n_k^j+n_j^i-\frac{\alpha_j^i}{2}} \in \mathbb{N}$$

is an integer. In other words, after Cycle-k, the number of upward steps of Sequence-D is smaller than the number of upward steps in Sequence-A

$$\frac{\alpha_j^i}{2}<n_k^j+n_j^i$$

If this is the case, then we could multiply the denominator and numerator of the second addend of Identity (3.12) by it and the fraction would remain an integer.

Therefore, to prove that

$$\frac{\alpha_j^i}{2}<n_k^j+n_j^i$$

we first apply Lemma A9.2 of Appendix 9. This lemma states that if Sequence-A does not converge, in other words

$$a_F^i \geq a_O^1$$

then

$$(A9.4) \qquad 2^{\sum_1^i \alpha_t} < \left(1+\frac{1}{a_O^1}\right)^{\sum_1^i n_t} \cdot 3^{\sum_1^i n_t}$$

Since Sequence-A does not converge, it does not converge after Cycle-k either. As a result, we can apply Lemma 9.2 for Sequence-A after this cycle. To simplify the expressions, we also introduce the Nomenclature stated above for Inequality (A9.4). Since in previous Stage 1 we have matched the number of downward steps for both sequences A and B, then the number of downward steps for Block 2 and Block 3 are:

$$\alpha_k^j = \Delta$$

$$\alpha_j^i = 2 \cdot (m-k)$$

and obtain

$$(3.13) \quad 2^{\Delta+\alpha_j^i} < \left(1+\frac{1}{a_O^{k+1}}\right)^{n_k^j+n_j^i} \cdot 3^{n_k^j+n_j^i}$$

Please, be reminded that $a_O^{k+1}$ is the first element of Cycle "k+1" of Sequence-A.

We solve for $\Delta + \alpha_j^i$

$$(3.13) \quad \Delta + \alpha_j^i < \frac{\ln\left[\left(1+\frac{1}{a_O^{k+1}}\right)\cdot 3\right]}{\ln(2)} \cdot \left(n_k^j + n_j^i\right)$$

We solve for $\alpha_j^i$

$$\alpha_j^i < \frac{\ln\left[\left(1+\frac{1}{a_O^{k+1}}\right)\cdot 3\right]}{\ln(2)} \cdot \left(n_k^j + n_j^i\right) - \Delta$$

We divide all the elements above by 2

$$(3.14) \qquad \frac{\alpha_j^i}{2} < \frac{1}{2} \cdot \frac{\ln\left[\left(1+\frac{1}{a_O^{k+1}}\right)\cdot 3\right]}{\ln(2)} \cdot \left(n_k^j + n_j^i\right) - \frac{\Delta}{2}$$

Since Sequence-A does not converge,

$$a_O^{k+1} > a_O^1$$

The largest numerator in the fraction above is for $a_O^{k+1} = 3$ . For this case, the fraction above is 1.

$$\frac{1}{2} \cdot \frac{\ln\left[\left(1+\frac{1}{3}\right) \cdot 3\right]}{\ln(2)} = \frac{1}{2} \cdot \frac{\ln(4)}{\ln(2)} = 1$$

Therefore, for any $a_O^{k+1}$ larger than 3

$$\frac{1}{2} \cdot \frac{\ln\left[\left(1+\frac{1}{a_O^{k+1}}\right) \cdot 3\right]}{\ln(2)} < 1$$

$$a_O^{k+1} > 3$$

Therefore, Inequality (3.14) can be expressed as:

$$\frac{\alpha_j^i}{2} < \left(n_k^j + n_j^i\right) - \frac{\Delta}{2}$$

As a result,

$$\frac{\alpha_j^i}{2} < n_k^j + n_j^i$$

Which is what we wanted to prove. Please, be reminded that $\frac{\alpha_j^i}{2}$ represents, after Cycle-k, the number of upward steps in Sequence-D and $n_k^j + n_j^i$ the number of upward steps in Sequence-A.

Therefore, $3^{n_k^j+n_j^i-\frac{\alpha_j^i}{2}}$ is an integer number and it is possible to multiply and divide the numerator and denominator of Identity (3.12) by

$$3^{n_k^j+n_j^i-\frac{\alpha_j^i}{2}} \in \mathrm{N}$$

and obtain an integer. Therefore, Identity (3.12) becomes:

$$(3.15)\quad d_O^1 = 2^\alpha \cdot Q_{O.m}^D - \frac{2^{w_m}+1}{3^{n_1^k+n_k^j+n_j^i}} \cdot \left[3^{n_k^j+n_j^i} \cdot \sum K + 2^{\alpha_1^k+\Delta} \cdot 3^{n_k^j+n_j^i-\frac{\alpha_j^i}{2}} \cdot \left(4^{\frac{\alpha_j^i}{2}} - 3^{\frac{\alpha_j^i}{2}}\right)\right]$$

Please note that now identities (3.11) and (3.15) now share the same power of 3 in the denominator of the fraction above.

$$n_1^k + n_k^j + n_j^i = n$$

$$\frac{2^w + 1}{3^{n_1^k + n_k^j + n_j^i}} = \frac{2^w + 1}{3^n}$$

Therefore, applying the above to identities (3.11) and (3.15) we obtain

$$(3.16)\quad a_O^1 = 2^\alpha \cdot Q_{O.i}^A - \frac{2^{w_i} + 1}{3^n} \cdot \left(3^{n_k^j} \cdot 3^{n_j^i} \cdot \sum K + 2^{\alpha_1^k} \cdot 3^{n_j^i} \sum J + 2^{\alpha_1^k} \cdot 2^\Delta \sum I\right)$$

$$(3.17)\quad d_O^1 = 2^\alpha \cdot Q_{O.m}^D - \frac{2^{w_m} + 1}{3^n} \cdot \left[3^{n_k^j + n_j^i} \cdot \sum K + 2^{\alpha_1^k + \Delta} \cdot 3^{n_k^j + n_j^i - \frac{\alpha_j^i}{2}} \cdot \left(4^{\frac{\alpha_j^i}{2}} - 3^{\frac{\alpha_j^i}{2}}\right)\right]$$

**Stage 3**

In this third stage we are going to match the power of 2 ($2^w$) in the numerator of the fraction. As it is now, Identity (3.16) for Sequence-A shows

$$\frac{2^{w_i} + 1}{3^n}$$

But Identity (3.17) for Sequence-D shows

$$\frac{2^{w_m} + 1}{3^n}$$

And

$$w_m \neq w_i$$

We use Lemma A9.1 from Appendix 9 to solve the problem. Lemma A9.1 from Appendix 9 states that we can always substitute a larger $w_B$ for $w_A$,

$$w_B > w_A$$

This substitution requires the following changes:

$$(A9.3)\quad Q_{O.i}^B = Q_{O.i}^A + 2^{w_A - \sum_{t=1}^{t=i} \alpha_t} \cdot \frac{2^{w_B - w_A} - 1}{3^{\sum_{t=1}^{t=i} n_t}} \cdot \sum_{s=1}^{i} 3^{\sum_{t=s+1}^{t=i} n_t} \cdot 2^{\sum_{t=1}^{t=s-1} \alpha_t} \cdot (3^{n_s} - 2^{n_s})$$

and the original odd number $a_O^1$ remains the same.

Since, as we just previously saw, the number of upward steps is larger for Sequence-A than for Sequence-D

$$\frac{\alpha_j^i}{2} < n_k^j + n_j^i$$

Then, the parameters $w_i$ and $w_m$ from Identities (3.16) and (3.15) are such (See Nomenclature)

$$w_i > w_m$$

Therefore, we apply Lemma A9.1 to the odd number $d_O^1$ of Sequence-D. As a result, the parameter $Q_{O.m}^{D}$ of Identity (3.15) becomes a new $Q_{O.m}^{DD}$:

$$(A9.3)\quad Q_{O.m}^{DD} = Q_{O.m}^{D} + 2^{w_m - [\sum_{t=1}^{t=k} \alpha_t + \Delta + \alpha_j^i]} \,.\, \frac{2^{w_i - w_m} - 1}{3^n} \cdot \left[ 3^{n_k^j + n_j^i} \cdot \sum K + 2^{\alpha_1^k + \Delta} \,.\, 3^{n_k^j + n_j^i - \frac{\alpha_j^i}{2}} \cdot \left( 4^{\frac{\alpha_j^i}{2}} - 3^{\frac{\alpha_j^i}{2}} \right) \right]$$

In this case, applying Lemma A9.1, we can express the odd number $d_O^1$ as:

$$(3.18)\qquad d_O^1 = 2^{\alpha} \cdot Q_{O.m}^{DD} - \frac{2^{w_i} + 1}{3^n} \cdot \left[ 3^{n_k^j + n_j^i} \cdot \sum K + 2^{\alpha_1^k + \Delta} \,.\, 3^{n_k^j + n_j^i - \frac{\alpha_j^i}{2}} \cdot \left( 4^{\frac{\alpha_j^i}{2}} - 3^{\frac{\alpha_j^i}{2}} \right) \right]$$

and identities (3.18) an (3.17) represent the same $d_O^1$.

Please note that Identities (3.18) for $d_O^1$ and (3.16) for $a_O^1$ they both now share the parameters $\frac{2^{w_i+1}}{3^n}$ and $\alpha$.

$$(3.16)\quad a_O^1 = 2^{\alpha} \cdot Q_{O.i}^{A} - \frac{2^{w_i} + 1}{3^n} \cdot \left( 3^{n_k^j} \cdot 3^{n_j^i} \cdot \sum K + 2^{\alpha_1^k} \cdot 3^{n_j^i} \sum J + 2^{\alpha_1^k} \cdot 2^{\Delta} \sum I \right)$$

$$(3.18)\quad d_O^1 = 2^{\alpha} \cdot Q_{O.m}^{DD} - \frac{2^{w_i} + 1}{3^n} \cdot \left[ 3^{n_k^j + n_j^i} \cdot \sum K + 2^{\alpha_1^k + \Delta} \,.\, 3^{n_k^j + n_j^i - \frac{\alpha_j^i}{2}} \cdot \left( 4^{\frac{\alpha_j^i}{2}} - 3^{\frac{\alpha_j^i}{2}} \right) \right]$$

## 4. Evolution of subtrahends.

We are trying to compare identities that represent the initial numbers of Sequence-A and Sequence-D. So far, we have matched the number of downward steps and the fraction $\frac{2^{w_i+1}}{3^n}$ for both identities. In this section we are going to assess how the summation of the identity that represents an odd number evolves as the number of cycles grows indefinitely. That way, in a later section, we will be able to assess if the initial number of Sequence-A can also remain indefinitely smaller than the initial number of Sequence-D.

### Lemma 4.1

Let $a_O^1$ be the initial element of a Sequence-A that does not converge

$$(4.1) \quad a_O^1 = 2^{\alpha} \cdot Q_{O.i}^{A} - \frac{2^{w_i}+1}{3^{n_1^i}} \cdot \left(\sum I\right)$$

$$a_O^1 \leq a_F^i$$

$$a_O^1 < 2^{K_{00}}$$

Let $d_O^1$ be the initial element of a Sequence-D that converges to 1 in Cycle-k.

$$(4.2) \quad d_O^1 = 2^{\alpha} \cdot Q_{O.m}^{DD} - \frac{2^{w_i}+1}{3^{n_1^k+m-k}} \cdot \left[3^{m-k} \cdot \sum K + 2^{\alpha_1^k} \cdot \left(4^{m-k} - 3^{m-k}\right)\right]$$

Both sequences have the same number of downward steps

$$(4.3) \quad \alpha = \alpha_1^i = \alpha_1^k + 2 \cdot (m-k)$$

As new cycles unfold indefinitely, the summations within the brackets in Identity (4.2) - Sequence-D –

$$\left[3^{m-k} \cdot \sum K + 2^{\alpha_1^k} \cdot \left(4^{m-k} - 3^{m-k}\right)\right]$$

grow faster than the summation within the brackets in Identity (4.1) – Sequence-A -

$$\left[\sum I\right]$$

In other words

$$\lim_{i \to \infty} \left[3^{m-k} \cdot \sum K + 2^{\alpha_1^k} \cdot \left(4^{m-k} - 3^{m-k}\right)\right] > \lim_{i \to \infty} \left[\sum I\right]$$

## Proof

As we already mention, so far, we have matched the number of downward steps and the fraction $\frac{2^{w_i+1}}{3^n}$ for both identities. That allows us to compare them later. Therefore, we first match the power of 3 in the denominator of the fraction in both Identity (4.1) and (4.2).

We multiply the numerator and denominator of the fraction in Identity (4.2) by

$$3^{n_1^i-(n_1^k+m-k)}$$

We obtain:

$$(4.4)\quad d_O^1=2^\alpha\cdot Q_{O.m}^{DD}-\frac{2^{w_i}+1}{3^{n_1^i}}\cdot\left[3^{n_1^i-n_1^k}\cdot\sum K+2^{\alpha_1^k}\cdot 3^{n_1^i-(n_1^k+m-k)}\cdot\left(4^{m-k}-3^{m-k}\right)\right]$$

Since both sequences have the same number of downward steps – see Identity (4.3) -, then

$$(4.5)\quad m-k=\frac{\alpha_1^i-\alpha_1^k}{2}$$

We apply the above to Identity (4.4)

$$d_O^1=2^\alpha\cdot Q_{O.m}^{DD}-\frac{2^{w_i}+1}{3^{n_1^i}}\cdot\left[3^{n_1^i-n_1^k}\cdot\sum K+2^{\alpha_1^k}\cdot 3^{n_1^i-\left(n_1^k+\frac{\alpha_1^i-\alpha_1^k}{2}\right)}\cdot\left(4^{\frac{\alpha_1^i-\alpha_1^k}{2}}-3^{\frac{\alpha_1^i-\alpha_1^k}{2}}\right)\right]$$

We introduce the factor

$$3^{\frac{\alpha_1^i-\alpha_1^k}{2}}$$

within the brackets to simplify the above and obtain

$$(4.6)\quad d_O^1=2^\alpha\cdot Q_{O.m}^{DD}-\frac{2^{w_i}+1}{3^{n_1^i}}\cdot\left[3^{n_1^i-n_1^k}\cdot\sum K+2^{\alpha_1^k}\cdot 3^{n_1^i-n_1^k}\cdot\left(\left(\frac{4}{3}\right)^{\frac{\alpha_1^i-\alpha_1^k}{2}}-1\right)\right]$$

Therefore, as the number of cycles "i" grow indefinitely, if we want to know which subtrahend grows faster, then we need to compare the subtrahend of both identities. We need to compare the subtrahend of Identity (4.1) – Sequence-A –

$$\sum I$$

to the subtrahend of Identity (4.6) – Sequence-D -

$$3^{n_1^i-n_1^k}\cdot\sum K+2^{\alpha_1^k}\cdot 3^{n_1^i-n_1^k}\cdot\left(\left(\frac{4}{3}\right)^{\frac{\alpha_1^i-\alpha_1^k}{2}}-1\right)$$

We compare them both by calculating the following limit and make the number cycles grow indefinitely.

$$(4.7)\quad L=\lim_{i\to\infty}\frac{3^{n_1^i-n_1^k}\cdot\sum K+2^{\alpha_1^k}\cdot 3^{n_1^i-n_1^k}\cdot\left(\left(\frac{4}{3}\right)^{\frac{\alpha_1^i-\alpha_1^k}{2}}-1\right)}{\sum I}\overset{???}{<>}1$$

A larger or smaller than 1 Limit L would tell us which elements are growing faster, and which ones are slower.

Please be reminded that both sequences unfold at different rates, but they do it in such a way that they both have the same number of downward steps. As Sequence-A unfolds new cycles, the number of downward steps - $\alpha_1^i$ – must remain within two limits. As we saw in Appendix 6, it must be at least one step more than the number of upward steps - $n_1^i$-. In other words, the minimum value for the downward steps is:

$$(4.8)\quad n_1^i+1\le\left(\alpha_1^i\right)_{min}$$

Also, since Sequence-A does not converge, as we have previously seen in Stage 2 of Section 3, Lemma A9.2 applies. Therefore $\alpha_1^i$ also has an upper limit. We use Identity (3.13).

$$(3.13)\quad\left(\alpha_1^i\right)_{MAX}<\frac{\ln\left[\left(1+\frac{1}{a_O^1}\right)\cdot 3\right]}{\ln(2)}\cdot n_1^i$$

We combine inequalities (4.8) and (3.13)

$$n_1^i+1\le\alpha_1^i<\frac{\ln\left[\left(1+\frac{1}{a_O^1}\right)\cdot 3\right]}{\ln(2)}\cdot n_1^i$$

Therefore, as long as Sequence-A does not converge, the number of downward steps of Sequence-A must remain within the limits shown above.

We want to prove that the subtrahend of Sequence-D grows faster than the subtrahend of Sequence-A. If this is the case, then Limit L must be larger than 1. If we prove that, even for the minimum value of Limit L, Limit L is still larger than 1, then we would have proved Lemma 4.1.

The minimum value for Limit (4.7) is obtained for the minimum value in the numerator and the maximum value in the denominator.

The minimum value in the numerator of Limit (4.7) is for:

$$n_1^i + 1 = \left(\alpha_1^i\right)_{min}$$

and the maximum value in the denominator of Limit (4.7)

$$\left[\sum I\right]_{MAX}$$

is obtained using Inequality (3.13), for the maximum $\alpha_1^i = \left(\alpha_1^i\right)_{MAX}$

$$\left(\alpha_1^i\right)_{MAX} = \frac{\ln\left[\left(1 + \frac{1}{a_O^1}\right) \cdot 3\right]}{\ln\left(2\right)} \cdot n_1^i$$

In this case, Lemma A9.3 from Appendix 9 gives the expression for the maximum value of the summation for Sequence-A.

$$(A9.9) \quad \left[\sum I\right]_{MAX} < 3^{\sum_{t=1}^{t=i} n_t} \cdot \left(a_O^1 + 1\right) \cdot \left[1 + \frac{1}{a_O^1}\right]^{\sum_1^i n_t}$$

To simplify the above, we apply the same Nomenclature as before.

$$(4.9) \quad \left[\sum I\right]_{MAX} = 3^{n_1^i} \cdot \left(a_O^1 + 1\right) \cdot \left[1 + \frac{1}{a_O^1}\right]^{n_1^i}$$

As we already stated above, since we want to know if the Limit (4.7) is smaller or larger than 1, if we apply the maximum value of the denominator - Inequality (4.9) – and the minimum value of the numerator – Inequality (4.8) - to Limit (4.7), we obtain a smaller value for limit L:

$$L > \lim_{i \to \infty} \frac{3^{n_1^i - n_1^k} \cdot \sum K + 2^{\alpha_1^k} \cdot 3^{n_1^i - n_1^k} \cdot \left( \left(\frac{4}{3}\right)^{\frac{n_1^i + 1 - \alpha_1^k}{2}} - 1 \right)}{3^{n_1^i} \cdot \left(a_O^1 + 1\right) \cdot \left[1 + \frac{1}{a_O^1}\right]^{n_1^i}}$$

min

MAX

We simplify the above limit by eliminating the common factor $3^{n_1^i}$ from the denominator and numerator. Also, since

$$\lim_{i\to\infty}\left(\frac{4}{3}\right)^{\frac{n_1^i+1-\alpha_1^k}{2}} \gg 1$$

we can eliminate the "1" within the brackets.

We obtain

$$L > \lim_{i\to\infty}\frac{\sum K + 2^{\alpha_1^k}\cdot\left(\frac{4}{3}\right)^{\frac{n_1^i+1-\alpha_1^k}{2}}}{3^{n_1^k}\cdot(a_O^1+1)\cdot\left[1+\frac{1}{a_O^1}\right]^{n_1^i}}$$

We split the limit in two

$$L > \lim_{i\to\infty}\frac{\sum K}{3^{n_1^k}\cdot(a_O^1+1)\cdot\left[1+\frac{1}{a_O^1}\right]^{n_1^i}} + \frac{2^{\alpha_1^k}\cdot\left(\frac{4}{3}\right)^{\frac{n_1^i+1-\alpha_1^k}{2}}}{3^{n_1^k}\cdot(a_O^1+1)\cdot\left[1+\frac{1}{a_O^1}\right]^{n_1^i}}$$

Since the factor

$$\left[1+\frac{1}{a_O^1}\right]^{n_1^i}$$

in the denominator of the first fraction keeps growing while the numerator remains constant, as the number of cycles grows indefinitely, this first fraction becomes irrelevant. Therefore, the limit becomes:

$$L > \lim_{i\to\infty}\frac{2^{\alpha_1^k}\cdot\left(\frac{4}{3}\right)^{\frac{n_1^i+1-\alpha_1^k}{2}}}{3^{n_1^k}\cdot(a_O^1+1)\cdot\left[1+\frac{1}{a_O^1}\right]^{n_1^i}}$$

We extract all the constant elements out of the limit

$$L > \frac{2^{\alpha_1^k}\cdot\left(\frac{4}{3}\right)^{\frac{1-\alpha_1^k}{2}}}{3^{n_1^k}\cdot(a_O^1+1)}\cdot\lim_{i\to\infty}\frac{\left(\frac{4}{3}\right)^{\frac{n_1^i}{2}}}{\left[1+\frac{1}{a_O^1}\right]^{n_1^i}}$$

The maximum value in the denominator for a number that does not converge to 1 in a single cycle is for is obtained for $a_O^1 = 7$

$$\left(1+\frac{1}{a_O^1}\right)_{MAX}=1+\frac{1}{7}=\frac{8}{7}$$

If we substitute the above in the limit, the inequality holds. We also simplify the constant fraction outside the limit.

$$L>\frac{2\cdot(3)^{\frac{\alpha_1^k-1}{2}-n_1^k}}{(a_O^1+1)}\cdot\lim_{i\to\infty}\frac{\left(\frac{4}{3}\right)^{\frac{n_1^i}{2}}}{\left[\frac{8}{7}\right]^{n_1^i}}$$

We simplify the above

$$L>\frac{2\cdot(3)^{\frac{\alpha_1^k-1}{2}-n_1^k}}{(a_O^1+1)}\cdot\lim_{i\to\infty}1.010^{n_1^i}\gg 1$$

Therefore, even for the minimum value of Limit L, the limit is still larger than one. As a result, as new cycles unfold indefinitely, the summations within the brackets in Identity (4.6) - Sequence-D

$$\left[3^{n_1^i-n_1^k}\cdot\sum K+2^{\alpha_1^k}\cdot 3^{n_1^i-n_1^k}\cdot\left(\left(\frac{4}{3}\right)^{\frac{n_1^i+1-\alpha_1^k}{2}}-1\right)\right]$$

grow faster than the summations within the brackets in Identity (4.1) – Sequence-A -

$$\left[\sum I\right]$$

and eventually become larger.

$$\lim_{i\to\infty}\left[3^{m-k}\cdot\sum K+2^{\alpha_1^k}\cdot\left(4^{m-k}-3^{m-k}\right)\right]>\lim_{i\to\infty}\left[\sum I\right]$$

Which is what we wanted to prove.

In other words, even though initially the subtrahend in Identity (4.2) may be smaller than the subtrahend in Identity (4.1), as new cycles unfold, eventually, the opposite is true. In other words, as Sequence-A unfolds new cycles and Sequence-D indefinitely unfolds the trivial cycle, at some point, the subtrahend in Identity (4.2) becomes larger than the one in Identity (4.1).

# 5. $Q^{A}_{O.i}$ and $Q^{DD}_{O.m}$ relationship

We want to prove that Sequence-A converges to 1. In order to do it, we have generated another Sequence-D that does converge to 1. In Section 3 we have adapted the identities defining the initial odd numbers $a^{1}_{O}$ and $d^{1}_{O}$ of each sequence to able to compare them. - See Identity (3.16) and Identity (3.18) –. In the previous section, we saw how the subtrahends in Identity (3.16) and Identity (3.18) evolve as the number of cycles grows indefinitely. The last variables pending comparison are $Q^{A}_{O.i}$ and $Q^{DD}_{O.m}$. In this section we are going to see the constraints between $Q^{A}_{O.i}$ and $Q^{DD}_{O.m}$.

## Lemma 5.1

Let $a^{1}_{O}$ be the initial element of a Sequence-A

$$(5.1) \qquad a^{1}_{O} = 2^{\alpha} \cdot Q^{A}_{O.i} - \frac{2^{w_i}+1}{3^{n}} \cdot \sum A$$

Let $d^{1}_{O}$ be the initial element of a Sequence-D that converges to 1.

$$(5.2) \qquad d^{1}_{O} = 2^{\alpha} \cdot Q^{DD}_{O.m} - \frac{2^{w_i}+1}{3^{n}} \cdot \sum D$$

$$d^{i}_{F} = 1$$

Please note that Identity (5.1) and Identity (5.2) share the same parameters $\alpha$ and $\frac{2^{w_i+1}}{3^{n}}$.

If $d^{1}_{O}$ is larger than $a^{1}_{O}$

$$d^{1}_{O} > a^{1}_{O}$$

Then

$$Q^{DD}_{O.m} < Q^{A}_{O.i}$$

## Proof

We apply Lemma A6.3 from Appendix 6.

$a^{1}_{O}$ is the initial element of a Sequence-A made of "$i$" cycles and $a^{i}_{F}$ the final element of this sequence after “i” cycles. According to Lemma A6.3, the following identity is true:

$$(A6.5) \qquad Q_{O.i} = \frac{2^{3^{\sum_{t=1}^{t=i} n_t - 1} \cdot j^{\delta}_{o.i}} + 1}{3^{\sum_{t=1}^{t=i} n_t}} \cdot a^{i}_{F} - 2^{3^{\sum_{t=1}^{t=i} n_t - 1} \cdot j^{\delta}_{o.i} - \sum_{t=1}^{t=i} \alpha_t} \cdot a^{1}_{O}$$

We apply the above to Sequence-A using the Nomenclature of Section 3 and obtain:

$$(5.3) \qquad Q_{O.i}^{A} = \frac{2^{w_i}+1}{3^n} \cdot \ a_F^i - 2^{w_i - \alpha} \cdot a_O^1$$

We do the same with Sequence-D. Please be aware that Sequence-D converges to 1 in the last cycle.

$$d_F^i = 1$$

Therefore:

$$(5.4) \qquad Q_{O.m}^{DD} = \frac{2^{w_i}+1}{3^n} - 2^{w_i - \alpha} \cdot d_O^1$$

On the other hand, if we compare Identity (5.3) and Identity (5.4), since the minuend in Identity (5.3) is larger than the one in Identity (5.4)

$$a_F^i > 1$$

and the subtrahend in Identity (5.3) is smaller than the one in Identity (5.4)

$$a_O^1 < d_O^1$$

Then, we infer that

$$Q_{O.m}^{DD} < Q_{O.i}^{A}$$

Which is what we wanted to prove.

# 6. Final Proof

As we have already seen, Sequence-A is a sequence that does not converge

$$a_O^1 \leq a_F^i$$

Also, its initial element is a finite number

$$a_O^1 < 2^{K_{00}}$$

We build a Sequence-D that follows a parallel path to Sequence-A during the first "k" cycles but converges to 1 in this very last Cycle-k. As we have seen before, since both sequences follow parallel paths during the first k cycles, then $d_O^1$ must be larger than $a_O^1$.

$$d_O^1 > a_O^1$$

Sequence-A unfolds in three blocks. The First Block represents the first "k" cycles until the number of downward steps is larger than $K_{00}$. The Second Block represents the following "j-k" cycles until both sequences, A and D, have the same number of downward steps. The final Third Block represents the following "i-j" cycles. Sequence-A continues unfolding indefinitely.

On the other hand, Sequence-D unfolds until it reaches Cycle-k, but once it reaches it, it converges to 1 in one single trajectory. Appendix 7 proves that it is always possible to find such Sequence-D. Logically, Sequence-D has a final downward trajectory of additional downward steps until it reaches 1.

Once Sequence-A is in the Third Block, both sequences unfold in such a way that the number of downward steps is always the same. Sequence-A unfolds its own cycles while Sequence-D unfolds the trivial cycle 4-2-1.

The initial element of Sequence-A can be expressed using Identity (3.16)

$$(3.16) \quad a_O^1 = 2^{\alpha_1^k+\Delta+2\cdot(m-k)} \cdot Q_{O.i}^A - \frac{2^{w_i}+1}{3^n} \cdot \left(3^{n_k^j} \cdot 3^{n_j^i} \cdot \sum K + 2^{\alpha_1^k} \cdot 3^{n_j^i} \sum J + 2^{\alpha_1^k} \cdot 2^{\Delta} \sum I\right)$$

And the initial element of Sequence-D can be expressed using Identity (3.18)

$$(3.18) \quad d_O^1 = 2^{\alpha_1^k+\Delta+\alpha_j^i} \cdot Q_{O.m}^{DD} - \frac{2^{w_i}+1}{3^n} \cdot \left[3^{n_k^j+n_j^i} \cdot \sum K + 2^{\alpha_1^k+\Delta} \cdot 3^{n_k^j+n_j^i-\frac{\alpha_j^i}{2}} \cdot \left(4^{\frac{\alpha_j^i}{2}} - 3^{\frac{\alpha_j^i}{2}}\right)\right]$$

Also, as we saw in Lemma 5.1 of Section 5,

$$Q_{O.m}^{DD} < Q_{O.i}^A$$

Please note that Sequence-A and Sequence-D unfold at different rates. Even though they both have the same number of downward steps, they have a different number of upward steps. Sequence-A unfolds "i" cycles with "$n_1^k + n_k^j + n_j^i$" upward steps while Sequence-D unfolds "m" cycles with "$n_1^k + \frac{\alpha_j^i}{2}$" upward steps.

Initially, the subtrahend of Identity (3.16) is larger than the subtrahend of Identity (3.18), nevertheless, as new cycles unfold, according to Lemma 4.1 - as new cycles unfold indefinitely - the summations within the brackets in Identity (3.18) - Sequence-D –

$$\left[3^{n_k^j+n_j^i}\cdot\sum K + 2^{\alpha_1^k+\Delta}\cdot 3^{n_k^j+n_j^i-\frac{\alpha_j^i}{2}}\cdot\left(4^{\frac{\alpha_j^i}{2}} - 3^{\frac{\alpha_j^i}{2}}\right)\right]$$

grow faster than the summations within the brackets in Identity (3.16) – Sequence-A -

$$\left[3^{n_k^j}\cdot 3^{n_j^i}\cdot\sum K + 2^{\alpha_1^k}\cdot 3^{n_j^i}\sum J + 2^{\alpha_1^k}\cdot 2^{\Delta}\sum I\right]$$

Therefore, it must be a Cycle-m for Sequence-D and a Cycle-i for Sequence-A for which the summations within the brackets in Sequence-D are larger than Sequence-A.

$$\left[3^{n_k^j+n_j^i}\cdot\sum K + 2^{\alpha_1^k+\Delta}\cdot 3^{n_k^j+n_j^i-\frac{\alpha_j^i}{2}}\cdot\left(4^{\frac{\alpha_j^i}{2}} - 3^{\frac{\alpha_j^i}{2}}\right)\right] >$$

$$> \left[3^{n_k^j}\cdot 3^{n_j^i}\cdot\sum K + 2^{\alpha_1^k}\cdot 3^{n_j^i}\sum J + 2^{\alpha_1^k}\cdot 2^{\Delta}\sum I\right]$$

As a result, since, as we just saw in Lemma 5.1

$$Q_{O.m}^{DD} < Q_{O.i}^{A}$$

then, if both sequences continue unfolding indefinitely, Identity (3.18) must produce a smaller number than Identity (3.16) since the minuend in Identity (3.18) is smaller and the subtrahend is larger than the ones in Identity (3.16). As a result, we must conclude that $a_O^1$ is larger than $d_O^1$:

$$a_O^1 > d_O^1$$

Which is a contradiction since $a_O^1$ was smaller than $d_O^1$ - see Identity (3.1) -. This is the case since both sequences follow parallel paths during the first "k" cycles and, according to Lemma A8.2 from Appendix 8, only one number can be smaller than $2^{k_{00}}$ which is $a_O^1$.

As a result, Sequence-A cannot unfold without converging to 1 indefinitely. Finally, since we proved in Appendix 7 that for any combination of cycles made of upward and downward steps, there is always a sequence with the same combination of cycles that converges to 1, at some point, before the number of cycles is too large and the contradiction above occurs, one of these sequences must be Sequence-A, and therefore, Sequence-A converges to 1.

**All sequences converge to 1.**

# Appendix 1 – Odd Numbers

In this appendix, we will prove that, for all n, any odd number can be written as follows:

$$(A1.1) \quad a_O = 2^{\alpha} \cdot Q_O - \frac{2^{3^{n-1} \cdot j_o^{\beta}} + 1}{3^n} \cdot (3^n - 2^n)$$

$$j_o^{\beta}, Q_O \in \{odd\}$$

$$\alpha, n \in \mathbb{N}$$

$$3^{n-1} \cdot j_o^{\beta} > \alpha$$

$$\forall\, n \in \mathbb{N}$$

For example, let's take the odd number

$$a_O = 739$$

In this case, using identity (A1.1):

| n | 2 |
|---|---|
| $\alpha$ | 9 |
| $K_O$ | 37 |
| $j_o^{\beta}$ | 5 |

Therefore, this number can also be written as:

$$a_O = 2^9 \cdot 37 - \frac{2^{3^{(2-1)} \cdot 5} + 1}{3^2} \cdot (3^2 - 2^2) = 739$$

This is important because, as we shall see, it allows us to adapt the format of an odd number to the nature of the problem. In the long run, we will be able to create a formula that allows us to know the elements of a sequence as it unfolds from the initial odd number.

## Proof

**Lemma A1.1 (*)**

$$(A1.2)\quad \frac{2^{2\cdot 3^{n-1}\cdot j}-1}{3^n}=\frac{2^{2\cdot 3^{n-1}}-1}{3^n}\cdot\left[1+2^{2\cdot 3^{n-1}}+2^{2\cdot 3^{n-1}\cdot 2}+\cdots+2^{2\cdot 3^{n-1}\cdot(j-1)}\right]$$

(*) In Appendix 4 we prove that $(2^{2\cdot 3^{n-1}\cdot j}-1)$ is evenly divisible by $3^n$. In other words:

$$\left\{\frac{2^{2\cdot 3^{n-1}\cdot j}-1}{3^n}\right\}\in\mathbb{N}$$

**Proof**

Let $a_j$ be

$$(A1.3)\quad a_j=\frac{2^{2\cdot 3^{n-1}}-1}{3^n}\cdot\left[1+2^{2\cdot 3^{n-1}}+2^{2\cdot 3^{n-1}\cdot 2}+2^{2\cdot 3^{n-1}\cdot 3}+\cdots+2^{2\cdot 3^{n-1}\cdot(j-1)}\right]$$

Therefore, we can rewrite the above as:

$$(A1.4)\quad a_j=\frac{2^{2\cdot 3^{n-1}}-1}{3^n}\cdot S$$

and "S" is a simple geometric series

$$S=1+2^{2\cdot 3^{n-1}}+2^{2\cdot 3^{n-1}\cdot 2}+2^{2\cdot 3^{n-1}\cdot 3}+\cdots+2^{2\cdot 3^{n-1}\cdot(j-1)}$$

or

$$S=1+2^{2\cdot 3^{n-1}}+\left(2^{2\cdot 3^{n-1}}\right)^2+\left(2^{2\cdot 3^{n-1}}\right)^3+\cdots+\left(2^{2\cdot 3^{n-1}}\right)^{j-1}$$

We solve the simple geometric series

$$S=\frac{\left(2^{2\cdot 3^{n-1}}\right)^j-1}{2^{2\cdot 3^{n-1}}-1}$$

or

$$S=\frac{2^{2\cdot 3^{n-1}\cdot j}-1}{2^{2\cdot 3^{n-1}}-1}$$

We apply the above to identity (A1.4):

$$a_j = \frac{2^{2\cdot 3^{n-1}} - 1}{3^n} \cdot \frac{2^{2\cdot 3^{n-1}\cdot j} - 1}{2^{2\cdot 3^{n-1}} - 1}$$

We simplify the above by eliminating $\left(2^{2\cdot 3^{n-1}} - 1\right)$ from the numerator and the denominator

$$(A1.5) \quad a_j = \frac{2^{2\cdot 3^{n-1}\cdot j} - 1}{3^n}$$

We compare Identity (A1.5) and Identity (A1.3)

$$a_j = \frac{2^{2\cdot 3^{n-1}\cdot j} - 1}{3^n} = \frac{2^{2\cdot 3^{n-1}} - 1}{3^n} \cdot \left[1 + 2^{2\cdot 3^{n-1}} + 2^{2\cdot 3^{n-1}\cdot 2} + 2^{2\cdot 3^{n-1}\cdot 3} + \cdots + 2^{2\cdot 3^{n-1}\cdot (j-1)}\right]$$

Which is what we wanted to prove.

**Lemma A1.2**

$$\forall n \in \mathbb{N}$$

all odd numbers can be divided in $j$ subsets:

$$(A1.6) \quad B_j = \left\{2^{2\cdot 3^{n-1}(j-1)+k_n} \cdot K_O - 1\right\}$$

$$k_n \in \{1; 2; \ldots; 2\cdot 3^{n-1}\}$$

$$K_O \in \{odd\}$$

$$j \in \mathbb{N}$$

and the union of the $j$ subsets $B_j$ makes the set of all odd numbers:

$$\{odd\} = \bigcup_{j=1}^{j=\infty} B_j$$

**Proof**

Any odd number can be expressed as

$$(A1.7) \quad a_O = 2^t \cdot K_O - 1$$

$$t \in \mathbb{N}$$

$$K_O \in \{odd\}$$

Since all the Positive Integers

$$t \in \{1; 2; 3; \dots \dots.\}$$

can be divided in “j” subset as follows

| t | | | | | |
|---|---|---|---|---|---|
| j=1 | j=2 | j=3 | … | j | |
| 1 | $2 \cdot 3^{n-1} + 1$ | $(2 \cdot 3^{n-1}) \cdot 2 + 1$ | …. | $(2 \cdot 3^{n-1}) \cdot (j-1) + 1$ | …. |
| 2 | $2 \cdot 3^{n-1} + 2$ | $(2 \cdot 3^{n-1}) \cdot 2 + 2$ | …. | $(2 \cdot 3^{n-1}) \cdot (j-1) + 2$ | …. |
| 3 | $2 \cdot 3^{n-1} + 3$ | $(2 \cdot 3^{n-1}) \cdot 2 + 3$ | …. | $(2 \cdot 3^{n-1}) \cdot (j-1) + 3$ | …. |
| … | … | … | | … | |
| $k_n$ | $2 \cdot 3^{n-1} + k_n$ | $(2 \cdot 3^{n-1}) \cdot 2 + k_n$ | | $(2 \cdot 3^{n-1}) \cdot (j-1) + k_n$ | |
| … | … | … | | … | |
| $2 \cdot 3^{n-1}$ | $(2 \cdot 3^{n-1}) \cdot 2$ | $(2 \cdot 3^{n-1}) \cdot 3$ | | $(2 \cdot 3^{n-1}) \cdot (j)$ | |

We can therefore express “t” as:

$$t = 2 \cdot 3^{n-1}(j-1) + k_n$$

$$\forall\, n \in \mathbb{N}$$

$$k_n \in \{1,2, \dots ,2 \cdot 3^{n-1}\}$$

$$j \in \mathbb{N}$$

Therefore, if we apply the above to identity (A1.7) we obtain

$$b_O = 2^{2 \cdot 3^{n-1}(j-1)+k_n} \cdot K_O - 1$$

then we can state that all odd numbers can be divided in $j$ subsets:

$$(A1.8) \quad B_j = \left\{2^{2 \cdot 3^{n-1}(j-1)+k_n} \cdot K_O - 1\right\}$$

And the union of the $j$ subsets $B_j$ makes the set of all odd numbers:

$$\{odd\} = \bigcup_{j=1}^{j=\infty} B_j$$

**Lemma A1.3**

If we subtract the same EVEN number to the set of all ODD numbers, the new set also contains all the positive ODD numbers.

**Proof**

Let A be the set of odd numbers

$$A = \{1,3,5,7, \dots, K_O, \dots\}$$

$$K_O \in \{odd\}$$

We subtract $K_E$

$$K_E \in \{even\}$$

to all numbers in the set A above and obtain:

$$A_1 = \{1 - K_E; 3 - K_E; \dots : -1; 1; 3; 5; 7; \dots; K_O; \dots.\}$$

The intersection of the two sets gives the set of positive odd numbers

$$A \cap A_1 = A = \{odd\}$$

which is what we wanted to prove

For example, let A be the set of all odd numbers

$$A = \{1; 3; 5; 7; \dots; K_O; \dots.\}$$

If we subtract the same EVEN number to all numbers,

$$K_E = 4$$

we obtain

$$A_1 = \{-3; -1; 1; 3; 5; 7; \dots; K_O,; \dots.\}$$

$$K_O \in \{odd\}$$

We extract the positive odd numbers

$$A \cap A_1 = A = \{1; 3; 5; 7; \dots; K_O; \dots.\}$$

**Lemma A1.4**

Let $A_0$ be a set of odd numbers

$$A_0 = \left\{2^{2\cdot 3^{n-1}(j_a-1)+k_n^a} \cdot K_O^a - 1 - 2^n \cdot \frac{2^{2\cdot 3^{n-1}\cdot j_a} - 1}{3^n}\right\}$$

and let $B_0$ be another set of odd numbers

$$B_0 = \left\{2^{2\cdot 3^{n-1}(j_b-1)+k_n^b} \cdot K_O^b - 1 - 2^n \cdot \frac{2^{2\cdot 3^{n-1}\cdot j_b} - 1}{3^n}\right\}$$

and

$$j_a \neq j_b$$

<u>No odd number $a_0$ can be a member of two subsets at the same time.</u> In other words

if

$$a_0 \in A_0$$

then

$$a_0 \notin B_0$$

or

$$A_0 \cap B_0 = \emptyset$$

**Proof**

Let $a_0$ and $b_0$ be

$$(A1.9) \quad a_0 = 2^{2\cdot 3^{n-1}(j_a-1)+k_n^a} \cdot K_O^a - 1 - 2^n \cdot \frac{2^{2\cdot 3^{n-1}\cdot j_a} - 1}{3^n}$$

$$(A1.10) \quad b_0 = 2^{2\cdot 3^{n-1}(j_b-1)+k_n^b} \cdot K_O^b - 1 - 2^n \cdot \frac{2^{2\cdot 3^{n-1}\cdot j_b} - 1}{3^n}$$

$a_0$ is a member of the subset:

$$a_0 \in A_0$$

and $b_0$ is a member of the subset:

$$b_0 \in B_0$$

We will prove that $a_0$ and $b_0$ cannot be the same number, in other words:

$$a_0 \neq b_0$$

Let's assume the opposite is true, in other words:

$$a_0 = b_0$$

Therefore, equating (A1.9) and (A1.10)

$$2^{2\cdot 3^{n-1}(j_a-1)+k_n^a} \cdot K_O^a - 1 - 2^n \cdot \frac{2^{2\cdot 3^{n-1}\cdot j_a} - 1}{3^n} = 2^{2\cdot 3^{n-1}(j_b-1)+k_n^b} \cdot K_n^b - 1 - 2^n \cdot \frac{2^{2\cdot 3^{n-1}\cdot j_b} - 1}{3^n}$$

We eliminate the "-1" in both sides of the identity

$$2^{2\cdot 3^{n-1}(j_a-1)+k_n^a} \cdot K_O^a - 2^n \cdot \frac{2^{2\cdot 3^{n-1}\cdot j_a} - 1}{3^n} = 2^{2\cdot 3^{n-1}(j_b-1)+k_n^b} \cdot K_n^b - 2^n \cdot \frac{2^{2\cdot 3^{n-1}\cdot j_b} - 1}{3^n}$$

We apply Lemma A1.1 to both sides of the identity above.

$$2^n \cdot \frac{2^{2\cdot 3^{n-1}\cdot j} - 1}{3^n} = 2^n \cdot \frac{2^{2\cdot 3^{n-1}} - 1}{3^n} \cdot [1 + 2^{2\cdot 3^{n-1}} + 2^{2\cdot 3^{n-1}\cdot 2} + \cdots + 2^{2\cdot 3^{n-1}\cdot (j-1)}]$$

Therefore

$$(A1.11) \quad 2^{2\cdot 3^{n-1}(j_a-1)+k_n^a} \cdot K_O^a - 2^n \cdot \frac{2^{2\cdot 3^{n-1}} - 1}{3^n} \cdot \left[1 + 2^{2\cdot 3^{n-1}} + 2^{2\cdot 3^{n-1}\cdot 2} + \cdots + 2^{2\cdot 3^{n-1}\cdot (j_a-1)}\right] =$$

$$= 2^{2\cdot 3^{n-1}(j_b-1)+k_n^b} \cdot K_n^b - 2^n \cdot \frac{2^{2\cdot 3^{n-1}} - 1}{3^n} \cdot \left[1 + 2^{2\cdot 3^{n-1}} + 2^{2\cdot 3^{n-1}\cdot 2} + \cdots + 2^{2\cdot 3^{n-1}\cdot (j_b-1)}\right]$$

Let's assume

$$j_a > j_b$$

then we can eliminate the common addends on both sides of the identity above

$$1 + 2^{2\cdot 3^{n-1}} + 2^{2\cdot 3^{n-1}\cdot 2} + \cdots + 2^{2\cdot 3^{n-1}\cdot (j_b-1)}$$

and rewrite the identity (A1.11) above as follows

$$2^{2\cdot 3^{n-1}(j_a-1)+k_n^a}\cdot K_O^a - 2^n\cdot\frac{2^{2\cdot 3^{n-1}}-1}{3^n}\cdot\left[2^{2\cdot 3^{n-1}\cdot(j_b)}+2^{2\cdot 3^{n-1}\cdot(j_b+1)}+\cdots+2^{2\cdot 3^{n-1}\cdot(j_a-1)}\right]$$
$$= 2^{2\cdot 3^{n-1}(j_b-1)+k_n^b}\cdot K_O^b$$

We simplify the above

$$2^{2\cdot 3^{n-1}(j_a-1)+k_n^a}\cdot K_O^a - 2^n\cdot\frac{2^{2\cdot 3^{n-1}}-1}{3^n}\cdot\sum_{t=j_b}^{t=j_a-1}2^{2\cdot 3^{n-1}\cdot t} = 2^{2\cdot 3^{n-1}(j_b-1)+k_n^b}\cdot K_O^b$$

We extract

$$2^{2\cdot 3^{n-1}\cdot j_b}$$

from the summation

$$2^{2\cdot 3^{n-1}(j_a-1)+k_n^a}\cdot K_O^a - 2^n\cdot\frac{2^{2\cdot 3^{n-1}}-1}{3^n}\cdot 2^{2\cdot 3^{n-1}\cdot j_b}\sum_{t=0}^{t=j_a-j_b-1}2^{2\cdot 3^{n-1}\cdot t} = 2^{2\cdot 3^{n-1}(j_b-1)+k_n^b}\cdot K_O^b$$

We divide all the members of the above by

$$2^{2\cdot 3^{n-1}(j_b-1)+k_n^b}$$

and obtain

$$2^{2\cdot 3^{n-1}(j_a-j_b)+k_n^a-k_n^b}\cdot K_O^a - 2^n\cdot\frac{2^{2\cdot 3^{n-1}}-1}{3^n}\cdot 2^{2\cdot 3^{n-1}-k_n^b}\cdot\sum_{t=0}^{t=j_a-j_b-1}2^{2\cdot 3^{n-1}\cdot t} = K_O^b$$

Which cannot be since the left side of the identity is an EVEN number and the right side of the identity is an ODD number. Therefore

$$a_0 \neq b_0$$

Which is what we wanted to prove

**Lemma A1.5**

Let $A_{j_a}^{j_a}$ be a set of odd numbers

$$A_{j_a}^{j_a} = \left\{2^{2\cdot 3^{n-1}(j_a-1)+k_n} \cdot K_O^a - 1 - 2^n \cdot \frac{2^{2\cdot 3^{n-1}\cdot \boldsymbol{j_a}} - 1}{3^n}\right\}$$

and let $B_{j_a}^{j_b}$ be another set of odd numbers

$$B_{j_a}^{j_b} = \left\{2^{2\cdot 3^{n-1}(j_a-1)+k_n} \cdot K_O^b - 1 - 2^n \cdot \frac{2^{2\cdot 3^{n-1}\cdot \boldsymbol{j_b}} - 1}{3^n}\right\}$$

Note that

$$\boldsymbol{j_b \neq j_a}$$

**ONLY** in the exponent of the numerator of the fraction and

$$j_b > j_a$$

Both sets represent the same set of odd numbers. In other words, if

$$a_0 \in A_{j_a}^{j_a}$$

then

$$a_0 \in B_{j_a}^{j_b}$$

or

$$A_{j_a}^{j_a} \cap B_{j_a}^{j_b} = A_{j_a}^{j_a} \cap \{odd\} = B_{j_a}^{j_b} \cap \{odd\}$$

**Proof**

Let $a_0$ be a member of the set $B_{j_a}^{j_b}$

$$(A1.12) \qquad a_0 = 2^{2\cdot 3^{n-1}(j_a-1)+k_n} \cdot K_O^b - 1 - 2^n \cdot \frac{2^{2\cdot 3^{n-1}\cdot j_b} - 1}{3^n}$$

Since we want to prove that $a_0$ is also a member of $A_{j_a}^{j_a}$

$$a_0 \in A_{j_a}^{j_a}$$

then we have to find a $K_O^a$ such that

$$(A1.13)\quad a_0 = 2^{2\cdot 3^{n-1}(j_a-1)+k_n} \cdot K_O^a - 1 - 2^n \cdot \frac{2^{2\cdot 3^{n-1}\cdot j_a} - 1}{3^n}$$

$$K_O^a \in \{odd\}$$

We apply Lemma A1.1 to (A1.12), therefore:

$$a_0 = 2^{2\cdot 3^{n-1}(j_a-1)+k_n} \cdot K_O^b - 1 - 2^n \cdot \frac{2^{2\cdot 3^{n-1}} - 1}{3^n} \cdot [1 + 2^{2\cdot 3^{n-1}} + 2^{2\cdot 3^{n-1}\cdot 2} + \cdots + 2^{2\cdot 3^{n-1}\cdot(j_b-1)}]$$

Since

$$j_b > j_a$$

we can divide the elements within the brackets in two groups and rewrite the above as

$$a_0 = 2^{2\cdot 3^{n-1}(j_a-1)+k_n} \cdot K_O^b - 1 - 2^n \cdot \frac{2^{2\cdot 3^{n-1}} - 1}{3^n} \cdot [1 + 2^{2\cdot 3^{n-1}} + 2^{2\cdot 3^{n-1}\cdot 2} + \cdots + 2^{2\cdot 3^{n-1}\cdot(j_a-1)}]$$

$$-2^n \cdot \frac{2^{2\cdot 3^{n-1}} - 1}{3^n} \cdot [2^{2\cdot 3^{n-1}\cdot(j_a)} + 2^{2\cdot 3^{n-1}\cdot(j_a+1)} + \cdots + 2^{2\cdot 3^{n-1}\cdot(j_b-1)}]$$

We apply Lemma A1.1 to the first group and rewrite the above as

$$a_0 = 2^{2\cdot 3^{n-1}(j_a-1)+k_n} \cdot K_O^b - 1 - 2^n \cdot \frac{2^{2\cdot 3^{n-1}\cdot j_a} - 1}{3^n}$$
$$- \left[2^n \cdot \frac{2^{2\cdot 3^{n-1}} - 1}{3^n} \cdot \left(2^{2\cdot 3^{n-1}\cdot(j_a)} + 2^{2\cdot 3^{n-1}\cdot(j_a+1)} + \cdots + 2^{2\cdot 3^{n-1}\cdot(j_b-1)}\right)\right]$$

We rewrite the second group as

$$a_0 = 2^{2\cdot 3^{n-1}(j_a-1)+k_n} \cdot K_O^b - 1 - \left[2^n \cdot \frac{2^{2\cdot 3^{n-1}\cdot j_a} - 1}{3^n} + 2^n \cdot \frac{2^{2\cdot 3^{n-1}} - 1}{3^n} \cdot \sum_{t=2\cdot 3^{n-1}\cdot j_a}^{t=2\cdot 3^{n-1}\cdot(j_b-1)} 2^t\right]$$

We extract

$$2^{2\cdot 3^{n-1}(j_a-1)+k_n}$$

from all the members of the summation

$$a_0 = 2^{2\cdot 3^{n-1}(j_a-1)+k_n} \cdot K_O^b - 1 - 2^n \cdot \frac{2^{2\cdot 3^{n-1}\cdot j_a} - 1}{3^n} +$$

$$-2^n \cdot \frac{2^{2\cdot 3^{n-1}} - 1}{3^n} \cdot 2^{2\cdot 3^{n-1}(j_a-1)+k_n} \cdot \sum_{t=2\cdot 3^{n-1}-k_n}^{t=2\cdot 3^{n-1}\cdot(j_b-j_a)-k_n} 2^t$$

We rearrange the above by combining all elements with the factor:

$$2^{2\cdot 3^{n-1}(j_a-1)+k_n}$$

and obtain

$$a_0 = 2^{2\cdot 3^{n-1}(j_a-1)+k_n} \cdot \left[ K_O^b - 2^n \cdot \frac{2^{2\cdot 3^{n-1}} - 1}{3^n} \cdot \sum_{t=2\cdot 3^{n-1}-k_n}^{t=2\cdot 3^{n-1}\cdot(j_b-j_a)-k_n} 2^t \right] - 1 - 2^n \cdot \frac{2^{2\cdot 3^{n-1}\cdot j_a} - 1}{3^n}$$

We make the following change in the identity above

$$K_O^a = K_O^b - 2^n \cdot \frac{2^{2\cdot 3^{n-1}} - 1}{3^n} \cdot \sum_{t=2\cdot 3^{n-1}-k_n}^{t=2\cdot 3^{n-1}\cdot(j_b-j_a)-k_n} 2^t$$

and obtain identity (A1.13),

$$(A1.13) \quad a_0 = 2^{2\cdot 3^{n-1}(j_a-1)+k_n} \cdot K_O^a - 1 - 2^n \cdot \frac{2^{2\cdot 3^{n-1}\cdot j_a} - 1}{3^n}$$

Therefore

$$a_0 \in A_{j_a}^{j_a}$$

which is what we wanted to prove.

**Lemma A1.6**

$$\forall n \in \mathbb{N}$$

all ODD numbers can be divided in $j$ subsets:

$$(A1.14) \qquad D_j = \left\{ 2^{2\cdot 3^{n-1}(j-1)+k_n} \cdot K_O - 1 - 2^n \cdot \frac{2^{2\cdot 3^{n-1}\cdot j} - 1}{3^n} \right\}$$

$$k_n \in \{1,2,\dots,2\cdot 3^{n-1}\}$$

$$K_O \in \{odd\}$$

$$j \in \mathbb{N}$$

and the union of the $j$ subsets $D_j$ makes the set of all odd numbers:

$$\{odd\} = \bigcup_{j=1}^{j=\infty} D_j$$

**Proof**

Lemma A1.2 proves that the union of the $j$ subsets $B_j$

$$(A1.6) \qquad B_j = \left\{ 2^{2\cdot 3^{n-1}(j-1)+k_n} \cdot K_O - 1 \right\}$$

makes the set of all odd numbers:

$$\{odd\} = \bigcup_{j=1}^{j=\infty} B_j$$

We apply Lemma A1.3 to the above. As a result, if we subtract the same EVEN number

$$K_E = 2^n \cdot \frac{2^{2\cdot 3^{n-1}\cdot w} - 1}{3^n}$$

to all members, then the union of all the $j$ subsets $D_j^w$

$$D_j^w = \left\{ 2^{2\cdot 3^{n-1}(j-1)+k_n} \cdot K_O - 1 - 2^n \cdot \frac{2^{2\cdot 3^{n-1}\cdot w} - 1}{3^n} \right\}$$

makes the set of all odd numbers:

$$\{odd\} = \bigcup_{j=1}^{j=\infty} D_j^w$$

We divide the above in the union of infinite subsets

$$\{odd\} = D_1^w \bigcup D_2^w \bigcup D_3^w \dots\dots \bigcup D_w^w \dots\dots \bigcup D_j^w \bigcup D_{j+1}^w \bigcup D_{j+2}^w \dots.$$

We apply Lemma A1.5 for all

$$w > j$$

Then the subsets $D_j^w$ and $D_j^j$ represent the same set of odd numbers,

$$D_j^w = D_j^j \;\; ; \quad \forall\, w > j$$

Therefore, we can substitute

$$D_j^j = \left\{ 2^{2\cdot 3^{n-1}(j-1)+k_n} \cdot K_O - 1 - 2^n \cdot \frac{2^{2\cdot 3^{n-1}\cdot j} - 1}{3^n} \right\}$$

for

$$D_j^w = \left\{ 2^{2\cdot 3^{n-1}(j-1)+k_n} \cdot K_O - 1 - 2^n \cdot \frac{2^{2\cdot 3^{n-1}\cdot w} - 1}{3^n} \right\}$$

As a result, the odd numbers can be divided into the union of the following subsets:

$$\{odd\} = D_1^1 \bigcup D_2^2 \bigcup D_3^3 \dots. \bigcup D_j^j \bigcup D_{j+1}^w \bigcup D_{j+2}^w \dots.$$

If we bring w to infinity; then

$$\{odd\} = \lim_{w\to\infty} D_1^1 \bigcup D_2^2 \bigcup D_3^3 \dots\dots \bigcup D_j^j \bigcup D_{j+1}^{j+1} \dots.$$

then

$$\{odd\} = \bigcup_{j=1}^{j=\infty} \left\{ 2^{2\cdot 3^{n-1}(j-1)+k_n} \cdot K_O - 1 - 2^n \cdot \frac{2^{2\cdot 3^{n-1}\cdot j} - 1}{3^n} \right\}$$

or

$$(A1.15) \quad D_j^j = \left\{ 2^{2\cdot 3^{n-1}(j-1)+k_n} \cdot K_O - 1 - 2^n \cdot \frac{2^{2\cdot 3^{n-1}\cdot j} - 1}{3^n} \right\}$$

$$\{odd\} = \bigcup_{j=1}^{j=\infty} D_j^j$$

Which is what we wanted to prove

Lemma A1.4 states that all the above subsets are disjoint,

$$D_k^k \cap D_l^l = \emptyset$$
$$\forall k \neq l$$

In other words, no odd number is repeated in two subsets $D_j^j$ – see (A1.13) - at the same time. Therefore, we can state that, for all $n \in \mathbb{N}$, the odd numbers can unambiguously be written as:

$$a_O = 2^{2\cdot 3^{n-1}(j-1)+k_n} \cdot K_O - 2^n \cdot \frac{2^{2\cdot 3^{n-1}\cdot j} - 1}{3^n} - 1$$

Therefore, any odd number can be written as:

$$(A1.16) \quad a_O = 2^{2\cdot 3^{n-1}(j-1)+k_n} \cdot K_O - 2^n \cdot \frac{2^{2\cdot 3^{n-1}\cdot j} - 1}{3^n} - 1$$
$$\forall n, j \in \mathbb{N}$$
$$k_n \in \{1,2, \ldots ,2 \cdot 3^{n-1}\}$$
$$K_O \in \{odd\}$$

**Lemma A1.7**

$$\forall n \in \mathbb{N}$$

All ODD numbers can be expressed as:

$$(A1.1) \quad a_O = 2^{\alpha} \cdot Q_O - \frac{2^{3^{n-1} \cdot j_o^{\beta}} + 1}{3^n} \cdot (3^n - 2^n)$$

$$j_o^{\beta}, Q_O \in \{odd\}$$

$$\alpha \in \mathbb{N}$$

$$3^{n-1} \cdot j_o^{\beta} > \alpha$$

**Proof**

We have already proved that all odd numbers can be expressed using Identity (A1.16)

$$(A1.16) \quad a_O = 2^{2 \cdot 3^{n-1}(j-1)+k_n} \cdot K_O - 2^n \cdot \frac{2^{2 \cdot 3^{n-1} \cdot j} - 1}{3^n} - 1$$

Let

$$\alpha \in \mathbb{N}$$

be

$$\alpha = 2 \cdot 3^{n-1}(j-1) + k_n$$

and apply the above to Identity (A1.16) above and obtain

$$a_O = 2^{\alpha} \cdot K_O - 2^n \cdot \frac{2^{2 \cdot 3^{n-1} \cdot j} - 1}{3^n} - 1$$

We add and subtract to the identity above

$$\frac{2^{3^{n-1} \cdot j_o^{\beta}} + 1}{3^n} \cdot (3^n - 2^n)$$

Please note that $j_o^{\beta}$ can be any odd number that complies with the following:

$$3^{n-1} \cdot j_o^{\beta} > \alpha$$

and obtain

$$a_O = 2^\alpha \cdot K_O - 2^n \cdot \frac{2^{2\cdot 3^{n-1}\cdot j} - 1}{3^n} - 1 + \frac{2^{3^{n-1}\cdot j_o^\beta} + 1}{3^n} \cdot (3^n - 2^n) - \frac{2^{3^{n-1}\cdot j_o^\beta} + 1}{3^n} \cdot (3^n - 2^n)$$

We break the brackets

$$a_O = 2^\alpha \cdot K_O - 2^n \cdot \frac{2^{2\cdot 3^{n-1}\cdot j} - 1}{3^n} - 1 + 3^n \cdot \frac{2^{3^{n-1}\cdot j_o^\beta} + 1}{3^n} - 2^n \cdot \frac{2^{3^{n-1}\cdot j_o^\beta} + 1}{3^n}$$
$$- \frac{2^{3^{n-1}\cdot j_o^\beta} + 1}{3^n} \cdot (3^n - 2^n)$$

We simplify the above

$$a_O = 2^\alpha \cdot K_O - 2^n \cdot \frac{2^{3^{n-1}\cdot j_o^\beta} + 2^{2\cdot 3^{n-1}\cdot j}}{3^n} + 2^{3^{n-1}\cdot j_o^\beta} - \frac{2^{3^{n-1}\cdot j_o^\beta} + 1}{3^n} \cdot (3^n - 2^n)$$

Since $j_o^\beta$ can be any odd number, we can extract $2^{2\cdot 3^{n-1}\cdot j}$ from the fraction in the identity above

$$(A1.17) \quad a_O = 2^\alpha \cdot K_O - 2^{2\cdot 3^{n-1}\cdot j + n} \frac{2^{3^{n-1}\cdot\left(j_o^\beta - 2j\right)} + 1}{3^n} + 2^{3^{n-1}\cdot j_o^\beta} - \frac{2^{3^{n-1}\cdot j_o^\beta} + 1}{3^n} \cdot (3^n - 2^n)$$

Since

$$k_n \in \{1; 2; \dots; 2 \cdot 3^{n-1}\}$$

then

$$2 \cdot 3^{n-1} \cdot j > 2 \cdot 3^{n-1}(j-1) + k_n = \alpha$$

Also, since $j_o^\beta$ can be any odd number that complies with

$$3^{n-1} \cdot j_o^\beta > \alpha$$

then we can extract $2^\alpha$ from the first three elements in the right side of Identity (A1.17)

$$a_O = 2^{\alpha} \cdot \left[ K_O - 2^{2\cdot 3^{n-1}\cdot j+n-\alpha} \frac{2^{3^{n-1}\cdot\left(j_o^{\beta}-2j\right)}+1}{3^n} + 2^{3^{n-1}\cdot j_o^{\beta}-\alpha} \right] - \frac{2^{3^{n-1}\cdot j_o^{\beta}}+1}{3^n} \cdot (3^n - 2^n)$$

We make the following changes:

$$Q_O = K_O - 2^{2\cdot 3^{n-1}\cdot j+n-\alpha} \frac{2^{3^{n-1}\cdot\left(j_o^{\beta}-2j\right)}+1}{3^n} + 2^{3^{n-1}\cdot j_o^{\beta}-\alpha}$$

$$Q_O \in \{odd\}$$

and obtain

$$(A1.1) \quad a_O = 2^{\alpha} \cdot Q_O - \frac{2^{3^{n-1}\cdot j_o^{\beta}}+1}{3^n} \cdot (3^n - 2^n)$$

$$j_o^{\beta}, Q_O \in \{odd\}$$

$$\alpha, n \in \mathbb{N}$$

$$3^{n-1} \cdot j_o^{\beta} > \alpha$$

Which is what we wanted to prove.

Since we are dealing with positive odd numbers, $Q_O$ must be large enough to make $a_O$ a positive integer.

Please note that $j_o^{\beta}$ can be <u>any odd number</u> as long as it complies with the restriction:

$$3^{n-1} \cdot j_o^{\beta} > \alpha$$

We obtain the same odd number $a_O$ by "adapting" $Q_O$ to the new value $j_o^{\beta}$.

# Appendix 2. ODD Number Grade

In the previous Appendix 1, we prove that, for all n, with no restrictions, we can express any odd number with Identity (A1.1)

In this appendix, we restrict the use of the parameters $\alpha$ and n in Identity (A1.1):

$$\alpha > n$$

and still retain all odd numbers. We do so to adapt Identity (A1.1) to the rules of how the sequence is generated. Doing so, it allows us to break down a sequence in many cycles, where the beginning and the end of each cycle is an odd number.

In this appendix and the following Appendix 3 we prove that any sequence can be divided in cycles. Each cycle begins with an odd number and climbs up alternating upward and downward steps until it reaches an upper bound. The upper bound is reached when the sequence meets two or more consecutive downward steps. At this stage, the sequence descends following a series of even numbers – downward steps – until it reaches another odd number. This marks the end of the cycle.

**Definition A2.1**

We define the grade of an ODD number – $n_g$ – as the number of upward steps until a cycle finds an upper bound. The upper bound is reached when the sequence meets two or more consecutive downward steps.

For any ODD number $a_O$, the grade $n_g$ is given by the following formula:

$$(A2.1) \qquad a_O = 2^{n_g} \cdot K_O - 1$$

$$K_O \in \{odd\}$$

The upper bound of the cycle is

$$a_{u.b} = 3^{n_g} \cdot 2 \cdot K_O - 2$$

This upper bound is found after $n_g$ upward steps

Let $a_O$ be the initial number of a sequence

$$a_O = 2^{n_g} \cdot K_O - 1$$

$$K_O \in \{odd\}$$

The sequence generated by $a_0$ is as follows:

$$(A2.1) \;\; a_O = 2^{n_g} \cdot K_O - 1 \qquad ODD \text{ (Initial value)}$$

since $a_O = ODD$, then $a_1 = 3 \cdot a_O + 1$

$$a_1 = 3 \cdot 2^{n_g} \cdot K_O - 2 \quad EVEN$$

since $a_1 = EVEN$, then $a_2 = \frac{a_1}{2}$

$$a_2 = 3 \cdot 2^{n_g - 1} \cdot K_O - 1 \quad ODD$$

since $a_2 = ODD$, then $a_3 = 3 \cdot a_2 + 1$

$$a_3 = 3^2 \cdot 2^{n_g - 1} \cdot K_O - 2 \quad EVEN$$

since $a_3 = EVEN$, then $a_4 = \frac{a_3}{2}$

$$a_4 = 3^2 \cdot 2^{n_g - 2} \cdot K_O - 1 \quad ODD$$

we repeat the process $n_g$ times until we find two consecutive EVEN numbers. The first of these two consecutive EVEN numbers is the upper bound of the cycle generated by the initial ODD number.

$$(A2.2) \qquad a_{u.b} = 3^{n_g} \cdot 2 \cdot K_O - 2 \quad EVEN - UPPER\ BOUND$$

Now, at least two downward steps follow until we find another ODD number.

First downward step

$$a_{u.b+1} = \frac{a_{u.b}}{2} = 3^{n_g} \cdot K_O - 1$$

Second downward step

$$a_{u.b+2} = \frac{a_{u.b+1}}{2} = \frac{3^{n_g} \cdot K_O - 1}{2}$$

More downward steps may follow, depending on $K_O$.

For example, let's take the number

$$a_0 = 739$$

using (A2.1), we know the grade of this number is 2

$$a_O = 2^2(2 \cdot 92 + 1) - 1 = 739$$

using (A2.2) we know the upper bound is

$$a_{u.b} = 3^2 \cdot 2 \cdot (2 \cdot 92 + 1) - 2 = 3328$$

which is an EVEN number. Now, a downward path follows until we find an ODD number again which is the final element of the cycle:

$$a_F = 13$$

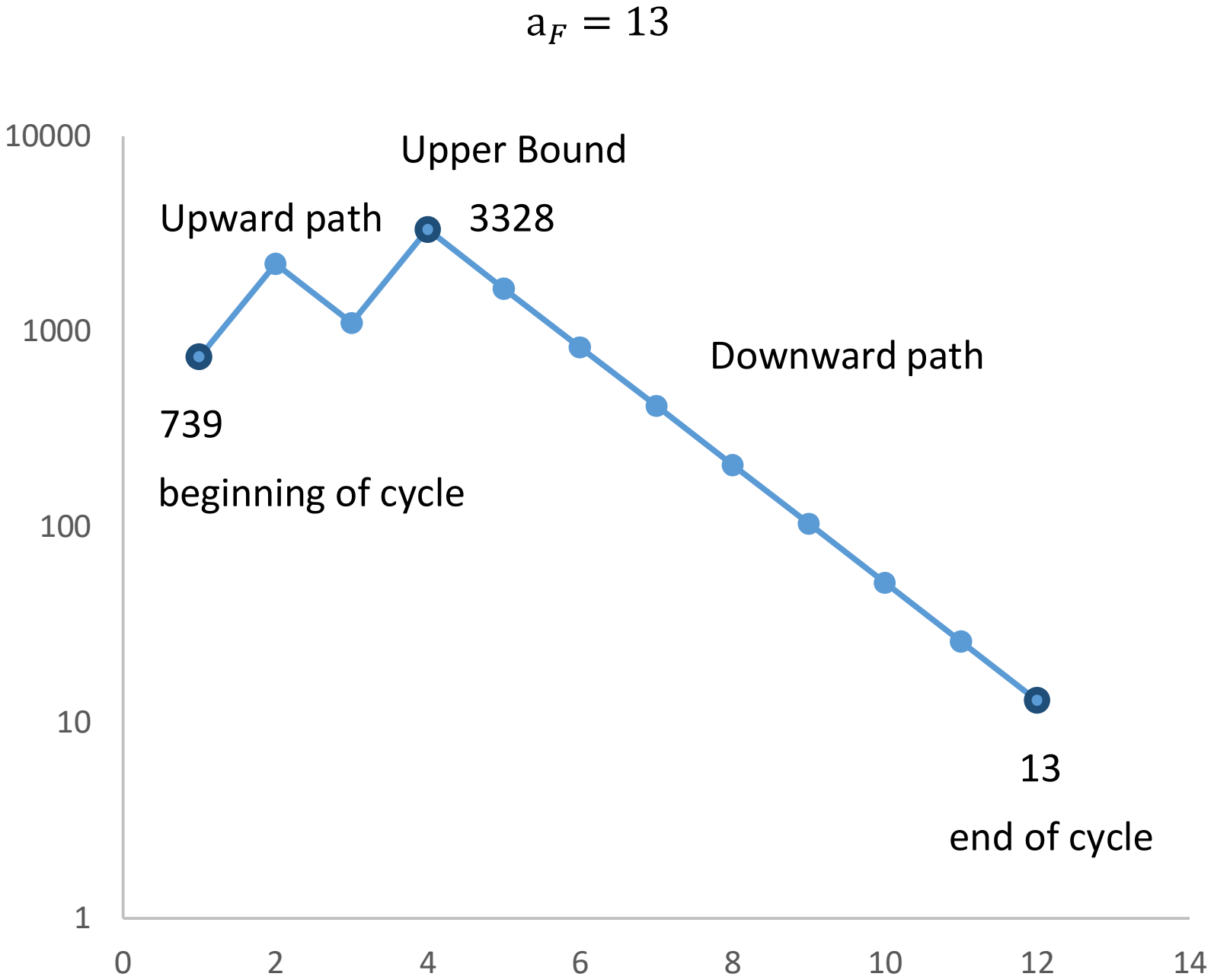

Note: The graph is in logarithm scale to facilitate the visualization of upward and downward steps.

From Appendix 1 we know that we can express any odd number as:

$$(A1.1) \quad a_O = 2^{\alpha} \cdot Q_O - \frac{2^{3^{n-1} \cdot j_o^{\beta}} + 1}{3^n} \cdot (3^n - 2^n)$$

$$3^{n-1} \cdot j_o^{\beta} > \alpha$$

**Lemma A2.2**

If

$$\alpha < n$$

then the grade of the number $a_O$ is $\alpha$

$$n_g = \alpha$$

**Proof**

We use Identity (A1.1) from Appendix 1. Any ODD number can be written as:

$$(A1.1) \quad a_O = 2^{\alpha} \cdot Q_O - \frac{2^{3^{n-1} \cdot j_o^{\beta}} + 1}{3^n} \cdot (3^n - 2^n)$$

$$Q_O \in \{odd\}$$

We break the brackets and rearrange the identity above

$$(A2.1) \quad a_O = 2^{\alpha} \cdot Q_O - 2^{3^{n-1} \cdot j_o^{\beta}} + 2^n \cdot \frac{2^{3^{n-1} \cdot j_o^{\beta}} + 1}{3^n} - 1$$

since

$$\alpha < n$$

and

$$3^{n-1} \cdot j_o^{\beta} > \alpha$$

we can extract $2^{\alpha}$ from the first three elements of Identity (A2.1) above

$$a_O = 2^{\alpha} \cdot \left[ Q_O - 2^{3^{n-1} \cdot j_o^{\beta} - \alpha} + 2^{n-\alpha} \cdot \frac{2^{3^{n-1} \cdot j_o^{\beta}} + 1}{3^n} \right] - 1$$

then the element within the brackets of the above expression is an ODD number. Therefore, the above identity has a similar structure as (A2.1)

$$(A2.1) \qquad a_O = 2^{n_g} \cdot K_O - 1$$

thus, as we have seen in Definition A2.1, $\alpha$ is the grade of $a_O$,

$$n_g = \alpha$$

which is what we wanted to prove.

**Lemma A2.3**

If

$$\alpha = n$$

then the grade of the number $a_O$ is larger than $\alpha$

$$n_g > \alpha$$

**Proof**

We use Identity (A2.1) above

$$(A2.1) \quad a_O = 2^{\alpha} \cdot Q_O - 2^{3^{n-1} \cdot j_o^{\beta}} + 2^n \cdot \frac{2^{3^{n-1} \cdot j_o^{\beta}} + 1}{3^n} - 1$$

Since

$$\alpha = n$$

We can extract $2^n$ from the first three elements of the expression above:

$$a_O = 2^n \cdot \left( Q_O - 2^{3^{n-1} \cdot j_o^{\beta} - n} + \frac{2^{3^{n-1} \cdot j_o^{\beta}} + 1}{3^n} \right) - 1$$

since $Q_O$ and $\frac{2^{3^{n-1}\cdot j_o^{\beta}}+1}{3^n}$ are both ODD, then the element within brackets of the above expression is an EVEN number. As a result, it can be rearranged as follows:

$$(A2.3) \qquad a_O = 2^{n+p} \cdot \frac{\left(Q_O + \frac{2^{3^{n-1}\cdot j_o^{\beta}}+1}{3^n} - 2^{3^{n-1}\cdot j_o^{\beta}-n}\right)}{2^p} - 1$$

$$p \geq 1$$

the above identity has a similar structure as (A2.1)

$$(A2.1) \qquad a_O = 2^{n_g} \cdot K_O - 1$$

thus, as we have seen in Definition A2.1, "n+p" is the grade of $a_O$

$$n + p = n_g$$

since

$$p \geq 1$$

then

$$n_g > n = \alpha$$

which is what we wanted to prove.

**Lemma A2.4**

If

$$\alpha > n$$

then the grade of the number $a_O$ is $n$

$$n_g = n$$

**Proof**

We use Identity (A2.1) above

$$(A2.1) \quad a_O = 2^{\alpha} \cdot Q_O - 2^{3^{n-1} \cdot j_o^{\beta}} + 2^n \cdot \frac{2^{3^{n-1} \cdot j_o^{\beta}} + 1}{3^n} - 1$$

Since

$$3^{n-1} \cdot j_o^{\beta} > \alpha > n$$

Then we can extract $2^n$ in the fraction of identity (A2.1) above

$$a_O = 2^n \cdot \left( 2^{\alpha - n} \cdot Q_O - 2^{3^{n-1} \cdot j_o^{\beta} - n} + \frac{2^{3^{n-1} \cdot j_o^{\beta}} + 1}{3^n} \right) - 1$$

since the element within the brackets of the expression above is ODD, then the above expression has a similar structure as (A2.1)

$$(A2.1) \qquad a_O = 2^{n_g} \cdot K_O - 1$$

thus, as we have seen in Definition A2.1, n is the grade of $a_O$,

$$n = n_g$$

which is what we wanted to prove.

As a result of the above, if we want identity (A1.1) to represent **Grade $n_g$** numbers ONLY, we have to apply the following restriction:

$$\alpha > n$$

$$(A1.1) \qquad a_O = 2^{\alpha} \cdot Q_O - \frac{2^{3^{n-1} \cdot j_o^{\beta}} + 1}{3^n} \cdot (3^n - 2^n)$$

$$3^{n-1} \cdot j_o^{\beta} > \alpha$$

$$Q_O \in \{odd\}$$

$$(A2.4) \qquad \alpha > n$$

# Appendix 3. ODD Number Evolution – End of Cycle

Appendix 1 proved that any ODD number can be written using Identity (A1.1). Also, as we have seen in Appendix 2, ODD numbers can be written using different n values but an $n_g$ always exists that represents the grade of said ODD number.

For all purposes **in the following appendix, and the rest of this paper for that matter**, whenever we use identity (A1.1), we assume that **the chosen value n is $n_g$**. Therefore, we will be using identity (A1.1) with its applicability limitations.

$$(A1.1) \qquad a_O = 2^{\alpha} \cdot Q_O - \frac{2^{3^{n-1} \cdot j_o^{\beta}} + 1}{3^n} \cdot (3^n - 2^n)$$

$$j_o^{\beta}, Q_O \in \{odd\}$$
$$\alpha, n \in \mathbb{N}$$
$$3^{n-1} \cdot j_o^{\beta} > \alpha > n$$

In Appendix 2 we also define the beginning and the end of a cycle within the sequence. This appendix, we prove that the final element of this cycle, after n upward steps and $\alpha$ downward steps, is

$$a_F = 3^n \cdot Q_O - 2^{3^{n-1} \cdot j_o^{\beta} - \alpha} \cdot (3^n - 2^n)$$

**Lemma A3.1**

Let $a_O$ be the initial element of a cycle

$$(A1.1) \qquad a_O = 2^{\alpha} \cdot Q_O - \frac{2^{3^{n-1} \cdot j_o^{\beta}} + 1}{3^n} \cdot (3^n - 2^n)$$

$$j_o^{\beta}, Q_O \in \{odd\}$$
$$\alpha, n \in \mathbb{N}$$
$$(A3.1) \qquad 3^{n-1} \cdot j_o^{\beta} > \alpha > n$$

The final element of this cycle is:

$$(A3.2) \qquad a_F = 3^n \cdot Q_O - 2^{3^{n-1} \cdot j_o^{\beta} - \alpha} \cdot (3^n - 2^n)$$

## Proof

We start with identity (A1.1).

$$a_O = 2^{\alpha} \cdot Q_O - \frac{2^{3^{n-1} \cdot j_o^{\beta}} + 1}{3^n} \cdot (3^n - 2^n)$$

To help visualize the evolution of the sequence, we break the brackets above and obtain:

$$a_O = 2^{\alpha} \cdot Q_O - 2^{3^{n-1} \cdot j_o^{\beta}} + 2^n \cdot \frac{2^{3^{n-1} \cdot j_o^{\beta}} + 1}{3^n} - 1$$

Since $a_O$ is ODD, the second element of the sequence is:

$$a_2 = 3 \cdot a_O + 1$$

Therefore

$$a_2 = 3 \cdot \left[ 2^{\alpha} \cdot Q_O - 2^{3^{n-1} \cdot j_o^{\beta}} + 2^n \cdot \frac{2^{3^{n-1} \cdot j_o^{\beta}} + 1}{3^{n-1}} - 1 \right] + 1$$

We break the bracket and multiply each element within the bracket by 3

$$a_2 = 3 \cdot 2^{\alpha} \cdot Q_O - 3 \cdot 2^{3^{n-1} \cdot j_o^{\beta}} + 2^n \cdot \frac{2^{3^{n-1} \cdot j_o^{\beta}} + 1}{3^{n-1}} - 3 + 1$$

We simplify the above by adding the last two elements of the expression above:

$$a_2 = 3 \cdot 2^{\alpha} \cdot Q_O - 3 \cdot 2^{3^{n-1} \cdot j_o^{\beta}} + 2^n \cdot \frac{2^{3^{n-1} \cdot j_o^{\beta}} + 1}{3^{n-1}} - 2$$

Since $a_2$ is EVEN, then the next element of the sequence is

$$a_3 = \frac{a_2}{2}$$

Therefore:

$$a_3 = \frac{3 \cdot 2^{\alpha} \cdot Q_O - 3 \cdot 2^{3^{n-1} \cdot j_o^{\beta}} + 2^n \cdot \frac{2^{3^{n-1} \cdot j_o^{\beta}} + 1}{3^{n-1}} - 2}{2}$$

thus

$$a_3 = 3 \cdot 2^{\alpha-1} \cdot Q_O - 3 \cdot 2^{3^{n-1} \cdot j_o^\beta - 1} + 2^{n-1} \cdot \frac{2^{3^{n-1} \cdot j_o^\beta} + 1}{3^{n-1}} - 1 = ODD$$

$a_3$ is ODD again, therefore, the following element of the sequence is:

$$a_4 = 3 \cdot a_3 + 1$$

therefore

$$a_4 = 3^2 \cdot 2^{\alpha-1} \cdot Q_O - 3^2 \cdot 2^{3^{n-1} \cdot j_o^\beta - 1} + 2^{n-1} \cdot \frac{2^{3^{n-1} \cdot j_o^\beta} + 1}{3^{n-2}} - 3 + 1$$

we do as before, we add the last two elements of the expression above:

$$a_4 = 3^2 \cdot 2^{\alpha-1} \cdot Q_O - 3^2 \cdot 2^{3^{n-1} \cdot j_o^\beta - 1} + 2^{n-1} \cdot \frac{2^{3^{n-1} \cdot j_o^\beta} + 1}{3^{n-2}} - 2 = EVEN$$

since $a_4$ is EVEN, the next element of the sequence is:

$$a_5 = \frac{a_4}{2} = 3^2 \cdot 2^{\alpha-2} \cdot Q_O - 3^2 \cdot 2^{3^{n-1} \cdot j_o^\beta - 2} + 2^{n-2} \cdot \frac{2^{3^{n-1} \cdot j_o^\beta} + 1}{3^{n-2}} - 1 = ODD$$

We repeat this process n-1 times:

$$a_{2n-1} = 3^{n-1} \cdot 2^{\alpha-(n-1)} \cdot Q_O - 3^{n-1} \cdot 2^{3^{n-1} \cdot j_o^\beta - (n-1)} + 2^{n-(n-1)} \cdot \frac{2^{3^{n-1} \cdot j_o^\beta} + 1}{3^{n-(n-1)}} - 1$$

we simplify the above:

$$a_{2n-1} = 3^{n-1} \cdot 2^{\alpha-(n-1)} \cdot Q_O - 3^{n-1} \cdot 2^{3^{n-1} \cdot j_o^\beta - (n-1)} + 2^1 \cdot \frac{2^{3^{n-1} \cdot j_o^\beta} + 1}{3^1} - 1$$

this is an odd number, therefore the next element of the sequence is:

$$a_{2n} = 3 \cdot a_{2n-1} + 1$$

$$a_{2n} = 3^n \cdot 2^{\alpha-(n-1)} \cdot Q_O - 3^n \cdot 2^{3^{n-1} \cdot j_o^\beta - (n-1)} + 2^1 \cdot \left(2^{3^{n-1} \cdot j_o^\beta} + 1\right) - 3 + 1$$

we simplify the above

$$(A3.3) \qquad a_{up} = 3^n \cdot 2^{\alpha-(n-1)} \cdot Q_O - 3^n \cdot 2^{3^{n-1}\cdot j_o^{\beta}-(n-1)} + 2^{3^{n-1}\cdot j_o^{\beta}+1} = EVEN$$

$a_{up}$ is an EVEN number. At this point we have reached the **upper bound** of the cycle.

As we have seen at the beginning of the appendix, for the purpose of this document, n represents the grade of the odd number:

$$n = n_g$$

therefore, identity (A3.1) applies:

$$(A2.5) \quad 3^{n-1} \cdot j_o^{\beta} > \alpha > n$$

$$\forall\, n$$

Therefore, we can extract $2^{\alpha-(n-1)}$ from all the elements of the right side of the identity and obtain:

$$a_{u.b} = 2^{\alpha-(n-1)} \cdot \left(3^n \cdot Q_O - 3^n \cdot 2^{3^{n-1}\cdot j_o^{\beta}-\alpha} + 2^{3^{n-1}\cdot j_o^{\beta}-\alpha+n}\right) = EVEN$$

We rearrange the elements within the brackets and extract the common factor $2^{3^{n-1}\cdot j_o^{\beta}-\alpha}$

$$(A3.4) \qquad a_{u.b} = 2^{\alpha-(n-1)} \cdot \left(3^n \cdot Q_O - 2^{3^{n-1}\cdot j_o^{\beta}-\alpha} \cdot (3^n - 2^n)\right) = EVEN$$

At this point, the cycle enters the downward trajectory. As a result, "$\alpha - (n-1)$" downward steps follow until we find another ODD number.

Since $a_{up}$ is an EVEN number, the following element of the sequence is obtained by dividing by 2. Therefore, we use (A3.4) above and divide by 2

$$a_{2n+1} = \frac{a_{up}}{2} = 2^{\alpha-n} \cdot \left[3^n \cdot Q_O - 2^{3^{n-1}\cdot j_o^{\beta}-\alpha} \cdot (3^n - 2^n)\right] = EVEN$$

Since the number above is EVEN, the following element of the sequence is again obtained by diving by 2.

$$a_{2n+2} = \frac{a_{2n+1}}{2} = 2^{\alpha-n-1} \cdot \left[3^n \cdot Q_O - 2^{3^{n-1} \cdot j_o^\beta - \alpha} \cdot (3^n - 2^n)\right] = EVEN$$

We repeat the process "$\alpha - (n-1)$" times until the sequence produces an odd number

$$(A3.2) \qquad a_F = 3^n \cdot Q_O - 2^{3^{n-1} \cdot j_o^\beta - \alpha} \cdot (3^n - 2^n) = ODD$$

which is what we wanted to prove. This is the end of the cycle.

Therefore, the evolution of any odd number always follows the same pattern:

- First, there is an upward trajectory; an upward step (3x+1) and a downward step (divided by 2) alternate n times until they reach an upper bound.

- Once they reach this upper bound, then the downward trajectory begins. There are "$\alpha - (n-1)$" consecutive downward steps (divided by two) until the sequence reaches an odd number again.

## Single-Cycle Sequences. Closing Cycle

There are sequences that converge to "1" in one single cycle. This is the case when the final element of the cycle is 1.

Identity (A3.2) above gives us the final element of any cycle

$$(A3.2) \qquad a_F = 3^n \cdot Q_O - 2^{3^{n-1} \cdot j_o^\beta - \alpha} \cdot (3^n - 2^n)$$

since we want this final element to be 1

$$(A3.8) \qquad a_F = 1$$

This always occurs for

$$(A3.5) \qquad \alpha = 3^{n-1} \cdot \left(j_o^\beta - q_E\right) + n$$

$$q_E \in \{even\}$$

and

$$(A3.6)^* \quad Q_O = 2^{3^{n-1} \cdot q_E - n} - \frac{2^{3^{n-1} \cdot q_E} - 1}{3^n}$$

(*) In Appendix 4 we prove that $2^{3^{n-1}\cdot q_E} - 1$ is always evenly divisible by $3^n$

The last and final cycle of any sequence is called the Closing Cycle. Any sequence has a closing cycle whose initial element is given by Identity (A3.9) below.

**Lemma A3.2**

Any sequence generated by the ODD number:

$$a_O = 2^n \cdot \frac{2^{3^{n-1}\cdot q_O} + 1}{3^n} - 1$$

Converges to 1 in one single cycle.

$$a_F = 1$$

## Proof

If we apply (A3.5) and (A3.6) to identity (A3.2)

$$a_F = 3^n \cdot \left(2^{3^{n-1}\cdot q_E - n} - \frac{2^{3^{n-1}\cdot q_E} - 1}{3^n}\right) - 2^{3^{n-1}\cdot j_o^{\beta} - \left(3^{n-1}\cdot\left(j_o^{\beta} - q_E\right) + n\right)} \cdot (3^n - 2^n)$$

We break the brackets and simplify the power of 2 above.

$$a_F = 3^n \cdot 2^{3^{n-1}\cdot q_E - n} - 2^{3^{n-1}\cdot q_E} + 1 - 3^n \cdot 2^{3^{n-1}\cdot q_E - n} + 2^{3^{n-1}\cdot q_E}$$

We simplify the above

$$a_F = 1$$

which is what we wanted to prove.

We apply (A3.5) and (A3.6) to identity (A1.1). The odd number that generates said single cycle is

$$a_O = 2^{3^{n-1}\cdot\left(j_o^{\beta} - q_E\right) + n} \cdot \left[2^{3^{n-1}\cdot q_E - n} - \frac{2^{3^{n-1}\cdot q_E} - 1}{3^n}\right] - \frac{2^{3^{n-1}\cdot j_o^{\beta}} + 1}{3^n} \cdot (3^n - 2^n)$$

We break the brackets above

$$a_O = 2^{3^{n-1}\cdot j_o^{\beta}} - 2^n \cdot \frac{2^{3^{n-1}\cdot j_o^{\beta}} - 2^{3^{n-1}\cdot\left(j_o^{\beta} - q_E\right)}}{3^n} - 2^{3^{n-1}\cdot j_o^{\beta}} - 1 + 2^n \cdot \frac{2^{3^{n-1}\cdot j_o^{\beta}} + 1}{3^n}$$

We simplify the above

$$a_O = 2^n \cdot \frac{2^{3^{n-1}\cdot\left(j_o^{\beta} - q_E\right)} + 1}{3^n} - 1$$

We make the following change

$$q_O = j_o^{\beta} - q_E$$

$$q_O \in \{odd\}$$

and simplify the identity above. We obtain:

$$(A3.7) \qquad a_O = 2^n \cdot \frac{2^{3^{n-1}\cdot q_O} + 1}{3^n} - 1$$

Which is what we wanted to prove.

The total number of downward steps throughout this cycle – upward and downward trajectories - is:

$$(A3.8) \qquad \alpha = 3^{n-1} \cdot q_O + n$$

The number of upward steps – upward trajectory - is:

$$n$$

## Example 1 – Simple Cycle

We use the same example shown of the introduction:

$$a_O = 739$$

In this case, using Identity (A1.1):

| n | 2 |
|---|---|
| $\alpha$ | 9 |
| $K_O$ | 37 |
| $j_O^{\beta}$ | 5 |

Therefore, this number can also be written as:

$$a_O = 2^9 \cdot 37 - \frac{2^{3^{(2-1)} \cdot 5} + 1}{3^2} \cdot (3^2 - 2^2) = 739$$

and using (A3.1) the final element of the cycle is:

$$\mathrm{a}_F = 3^2 \cdot 37 - 2^{3^{(2-1)} \cdot 5 - 9} \cdot (3^2 - 2^2) = 13$$

and the cycle follows the established pattern

- First there is the upward trajectory, (in this case n=2)

  739-2218-1109-3328

- Second, there is a downward trajectory of 7 steps down

  $$\alpha - n + 1 = 9 - 2 + 1 = 8$$

  1664-832-416-208-104-52-26-13

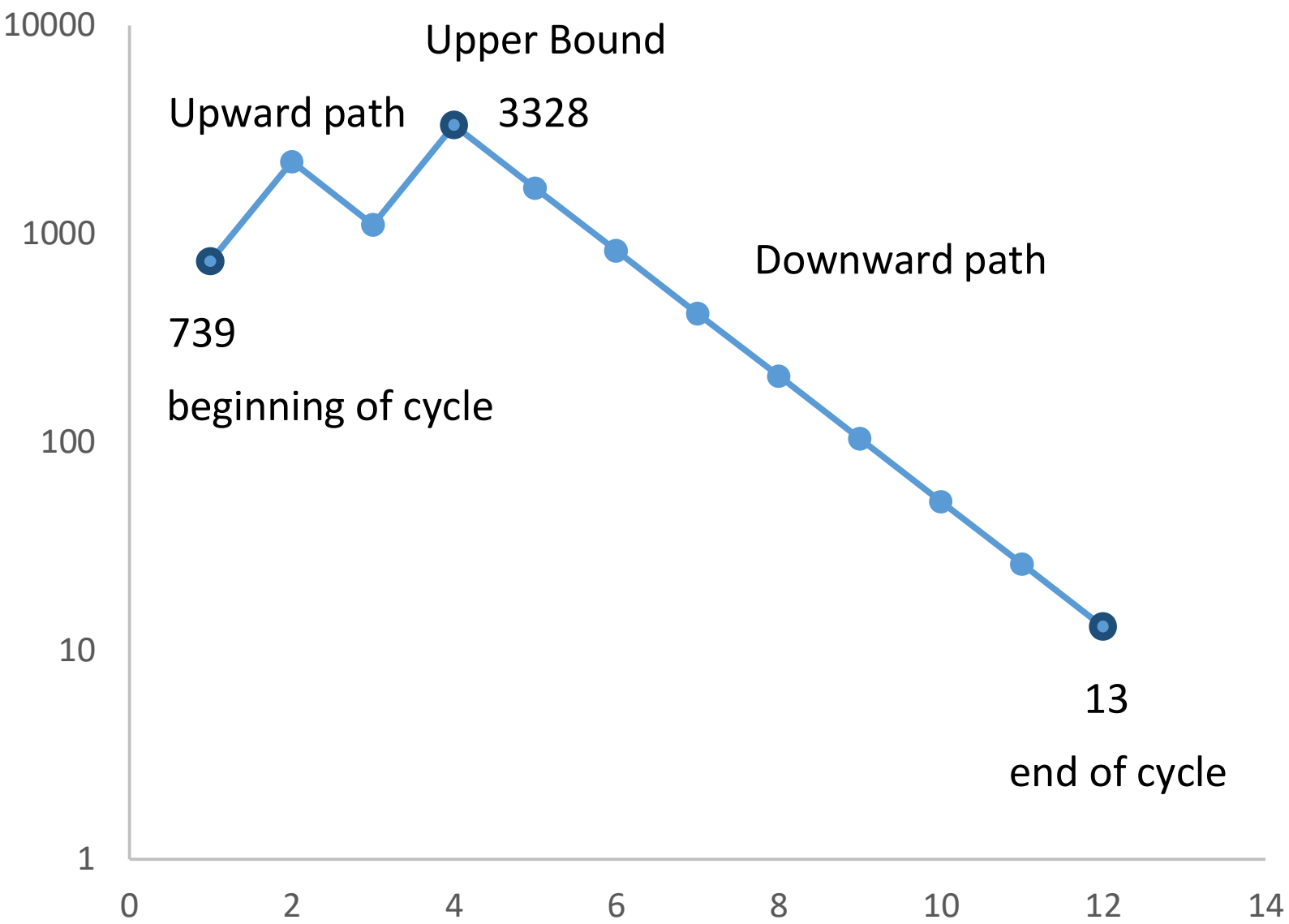

Note: The graph is in logarithm scale to facilitate the visualization of upward and downward steps.

## Example 2 – Closing Cycle

Let

$$q_O = 3 \qquad \text{and} \qquad n=3$$

we apply the above to Identity (A3.7)

$$a_O = 2^3 \cdot \frac{2^{3^{3-1}\cdot(3)} + 1}{3^3} - 1 = 8 \cdot \frac{2^{27} + 1}{27} - 1 = 39\,768\,215$$

The sequence generated by this number reaches "1" in one single cycle. The cycle follows an upward trajectory: one upward step (3x+1) and one downward step (divided by 2) alternate 3 times.

**39768215**
119304646
59652323
178956970
89478485
268435456
134217728

Once the upper bound is reach, the downward trajectory begins. There are 27 consecutive downward steps (divided by 2) until the sequence reaches 1.

67108864
33554432
16777216
8388608
4194304
2097152
1048576
524288
262144
131072
65536
32768
16384

8192
4096
2048
1024
512
256
128
64
32
16
8
4
2
1

The total number of downward steps – upward and downward trajectories – is 30:

3 + 27 = 30.

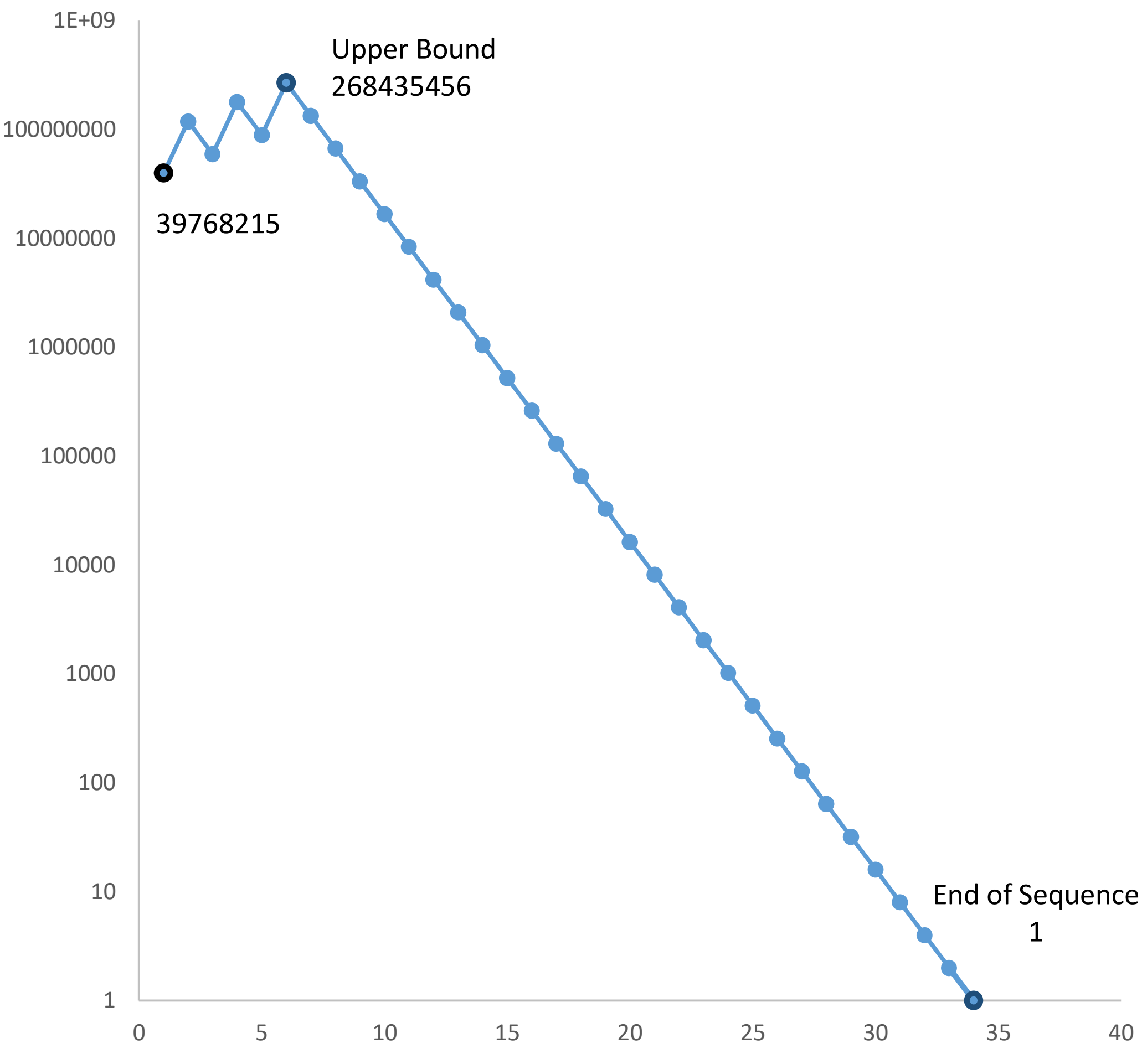

# Appendix 4 – Evenly divisible by $3^{n+1}$

In this appendix, we prove that the following expressions:

$$2^{3^n \cdot K_o} + 1 \ \text{ and } \ 2^{3^n \cdot K_E} - 1$$

are evenly divisible by $3^{n+1}$

In other words:

$$\left\{\frac{2^{3^n \cdot K_o} + 1}{3^{n+1}}\right\}; \left\{\frac{2^{3^n \cdot K_E} - 1}{3^{n+1}}\right\} \in \mathbb{N}$$

$$\forall\, n \in \mathbb{N}$$
$$\forall\, K_o \in \{ODD\} \in \mathbb{N}$$
$$\forall\, K_E \in \{EVEN\} \in \mathbb{N}$$

## Lemma A4.1

Let

$$n \in \mathbb{N}$$
$$K_o; i_o = ODD \in \mathbb{N}$$
$$i_E = EVEN \in \mathbb{N}$$

The following expression

$$(A4.1)\ \ \frac{2^{2\cdot 3^n \cdot K_o} - 2^{3^n \cdot K_o} + 1}{3} = \sum_{i_o = 3^n \cdot K_o}^{i_o = 2\cdot 3^n \cdot K_o - 1} 2^{i_o} - \sum_{i_o = 1}^{i_o = 3^n \cdot K_o} 2^{i_o} - \sum_{i_E = 3^n \cdot K_o + 1}^{i_E = 2\cdot 3^n \cdot K_o - 2} 2^{i_E} + \sum_{i_E = 0}^{i_E = 3^n \cdot K_o - 1} 2^{i_E}$$

is true. Therefore, the following expression

$$2^{2\cdot 3^n \cdot K} - 2^{3^n \cdot K} + 1$$

is evenly divisible by 3

In other words:

$$\left\{\frac{2^{2\cdot 3^n \cdot K} - 2^{3^n \cdot K} + 1}{3}\right\} \in$$

**Note:**

The denotation $i_E$ means that only the EVEN values are being used in the summation. Accordingly, the denotation $i_O$ means that only the ODD values are being used in the summation.

Example:

$$\sum_{i_E=2}^{i_E=10} 2^{i_E} = 2^2 + 2^4 + 2^6 + 2^8 + 2^{10}$$

## Proof

We want to prove that the identity below is true

$$\frac{2^{2\cdot 3^n\cdot K_o} - 2^{3^n\cdot K_o} + 1}{3} = \sum_{i_o=3^n\cdot K_o}^{i_o=2\cdot 3^n\cdot K_o-1} 2^{i_o} - \sum_{i_o=1}^{i_o=3^n\cdot K_o} 2^{i_o} - \sum_{i_E=3^n\cdot K_o+1}^{i_E=2\cdot 3^n\cdot K_o-2} 2^{i_E} + \sum_{i_E=0}^{i_E=3^n\cdot K_o-1} 2^{i_E}$$

We decompose "3" into the binomial "2 +1", and pass the denominator to the other side

$$(A4.2) \qquad 2^{2\cdot 3^n\cdot K_o} - 2^{3^n\cdot K_o} + 1 =$$

$$= (2+1)\cdot\left[\sum_{i_o=3^n\cdot K_o}^{i_o=2\cdot 3^n\cdot K_o-1} 2^{i_o} - \sum_{i_o=1}^{i_o=3^n\cdot K_o} 2^{i_o} - \sum_{i_E=3^n\cdot K_o+1}^{i_E=2\cdot 3^n\cdot K_o-2} 2^{i_E} + \sum_{i_E=0}^{i_E=3^n\cdot K_o-1} 2^{i_E}\right]$$

Therefore, if (A4.1) is true, (A4.2) must also be true.

We first multiply the first member of the binomial "2" by the 4 summations within the brackets, and then we do the same for the second member of the binomial "1".

$$2 \text{ x } [\ldots] = \sum_{i_E=3^n\cdot K_o+1}^{i_E=2\cdot 3^n\cdot K_o} 2^{i_E} - \sum_{i_E=2}^{i_E=3^n\cdot K_o+1} 2^{i_E} - \sum_{i_o=3^n\cdot K_o+2}^{i_o=2\cdot 3^n\cdot K_o-1} 2^{i} + \sum_{i_o=1}^{i_o=3^n\cdot K_o} 2^{i} +$$

$$1 \text{ x } [\ldots] = + \sum_{i_o=3^n\cdot K_o}^{i_o=2\cdot 3^n\cdot K_o-1} 2^{i_o} - \sum_{i_o=1}^{i_o=3^n\cdot K_o} 2^{i} - \sum_{i_E=3^n\cdot K_o+1}^{i_E=2\cdot 3^n\cdot K_o-2} 2^{i_E} + \sum_{i_E=0}^{i_E=3^n\cdot K_o-1} 2^{i_E}$$

We simplify the above. Note that almost all members of the summations annul each other.

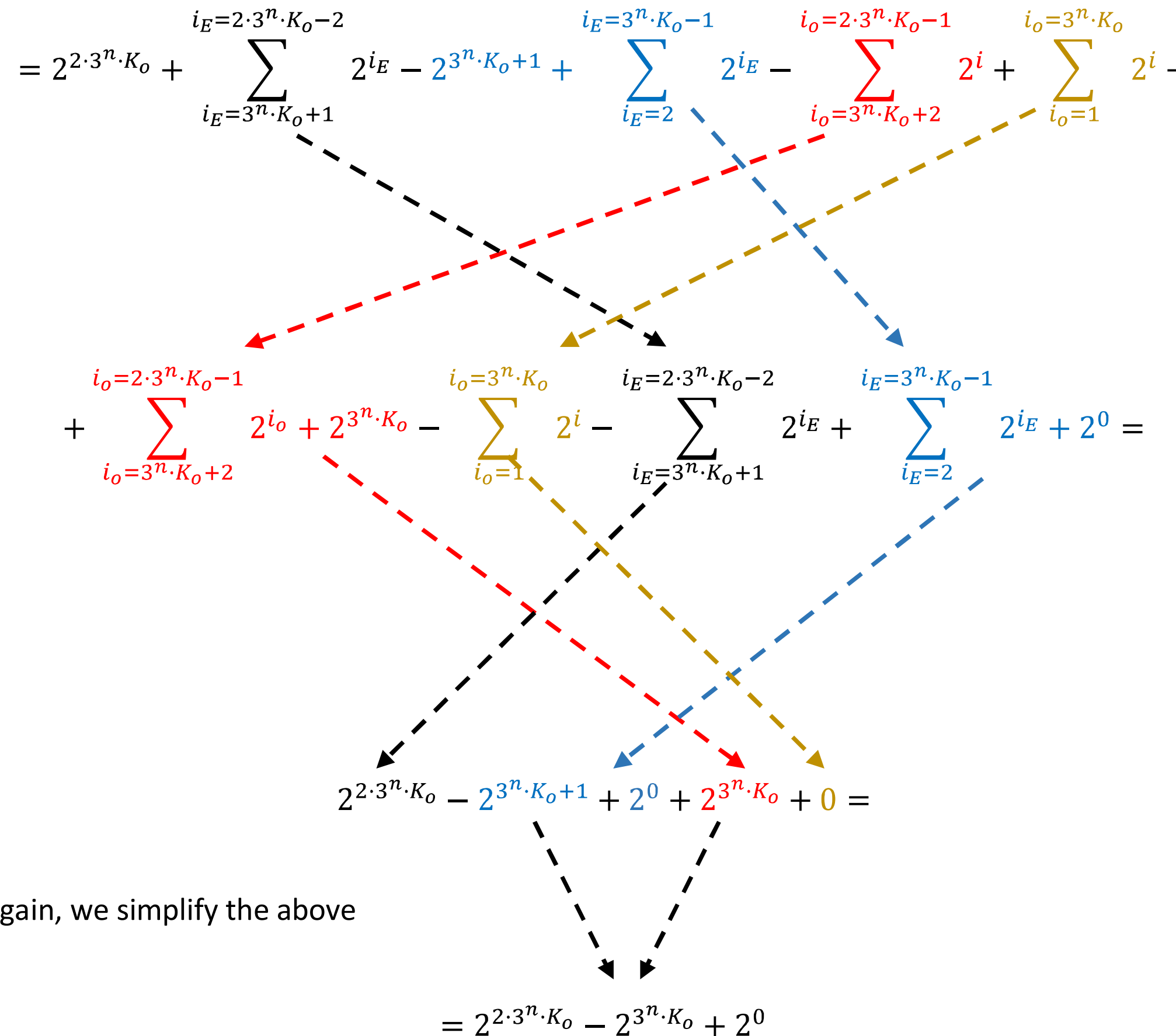

which is what we wanted to prove.

Therefore, the above number is evenly divisible by 3.

Lemma A4.2

The following expression

$$2^{K_o} + 1$$

is evenly divisible by 3

In other words:

$$\left\{\frac{2^{K_o} + 1}{3}\right\} \in \mathbb{N}$$

$$\forall\ K_o = ODD \in \mathbb{N}$$

## Proof

If we prove that the identity below is true, then Lemma A4.2 is true

$$\frac{2^{K_o} + 1}{3} = \frac{2^{K_o} + 1}{2 + 1} = \sum_{i_E=0}^{i_E=K_o-1} 2^{i_E} - \sum_{i_o=1}^{i_o=K_o-2} 2^{i_o}$$

we pass the binomial “2+1” to the other side

$$2^{K_o} + 1 = (2 + 1) \cdot [\sum_{i_E=0}^{i_E=K_o-1} 2^{i_E} - \sum_{i_o=1}^{i_o=K_o-2} 2^{i_o}]$$

We first multiply the first member of the binomial “2” by the 2 summations within the brackets, and then we do the same for the second member of the binomial “1”.

2 x [...] = $\sum_{i_O=1}^{i_O=K_o} 2^{i_O} - \sum_{i_E=2}^{i_E=K_o-1} 2^{i_E} +$

1 x [...] = $+ \sum_{i_E=0}^{i_E=K_o-1} 2^{i_E} - \sum_{i_o=1}^{i_o=K_o-2} 2^{i_o}$

We simplify the above. Note that almost all members of the summations annul each other.

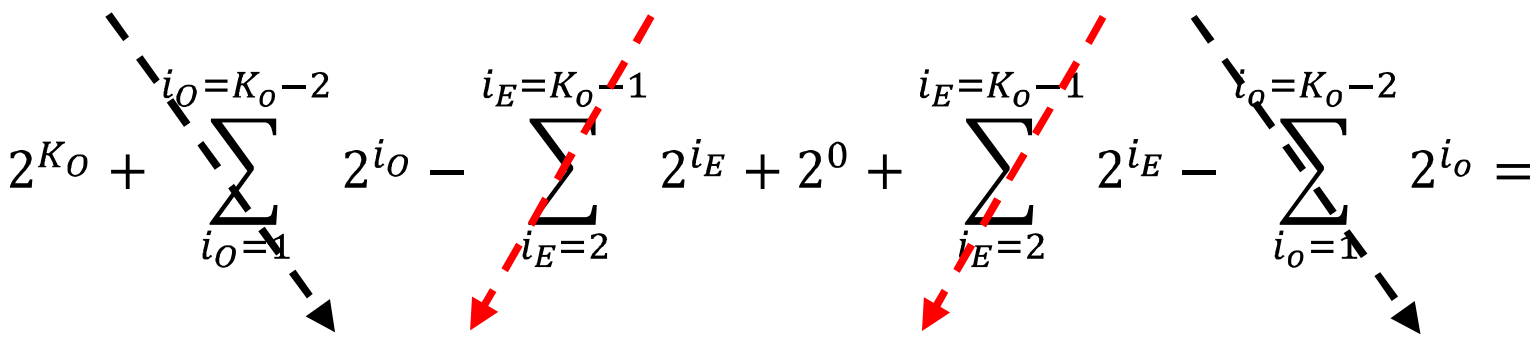

$$2^{K_O} + \sum_{i_O=1}^{i_O=K_O-2} 2^{i_O} - \sum_{i_E=2}^{i_E=K_O-1} 2^{i_E} + 2^0 + \sum_{i_E=2}^{i_E=K_O-1} 2^{i_E} - \sum_{i_o=1}^{i_o=K_o-2} 2^{i_o} =$$

$$2^{K_O} + 2^0 = 2^{K_O} + 1$$

which is what we wanted to prove.

**Lemma A4.3**

The following expression

$$2^{3^n \cdot K_o} + 1$$

is evenly divisible by $3^{n+1}$

In other words:

$$\left\{ \frac{2^{3^n \cdot K_o} + 1}{3^{n+1}} \right\} \in \mathbb{N}$$

$$\forall\, n \in \mathbb{N}$$

$$\forall\, K_o = ODD \in \mathbb{N}$$

For example: n=2; $K_O$=5

$$\frac{2^{5 \cdot 3^2} + 1}{3^3} = \frac{35{,}184{,}372{,}088{,}833}{27} = 1{,}303{,}124{,}892{,}179$$

## Proof

We take the binomial

$$(2^{3^n \cdot K_o} + 1)$$

This is the sum of two $n^{th}$ powers; therefore it can be factored in the following way:

$$2^{3^n \cdot K_o} + 1 = \left(2^{3^{n-1} \cdot K_o} + 1\right) \cdot (2^{2 \cdot 3^{n-1} \cdot K_o} - 2^{3^{n-1} \cdot K_o} + 1)$$

we repeat the same reasoning with the binomial $2^{3^{n-1} \cdot K} + 1$ and so on.

$$2^{3^{n-1} \cdot K_o} + 1 = \left(2^{3^{n-2} \cdot K_o} + 1\right) \cdot \left(2^{2 \cdot 3^{n-2} \cdot K_o} - 2^{3^{n-2} \cdot K_o} + 1\right)$$

$$2^{3^{n-2} \cdot K_o} + 1 = \left(2^{3^{n-3} \cdot K_o} + 1\right) \cdot \left(2^{2 \cdot 3^{n-3} \cdot K_o} - 2^{3^{n-3} \cdot K_o} + 1\right)$$

$$\cdots$$

$$\cdots$$

$$2^{3^2 \cdot K_o} + 1 = \left(2^{3^1 \cdot K_o} + 1\right) \cdot (2^{2 \cdot 3^1 \cdot K_o} - 2^{3^1 \cdot K_o} + 1)$$

$$2^{3^1 \cdot K_o} + 1 = \left(2^{3^0 \cdot K_o} + 1\right) \cdot \left(2^{2 \cdot 3^0 \cdot K_o} - 2^{3^0 \cdot K_o} + 1\right)$$

Therefore, we can rearrange the expression

$$2^{3^n \cdot K_o} + 1$$

as follows

$$(A4.3) \qquad 2^{3^n \cdot K_o} + 1 = (2^{K_o} + 1) \cdot \prod_{j=0}^{j=n-1} (2^{2\cdot 3^j \cdot K_o} - 2^{3^j \cdot K_o} + 1)$$

from Lemma A4.1 we know that

$$2^{2\cdot 3^j \cdot K_o} - 2^{3^j \cdot K_o} + 1$$

is evenly divisible by 3. Therefore (A4.3) can be rearranged as

$$(A4.4) \qquad 2^{3^n \cdot K_o} + 1 = (2^{K_o} + 1) \cdot \prod_{j=0}^{j=n-1} 3 \cdot \frac{(2^{2\cdot 3^j \cdot K_o} - 2^{3^j \cdot K_o} + 1)}{3}$$

from Lemma A4.2, we know that

$$2^{K_o} + 1$$

is evenly divisible by 3. Therefore (A4.4) can be rearranged as

$$2^{3^n \cdot K_o} + 1 = 3 \cdot \frac{2^{K_o} + 1}{3} \cdot \prod_{j=0}^{j=n-1} 3 \cdot \frac{(2^{2\cdot 3^j \cdot K_o} - 2^{3^j \cdot K_o} + 1)}{3}$$

we extract 3 from the factors of the product, therefore

$$(A4.5) \qquad 2^{3^n \cdot K_o} + 1 = 3^{n+1} \cdot \frac{2^{K_o} + 1}{3} \cdot \prod_{j=0}^{j=n-1} \frac{(2^{2\cdot 3^j \cdot K_o} - 2^{3^j \cdot K_o} + 1)}{3}$$

which is what we wanted to prove.

**Lemma A4.4**

The following expression

$$2^{3^n \cdot K_E} - 1$$

is evenly divisible by $3^{n+1}$

In other words:

$$\left\{\frac{2^{3^n \cdot K_E} - 1}{3^{n+1}}\right\} \in \mathbb{N}$$

$$\forall\, n \in \mathbb{N}$$

$$\forall\, K_E \in \{EVEN\}$$

## Proof

Let

$$K_E = 2^\alpha \cdot K_O$$

$$K_O = odd \in \mathbb{N}$$

then; we can rewrite

$$2^{3^n \cdot K_E} - 1 = 2^{3^n \cdot (2^\alpha \cdot K_O)} - 1$$

since the above is the difference of two squares

$$2^{3^n \cdot K_E} - 1 = 2^{3^n \cdot (2^\alpha \cdot K_O)} - 1 = \left(2^{3^n \cdot (2^{\alpha-1} \cdot K_O)} - 1\right) \cdot \left(2^{3^n \cdot 2^{\alpha-1} \cdot K_O} + 1\right)$$

we find again a difference of two squares. Then, we repeat the process until there is no more difference of two squares:

$$2^{3^n \cdot K_E} - 1 = \left(2^{3^n \cdot K_O} - 1\right) \cdot \left(2^{3^n \cdot K_O} + 1\right) \prod_{i=2}^{i=\alpha-1} \left(2^{3^n \cdot 2^i \cdot K_O} + 1\right)$$

We apply Lemma A4.3 to the binomial

$$2^{3^n \cdot K_O} + 1$$

therefore, the above is evenly divisible by $3^{n+1}$. We rearrange the above

$$2^{3^n \cdot K_E} - 1 = 3^{n+1} \cdot \frac{2^{3^n \cdot K_O} + 1}{3^{n+1}} \cdot \left(2^{3^n \cdot K_O} - 1\right) \cdot \prod_{2}^{\alpha-1} \left(2^{3^n \cdot 2^i \cdot K_O} + 1\right)$$

which is what we wanted to prove.

## Lemma A4.5

The following expression

$$\frac{2^{3^{n-1}\cdot K_o}+1}{3^n}$$

only becomes a pure power of 3

$$\frac{2^{3^{n-1}\cdot K_o}+1}{3^n}=3^t$$

if and only if

$$K_o=3$$

and

$$n=1$$

## Proof

We use identity (A4.5)

$$(A4.6)\qquad \frac{2^{3^{n-1}\cdot K_o}+1}{3^n}=\frac{2^{K_o}+1}{3}\cdot\prod_{j=0}^{j=n-2}\frac{\left(2^{2\cdot 3^j\cdot K_o}-2^{3^j\cdot K_o}+1\right)}{3}$$

Since the element

$$\frac{2^{3^{n-1}\cdot K_o}+1}{3^n}$$

is a product of factors, the expression above is a power of 3 if only and only if all the factors are the number 3.

The first factor

$$\frac{2^{K_o}+1}{3}=3$$

if and only if

$$K_o=3$$

The following factors are

$$\frac{2^{2\cdot 3^j\cdot K_o}-2^{3^j\cdot K_o}+1}{3}$$

Since we already know that $K_o = 3$, the expression above is

$$\frac{2^{2\cdot 3^{j}\cdot 3} - 2^{3^{j}\cdot 3} + 1}{3} \neq 3$$

$$\forall\, j$$

In reality is sufficient to prove that the expression above is NOT 3 for:

$$j = 0$$

$$\frac{2^{2\cdot 3^{j}\cdot 3} - 2^{3^{j}\cdot 3} + 1}{3} = \frac{2^{2\cdot 3^{0}\cdot 3} - 2^{3^{0}\cdot 3} + 1}{3} = 19 \neq 3$$

This is the case since this factor is the first element of multiplication, it appears for any

$$n \geq 2$$

As a result, if and only if

$$n = 1$$

is value that makes the identity (A4.6) a pure power of 3

$$K_o = 3$$

and

$$n = 1$$

$$\frac{2^{3^{n-1}\cdot K_o} + 1}{3^{n}} = \frac{2^{3^{0}\cdot 3} + 1}{3^{1}} = 3$$

Which is what we wanted to prove.

# Appendix 5 – Odd Numbers

In this appendix, we will prove that all odd numbers can be written using the following identity:

$$(A5.1) \quad a_O = 3^m \cdot K_O - 2^n \cdot \frac{2^{3^{n-1} \cdot q_o} + 1}{3^n}$$

$$\forall\, n, m \in \mathbb{N}$$

$$K_O \in \{Odd\}$$

$$q_o \in \{1\,, 3\,, 5, \dots, 2 \cdot 3^m - 1\}$$

## Proof

The set of all odd numbers could be divided into three subsets

| | 1 | 3 | 5 | 7 | … | … | $K_O$ |
|---|---|---|---|---|---|---|---|
| **S1** | 1 | 7 | 13 | 19 | … | … | $3 \cdot K_O - 2$ |
| **S2** | - | 3 | 9 | 15 | 21 | … | $3 \cdot K_O - 6$ |
| **S3** | - | 5 | 11 | 17 | 23 | … | $3 \cdot K_O - 4$ |

$$K_O \in \{Odd\}$$

In other words, the set of all odd numbers can be divided into three subsets as follows:

$$\forall\ a_O \in \{Odd\}$$

$$(A7.2)\ \ a_O \in\ \{3 \cdot K_O - K_E\}$$

$$K_E \in \{2;\ 4;\ 6\}$$

$$K_O \in \{Odd\}$$

and the union of the three subsets represents the set of all the odd numbers:

$$\bigcup_{i=1}^{i=3} S_i = \{Odd\}$$

And, obviously, the three subsets are disjoint.

$$S_i \cap S_j = \phi$$

$$\forall\ \ i\,;\ j \in \{1{,}2{,}3\}$$

Also, all odd numbers could be divided into nine subsets:

| **S1** | 1 | 19 | 37 | ... | $3^2 \cdot K_O - 8$ |
|---|---|---|---|---|---|
| **S2** | 3 | 21 | 39 | ... | $3^2 \cdot K_O - 6$ |
| **S3** | 5 | 23 | 41 | ... | $3^2 \cdot K_O - 4$ |
| **S4** | 7 | 25 | 43 | ... | $3^2 \cdot K_O - 2$ |
| **S5** | 9 | 27 | 45 | ... | $3^2 \cdot K_O - 18$ |
| **S6** | 11 | 29 | 47 | ... | $3^2 \cdot K_O - 16$ |
| **S7** | 13 | 31 | 49 | ... | $3^2 \cdot K_O - 14$ |
| **S8** | 15 | 33 | 51 | ... | $3^2 \cdot K_O - 12$ |
| **S9** | 17 | 35 | 53 | ... | $3^2 \cdot K_O - 10$ |

In other words, the set of all odd numbers can be divided in nine subsets as follows:

$$\forall \; a_O \in \{Odd\}$$

$$(A7.3) \;\; a_O \in \left\{3^2 \cdot K_O - K_E\right\}$$

$$K_E \in \{2; 4; 6; 8; 10; 12; 14; 16; 18\}$$

$$K_O \in \{Odd\}$$

and the union of the nine subsets represents the set of all the odd numbers:

$$\bigcup_{i=1}^{i=3^2} S_i = \{Odd\}$$

And, obviously, the nine subsets are disjoint.

$$S_i \cap S_j = \phi$$

$$\forall \; i \; ; \; j \in \{1,2,3,4,5,6,7,8,9\}$$

We could also divide all odd numbers into twenty-seven subsets:

| **S1** | 1 | 55 | 109 | ... | $3^3 \cdot K_O - 26$ |
|---|---|---|---|---|---|
| **S2** | 3 | 57 | 111 | ... | $3^3 \cdot K_O - 24$ |
| **S3** | 5 | 59 | 113 | ... | $3^3 \cdot K_O - 22$ |
| **S4** | 7 | 61 | 115 | ... | $3^3 \cdot K_O - 20$ |
| **S5** | 9 | 63 | 117 | ... | $3^3 \cdot K_O - 18$ |
| **S6** | 11 | 65 | 119 | ... | $3^3 \cdot K_O - 16$ |
| **S7** | 13 | 67 | 121 | ... | $3^3 \cdot K_O - 14$ |
| **S8** | 15 | 69 | 123 | ... | $3^3 \cdot K_O - 12$ |
| **S9** | 17 | 71 | 125 | ... | $3^3 \cdot K_O - 10$ |
| **S10** | 19 | 73 | 127 | ... | $3^3 \cdot K_O - 8$ |
| **S11** | 21 | 75 | 129 | ... | $3^3 \cdot K_O - 6$ |
| **S12** | 23 | 77 | 131 | ... | $3^3 \cdot K_O - 4$ |
| **S13** | 25 | 79 | 133 | ... | $3^3 \cdot K_O - 2$ |
| **S14** | 27 | 81 | 135 | ... | $3^3 \cdot K_O - 54$ |
| **S15** | 29 | 83 | 137 | ... | $3^3 \cdot K_O - 52$ |
| **S16** | 31 | 85 | 139 | ... | $3^3 \cdot K_O - 50$ |
| **S17** | 33 | 87 | 141 | ... | $3^3 \cdot K_O - 48$ |
| **S18** | 35 | 89 | 143 | ... | $3^3 \cdot K_O - 46$ |
| **S19** | 37 | 91 | 145 | ... | $3^3 \cdot K_O - 44$ |
| **S20** | 39 | 93 | 147 | ... | $3^3 \cdot K_O - 42$ |
| **S21** | 41 | 95 | 149 | ... | $3^3 \cdot K_O - 40$ |
| **S22** | 43 | 97 | 151 | ... | $3^3 \cdot K_O - 38$ |
| **S23** | 45 | 99 | 153 | ... | $3^3 \cdot K_O - 36$ |
| **S24** | 47 | 101 | 155 | ... | $3^3 \cdot K_O - 34$ |
| **S25** | 49 | 103 | 157 | ... | $3^3 \cdot K_O - 32$ |
| **S26** | 51 | 105 | 159 | ... | $3^3 \cdot K_O - 30$ |
| **S27** | 53 | 107 | 161 | ... | $3^3 \cdot K_O - 28$ |

In other words, the set of all odd numbers can be divided into twenty-seven subsets as follows:

$$\forall \ a_O \in \{Odd\}$$

$$(A7.4) \ \ a_O \in \left\{ 3^3 \cdot K_O - K_E \right\}$$

$$K_E \in \{2; 4; 6; 8; 10; \ldots . ; 54\}$$

$$K_O \in \{Odd\}$$

and the union of the twenty-seven subsets represents the set of all the odd numbers:

$$\bigcup_{i=1}^{i=3^3} S_i = \{Odd\}$$

And, obviously, the twenty-seven subsets are disjoint.

$$S_i \cap S_j = \phi$$

$$\forall \ i \ ; \ j \in \{1,2,3,4,5, \ldots ,26,27\}$$

Finally, if we follow this process "m" times, we could also divide all odd numbers in $3^m$ subsets. Therefore:

## Lemma A5.1

The set of all odd numbers can be divided in $3^m$ subsets as follows:

$$\forall \ a_O \in \{Odd\}$$

$$a_O \in S_i = \left\{3^m \cdot K_O - K_{E.i}\right\}$$
$$K_{E.i} \in \{2; 4; 6; 8; 10; \ldots.; 2 \cdot 3^m\}$$
$$K_O \in \{Odd\}$$
$$i \in \{1,2,3,\ldots,3^m\}$$
$$\forall \, m \in \mathbb{N}$$

and the union of all the subsets represents the set of all the odd numbers:

$$\bigcup_{i=1}^{i=3^m} S_i = \{Odd\}$$

And, obviously, the $3^m$ subsets are disjoint.

$$S_i \cap S_j = \phi$$
$$\forall \ i\,;\, j \in \{1,2,3,4,5,\ldots,3^m\}$$

## Lemma A5.2

Let the subset $S_1$ and subset $S_2$ of odd numbers be:

$$(A7.5) \qquad S_1 = \{3^m \cdot K_O^1 - K_E^1\}$$

$$(A7.6) \qquad S_2 = \{3^m \cdot K_O^2 - K_E^2\}$$

$$K_O^1; K_O^2 \in \{odd\}$$

$$K_E^1; K_E^2 \in \{even\}$$

$$m \in \mathbb{N}$$

Let $K_E^1$ and $K_E^2$ be even numbers such as:

$$(A7.7) \qquad K_E^2 = K_E^1 + k \cdot 2 \cdot 3^m$$

$$k \in \mathbb{N}$$

Both subsets represent the same subset of odd numbers

$$S_1 \cap S_2 = S_1 = S_2$$

## Proof

We use any odd number from subset $S_2$ using Identity (A5.6)

$$a_O^2 = 3^m \cdot K_O^2 - K_E^2$$

We apply Identity (A5.7) to the above

$$a_O^2 = 3^m \cdot K_O^2 - (K_E^1 + k \cdot 2 \cdot 3^m)$$

We combine all the elements with the factor $3^m$

$$a_O^2 = 3^m \cdot (K_O^2 - k \cdot 2) - K_E^1$$

We make the following change

$$K_O^1 = K_O^2 - k \cdot 2$$

and obtain

$$a_O^2 = 3^m \cdot K_O^1 - K_E^1$$

Therefore, $a_O^2$ is also an element of the subset $S_1$. As a result

$$S_1 \cap S_2 = S_1 = S_2$$

which is what we wanted to prove.

## Lemma A5.3

Let the subset $S_1$ and subset $S_2$ of odd numbers be:

$$(A7.8) \qquad S_1 = \{3^m \cdot K_O^1 - K_E^1\}$$

$$(A7.9) \qquad S_2 = \{3^m \cdot K_O^2 - K_E^2\}$$

$$K_O^1; K_O^2 \in \{odd\}$$

$$K_E^1; K_E^2 \in \{even\}$$

$$m \in \mathbb{N}$$

Let $K_E^1$ and $K_E^2$ be even numbers such as:

$$(A7.10) \qquad K_E^2 - K_E^1 \neq k \cdot 2 \cdot 3^m$$

$$k \in \mathbb{N}$$

Then, no number can belong to both subsets at the same time. In other words:

$$(A7.11) \qquad S_1 \cap S_2 = \phi$$

**Proof**

Let us assume the opposite is true, in other words, there is an element belonging to each subset that represent the same odd number:

$$a_O^1 \in S_1 \quad ; \quad a_O^2 \in S_2$$

$$a_O^1 = a_O^2$$

If this were the case, we could combine Identity (A5.8) and Identity (A5.9) and obtain

$$3^m \cdot K_O^1 - K_E^1 = 3^m \cdot K_O^2 - K_E^2$$

We put together all elements with the factor $3^m$ on one side of the identity:

$$(A7.12) \qquad 3^m \cdot (K_O^1 - K_O^2) = K_E^1 - K_E^2$$

Which cannot be. Inequality (A5.10) states that

$$K_E^2 - K_E^1 \neq k \cdot 2 \cdot 3^m$$

and, on the other hand, the element of left side of the Identity (A5.12) is evenly divisible by $3^m$. Therefore, there is a contradiction, as a result:

$$a_O^1 \neq a_O^2$$

and

$$S_1 \cap S_2 = \phi$$

Which is what we wanted to prove.

## Lemma A5.4

The set of all odd numbers can be divided in $3^{\mathrm{m}}$ subsets as follows

$$S'_i \in \left\{3^m \cdot K_O - Q_{E.i}\right\}$$

$$Q_{E.i} \in \{Q_{E.1}; Q_{E.2}; Q_{E.3}; \ldots.; Q_{E.3^m}\} \in \{even\}$$
$$K_O \in \{Odd\}$$
$$Q_{E..i} - Q_{E.j} \neq k \cdot 2 \cdot 3^m$$
$$\forall\ i; j \in \{1,2,3, \ldots., 3^m\}$$
$$\forall\, m \in \mathbb{N}$$

and the union of all the subsets represents the set of all the odd numbers:

$$\bigcup_{i=1}^{i=3^m} S'_i = \{Odd\}$$

and the $3^m$ subsets are disjoint.

$$S'_i \cap S'_j = \phi$$
$$\forall\ i\ ;\ j \in \{1,2,3,4,5, \ldots, 3^m\}$$

## Proof

Lemma A5.1 states that the set of all odd numbers can be divided in $3^m$ subsets as follows:

$$S_i \in \left\{3^m \cdot K_O - K_{E.i}\right\}$$
$$K_{E.i} \in \{2; 4; 6; 8; 10; \ldots.; 2 \cdot 3^m\}$$
$$K_O \in \{Odd\}$$

and the union of all the subsets represents the set of all the odd numbers

$$\bigcup_{i=1}^{i=3^m} S_i = \{Odd\}$$

and the $3^m$ subsets are disjoint

$$S_i \cap S_j = \phi$$
$$\forall\ i\ ;\ j \in \{1,2,3,4,5, \ldots, 3^m\}$$

Lemma A5.2 states, for each subset above, we can substitute the value $K_{E.i}$ for another even number $Q_{E.i}$ such as

$$(A7.13) \quad Q_{E.i} = K_{E.i} + k \cdot 2 \cdot 3^m$$

and the subset we obtain is exactly the same as the previous one. In other words

$$S_i \in \{3^m \cdot K_O - K_{E.i}\}$$
$$S'_i \in \{3^m \cdot K_O - Q_{E.i}\}$$

and

$$S_i \cap S'_i = S_i = S'_i$$

Therefore, the set of all numbers can be divided in $3^m$ subsets as follows:

$$\forall \ a_O \in \{Odd\}$$

$$S'_i \in \{3^m \cdot K_O - Q_{E.i}\}$$
$$Q_{E.i} \in \{Q_{E.1}; Q_{E.2}; Q_{E.3}; \ldots.; Q_{E.3^m}\} \in \{even\}$$
$$K_O \in \{Odd\}$$

and the union of all the subsets represents the set of all the odd numbers

$$\bigcup_{i=1}^{i=3^m} S'_i = \{Odd\}$$

and the $3^m$ subsets are disjoint

$$S'_i \cap S'_j = \phi$$
$$\forall \ i \ ; \ j \in \{1,2,3,4,5, \ldots, 3^m\}$$

Logically, since all the numbers within the subset $K_{E.i}$

$$K_{E.i} \in \{2; 4; 6; 8; 10; \ldots.; 2 \cdot 3^m\}$$

are smaller than $2 \cdot 3^m$, (except for the last one that it is equal to $2 \cdot 3^m$)

$$K_{E.i} < 2 \cdot 3^m$$

then, we can state that the difference between two values $K_{E.i}$ is always different than $k \cdot 2 \cdot 3^m$

$$(A7 - 14) \quad K_{E.i} - K_{E.j} \neq k \cdot 2 \cdot 3^m$$
$$\forall \, i, j \ \in \{1,2,3, \ldots, 3^m\}$$

As a result

$$Q_{E.i} - Q_{E.j}$$

also complies with

$$Q_{E.i} - Q_{E.j} \neq k \cdot 2 \cdot 3^m$$

This is the case since Identity (A5.13) states

$$Q_{E.i} = K_{E.i} + k \cdot 2 \cdot 3^m$$

then

$$Q_{E.i} - Q_{E.j} = K_{E.i} + k_i \cdot 2 \cdot 3^m - K_{E.j} - k_j \cdot 2 \cdot 3^m$$

We reorganized the above putting the elements with the factor $3^m$ together and obtain

$$(A7.15) \quad Q_{E.i} - Q_{E.j} = K_{E.i} - K_{E.j} + \left(k_i - k_j\right) \cdot 2 \cdot 3^m$$

Since Inequality (A5.14) above states that

$$K_{E.i} - K_{E.j} \neq k \cdot 2 \cdot 3^m$$

then, we apply the above to identity (A5.15) and obtain

$$Q_{E.i} - Q_{E.j} \neq k \cdot 2 \cdot 3^m$$

Which is what we wanted to prove. This is logical since any two subsets above are disjoint.

$$S_i \cap S_j = \phi$$

as a result, $Q_{E.i}$ and $Q_{E.j}$ must comply with Lemma A5.3

## Lemma A5.5

Let the subsets of odd numbers be as follows

$$S_i \in \left\{ 3^m \cdot K_O - 2^n \cdot \frac{2^{3^{n-1} \cdot q_{o.i}} + 1}{3^n} \right\}$$

$$q_{o.i} \in \{1; 3; 5; 7; 11; \ldots.; 2 \cdot 3^m - 1\}$$

$$K_O \in \{Odd\}$$

$$i \in \{1; 2; 3; \ldots.; 3^m\}$$

These subsets are disjoint. In other words, the same odd number cannot belong to two different subsets with different value $q_O^i$

$$S_i \cap S_j = \emptyset$$

$$\forall\, i \;;\; j \in \{1; 2; 3; \ldots.; 3^m\}$$

## Proof

Let $a_O^1$ and $a_O^2$ two odd numbers from subsets $S_1$ and $S_2$ be:

$$a_O^1 \in \left\{ 3^m \cdot K_O^1 - 2^n \cdot \frac{2^{3^{n-1} \cdot q_{o.1}} + 1}{3^n} \right\}$$

and

$$a_O^2 \in \left\{ 3^m \cdot K_O^2 - 2^n \cdot \frac{2^{3^{n-1} \cdot q_{o.2}} + 1}{3^n} \right\}$$

Let us assume that both numbers are the same:

$$a_O^1 = a_O^2$$

Therefore

$$3^m \cdot K_O^1 - 2^n \cdot \frac{2^{3^{n-1} \cdot q_{o.1}} + 1}{3^n} = 3^m \cdot K_O^2 - 2^n \cdot \frac{2^{3^{n-1} \cdot q_{o.2}} + 1}{3^n}$$

We put all elements with the factor $3^m$ on the left side of the identity:

$$3^m \cdot (K_O^1 - K_O^2) = 2^n \cdot \frac{2^{3^{n-1} \cdot q_{o.1}} + 1}{3^n} - 2^n \cdot \frac{2^{3^{n-1} \cdot q_{o.2}} + 1}{3^n}$$

We simplify the right side of the identity above:

$$(A7.16) \quad 3^m \cdot (K_O^1 - K_O^2) = 2^n \cdot \frac{2^{3^{n-1} \cdot q_{o.1}} - 2^{3^{n-1} \cdot q_{o.2}}}{3^n}$$

Let us assume that

$$q_{o.1} > q_{o.2}$$

Please note it would make no difference if it were the other way around. Therefore, since

$$q_{o.1} > q_{o.2}$$

then we can extract

$$2^{3^{n-1} \cdot q_{o.2}}$$

from the numerator of the Identity (A5.16) above, and obtain:

$$(A7.17) \quad 3^m \cdot (K_O^1 - K_O^2) = 2^n \cdot 2^{3^{n-1} \cdot q_O^2} \cdot \frac{2^{3^{n-1} \cdot (q_{o.1} - q_{o.2})} - 1}{3^n}$$

Since

$$q_{o.1} \ ; \ q_{o.2} \in \{1; 3; 5; 7; 11; \ldots .; 2 \cdot 3^m - 1\}$$

The difference between the two numbers must be smaller than $2 \cdot 3^m$

$$q_{o.1} - q_{o.2} < 2 \cdot 3^m$$

Therefore, if we apply Lemma A4.4 from Appendix 4, the fraction on the right side of identity (A5.16)

$$\frac{2^{3^{n-1} \cdot (q_{o.1} - q_{o.2})} - 1}{3^n}$$

is **not** evenly divisible by $3^m$. In other words:

$$\frac{2^{3^{n-1} \cdot (q_{o.1} - q_{o.2})} - 1}{3^n} mod(3^m) \not\equiv 0$$

Which generates a contradiction. This is the case since the left side of Identity (A5.17) is evenly divisible by $3^m$ and the right side is not. As a result

$$a_O^1 \neq a_O^2$$

In other words, the same odd number cannot belong to two different subsets. Therefore, the two subsets $S_1$ and $S_2$ are disjoint

$$S_1 \cap S_2 = \emptyset$$

which is what we wanted to prove.

## Lemma A5.6

The set of all odd numbers can be divided in $3^{\mathrm{m}}$ subsets as follows

$$S_i \in \left\{ 3^m \cdot K_O - 2^n \cdot \frac{2^{3^{n-1} \cdot q_{o.i}} + 1}{3^n} \right\}$$

$$q_{o.i} \in \{1; 3; 5; 7; 11; \ldots.; 2 \cdot 3^m - 1\}$$
$$K_O \in \{Odd\}$$
$$i \in \{1,2,3,4,5, \ldots, 3^m\}$$
$$\forall\, n\,; m \in \mathbb{N}$$

The union of all the subsets represents the set of all the odd numbers:

$$\bigcup_{i=1}^{i=3^m} S_i = \{Odd\}$$

and the $3^m$ subsets are disjoint:

$$S_i \cap S_j = \phi$$
$$\forall\;\; i\,;\, j \in \{1,2,3,4,5, \ldots, 3^m\}$$

## Proof

Lemma A5.4 states that the union of $3^{\mathrm{m}}$ subsets that comply with the following:

$$S_i \in \{3^m \cdot K_O - K_{E.i}\}$$

$$K_{E.i} \in \{K_{E.1}; K_{E.2}; K_{E.3}; \ldots.; K_{E.3^m}\} \in \{even\}$$
$$K_O \in \{Odd\}$$
$$K_{E..i} - K_{E.j} \neq k \cdot 2 \cdot 3^m$$
$$\forall\;\; i; j \in \{1,2,3, \ldots., 3^m\}$$

represents the set of all the odd numbers

$$\bigcup_{i=1}^{i=3^m} S_i = \{Odd\}$$

Therefore, we generate a group of $3^m$ subsets of odd numbers choosing the different values $K_{E.i}$ as follows:

$$K_{E.i} = 2^n \cdot \frac{2^{3^{n-1} \cdot q_{o.i}} + 1}{3^n}$$
$$q_{o.i} \in \{1; 3; 5; 7; 11; \ldots.; 2 \cdot 3^m - 1\}$$

and obtain $3^m$ different subsets

$$S_i \in \left\{ 3^m \cdot K_O - 2^n \cdot \frac{2^{3^{n-1} \cdot q_{o.i}} + 1}{3^n} \right\}$$

$$q_{o.i} \in \{1; 3; 5; 7; 11; \ldots.; 2 \cdot 3^m - 1\}$$
$$K_O \in \{Odd\}$$
$$i \in \{1; 2; 3; \ldots.; 3^m\}$$
$$\forall\, n\, ; m \in \mathbb{N}$$

Lemma A5.5 above states that the subsets above are disjoint. In other words, the same odd number cannot belong to two different subsets with different value $q_{o.i}$

$$S_i \cap S_j = \emptyset$$
$$\forall\, i \;\; ; \;\; j \in \{1; 2; 3; \ldots.; 3^m\}$$

Therefore, the subsets above comply with Lemma A5.4. As a result, the union of these subsets represents the set of all the odd numbers

$$\bigcup_{i=1}^{i=3^m} S_i = \{Odd\}$$

Which is what we wanted to prove.

Please note that we arbitrarily chose the parameters “m” and “n”.

# Appendix 6 – Sequence of "i" Consecutive Cycles

Lemma A3.1 of Appendix 3 proves that any odd number generates a cycle. Every cycle follows the same pattern:

- First, there is an upward trajectory; an upward step (3x+1) and a downward step (divided by 2) alternate n times until they reach an upper bound.

- Once they reach this upper bound, then the downward trajectory begins. There are "$\alpha - (n-1)$" consecutive downward steps (divided by two) until the sequence reaches an odd number again.

Let $a_O$ be the initial element of a cycle

$$(A1.1) \qquad a_O = 2^{\alpha} \cdot Q_O - \frac{2^{3^{n-1}\cdot j_o^{\beta}} + 1}{3^n} \cdot (3^n - 2^n)$$

$$j_o^{\beta}, Q_O \in \{odd\}$$

$$\alpha, n \in \mathbb{N}$$

$$3^{n-1} \cdot j_o^{\beta} > \alpha > n$$

The final element of this cycle is:

$$(A3.2) \qquad a_F = 3^n \cdot Q_O - 2^{3^{n-1}\cdot j_o^{\beta} - \alpha} \cdot (3^n - 2^n)$$

In this Appendix, we combine "i" cycles together. The way we proceed is to modify the way the mathematical expression below is presented to obtain an identity with the same format as Identity (A1.1). Once we reach this point, we can state that that's the initial element of a cycle, and using Identity (A3.2), we can calculate the final element of this cycle.

## Lemma A6.1

Let $a_O^1$ - the initial element of a sequence made of "$i$" cycles - be

$$(A6.1) \quad a_O^1 = 2^{\sum_{t=1}^{t=i}\alpha_t} \cdot Q_{O.i} - \frac{2^{3^{\sum_{t=1}^{t=i} n_t - 1} \cdot j_{o.i}^{\delta}} + 1}{3^{\sum_{t=1}^{t=i} n_t}} \cdot \sum_{s=1}^{i} 3^{\sum_{t=s+1}^{t=i} n_t} \cdot 2^{\sum_{t=1}^{t=s-1}\alpha_t} \cdot (3^{n_s} - 2^{n_s})$$

where the above parameters $n_s, \alpha_s$ are the standard parameters of any Cycle "s" as described in Appendix 3

$$Q_{O.i}; j_{o.i}^{\delta} \in \{odd\}$$

$$\alpha_t > n_t \qquad \forall t$$

$$(A6.2) \quad 3^{\sum_{t=1}^{t=i} n_t - 1} \cdot j_{o.i}^{\delta} > \sum_{t=1}^{t=i} \alpha_t$$

The final element of this sequence after "$i$" cycles is

$$(A6.3) \quad a_F^i = 3^{\sum_{t=1}^{t=i} n_t} \cdot Q_{O.i} - 2^{3^{\sum_{t=1}^{t=i} n_t - 1} \cdot j_{o.i}^{\delta} - \sum_{t=1}^{t=i}\alpha_t} \cdot \sum_{s=1}^{i} 3^{\sum_{t=s+1}^{t=i} n_t} \cdot 2^{\sum_{t=1}^{t=s-1}\alpha_t} \cdot (3^{n_s} - 2^{n_s})$$

Please be aware that we are using the following convention:

$$2^{\sum_{t=1}^{t=0}\alpha_t} = 1$$

$$3^{\sum_{t=i+1}^{t=i} n_t} = 1$$

## Proof

We start with Identity (A6.1)

$$(A6.1) \quad a_O^1 = 2^{\sum_{t=1}^{t=i}\alpha_t} \cdot Q_{O.i} - \frac{2^{3^{\sum_{t=1}^{t=i} n_t - 1} \cdot j_{o.i}^{\delta}} + 1}{3^{\sum_{t=1}^{t=i} n_t}} \cdot \sum_{s=1}^{i} 3^{\sum_{t=s+1}^{t=i} n_t} \cdot 2^{\sum_{t=1}^{t=s-1}\alpha_t} \cdot (3^{n_s} - 2^{n_s})$$

We are going to modify the identity above to obtain an expression whose format coincides with Identity (A1.1). Once this is done, then we can say it represents the initial element of the first cycle and, using Identity (A3.2), we can then calculate the final element of this first cycle.

We extract the element "s=1" of the summation

$$a_O^1 = 2^{\sum_{t=1}^{t=i}\alpha_t}\cdot Q_{O.i} - \frac{2^{3^{\sum_{t=1}^{t=i}n_t-1}\cdot j_{o.i}^{\delta}}+1}{3^{\sum_{t=1}^{t=i}n_t}}\cdot\sum_{s=2}^{s=i}3^{\sum_{t=s+1}^{t=i}n_t}\cdot 2^{\sum_{t=1}^{t=s-1}\alpha_t}\cdot(3^{n_s}-2^{n_s})$$
$$-\frac{2^{3^{\sum_{t=1}^{t=i}n_t-1}\cdot j_{o.i}^{\delta}}+1}{3^{\sum_{t=1}^{t=i}n_t}}\cdot 3^{\sum_{t=2}^{t=i}n_t}\cdot(3^{n_1}-2^{n_1})$$

We simplify the last element. We eliminate the factor $3^{\sum_{t=2}^{t=i}n_t}$ that it is also in the denominator

$$a_O^1 = 2^{\sum_{t=1}^{t=i}\alpha_t}\cdot Q_{O.i} - \frac{2^{3^{\sum_{t=1}^{t=i}n_t-1}\cdot j_{o.i}^{\delta}}+1}{3^{\sum_{t=1}^{t=i}n_t}}\cdot\sum_{s=2}^{s=i}3^{\sum_{t=s+1}^{t=i}n_t}\cdot 2^{\sum_{t=1}^{t=s-1}\alpha_t}\cdot(3^{n_s}-2^{n_s})$$
$$-\frac{2^{3^{\sum_{t=1}^{t=i}n_t-1}\cdot j_{o.i}^{\delta}}+1}{3^{n_1}}\cdot(3^{n_1}-2^{n_1})$$

We extract the common factor $2^{\alpha_1}$ from all the remaining elements of the summation

$$a_O^1 = 2^{\sum_{t=1}^{t=i}\alpha_t}\cdot Q_{O.i} - 2^{\alpha_1}\cdot\frac{2^{3^{\sum_{t=1}^{t=i}n_t-1}\cdot j_{o.i}^{\delta}}+1}{3^{\sum_{t=1}^{t=i}n_t}}\cdot\sum_{s=2}^{s=i}3^{\sum_{t=s+1}^{t=i}n_t}\cdot 2^{\sum_{t=2}^{t=s-1}\alpha_t}\cdot(3^{n_s}-2^{n_s})$$
$$-\frac{2^{3^{\sum_{t=1}^{t=i}n_t-1}\cdot j_{o.i}^{\delta}}+1}{3^{n_1}}\cdot(3^{n_1}-2^{n_1})$$

We also extract the common factor $2^{\alpha_1}$from the first element of the expression

$$a_O^1 = 2^{\alpha_1}\cdot\left[2^{\sum_{t=2}^{t=i}\alpha_t}\cdot Q_{O.i} - \frac{2^{3^{\sum_{t=1}^{t=i}n_t-1}\cdot j_{o.i}^{\delta}}+1}{3^{\sum_{t=1}^{t=i}n_t}}\cdot\sum_{s=2}^{s=i}3^{\sum_{t=s+1}^{t=i}n_t}\cdot 2^{\sum_{t=2}^{t=s-1}\alpha_t}\cdot(3^{n_s}-2^{n_s})\right]$$
$$-\frac{2^{3^{\sum_{t=1}^{t=i}n_t-1}\cdot j_{o.i}^{\delta}}+1}{3^{n_1}}\cdot(3^{n_1}-2^{n_1})$$

We extract $3^{n_1-1}$ from the power in the numerator of the last fraction

$$(A6.4)\qquad a_O^1 = 2^{\alpha_1}\cdot\left[2^{\sum_{t=2}^{t=i}\alpha_t}\cdot Q_{O.i} - \frac{2^{3^{\sum_{t=1}^{t=i}n_t-1}\cdot j_{o.i}^{\delta}}+1}{3^{\sum_{t=1}^{t=i}n_t}}\cdot\sum_{s=2}^{s=i}3^{\sum_{t=s+1}^{t=i}n_t}\cdot 2^{\sum_{t=2}^{t=s-1}\alpha_t}\cdot(3^{n_s}-2^{n_s})\right]$$
$$-\frac{2^{3^{n_1-1}\cdot\left(3^{\sum_{t=2}^{t=i}n_t}\cdot j_{o.i}^{\delta}\right)}+1}{3^{n_1}}\cdot(3^{n_1}-2^{n_1})$$

Since $j_{o.i}^{\delta}$ complies with Inequality (A6.2)

$$(A6.2) \qquad 3^{\sum_{t=1}^{t=i} n_t - 1} \cdot j_{o.i}^{\delta} > \sum_{t=1}^{t=i} \alpha_t$$

It is also true

$$3^{\sum_{t=1}^{t=i} n_t - 1} \cdot j_{o.i}^{\delta} > \alpha_1$$

In other words, the Identity (A6.4) above follows the same format as Identity (A1.1) and represents the initial odd number of Cycle 1. Indeed, if we compare Identity (A1.1) and Identity (A6.4) we obtain:

$$Q_{O.1} = 2^{\sum_{t=2}^{t=i} \alpha_t} \cdot Q_{O.i} - \frac{2^{3^{\sum_{t=1}^{t=i} n_t - 1} \cdot j_{o.i}^{\delta}} + 1}{3^{\sum_{t=1}^{t=i} n_t}} \cdot \sum_{s=2}^{s=i} 3^{\sum_{t=s+1}^{t=i} n_t} \cdot 2^{\sum_{t=2}^{t=s-1} \alpha_t} \cdot (3^{n_s} - 2^{n_s})$$

$$j_{o.1}^{\beta} = 3^{\sum_{t=2}^{t=i} n_t} \cdot j_{o.i}^{\delta}$$

Therefore, we can obtain the final element of Cycle 1 using Identity (A3.2).

$$a_F^1 = 3^{n_1} \cdot \left[ 2^{\sum_{t=2}^{t=i} \alpha_t} \cdot Q_{O.i} - \frac{2^{3^{\sum_{t=1}^{t=i} n_t - 1} \cdot j_{o.i}^{\delta}} + 1}{3^{\sum_{t=1}^{t=i} n_t}} \cdot \sum_{s=2}^{s=i} 3^{\sum_{t=s+1}^{t=i} n_t} \cdot 2^{\sum_{t=2}^{t=s-1} \alpha_t} \cdot (3^{n_s} - 2^{n_s}) \right] - 2^{3^{\sum_{t=1}^{t=i} n_t - 1} \cdot j_{o.i}^{\delta} - \alpha_1} \cdot (3^{n_1} - 2^{n_1})$$

We need now to obtain the initial element of the second cycle - $a_O^2$. We follow the same procedure as before and modify the identity above to obtain an expression whose format coincides with Identity (A1.1). To do so, we first insert the factor $3^{n_1}$ inside the brackets

$$a_F^1 = 3^{n_1} \cdot 2^{\sum_{t=2}^{t=i} \alpha_t} \cdot Q_{O.i} - \frac{2^{3^{\sum_{t=1}^{t=i} n_t - 1} \cdot j_{o.i}^{\delta}} + 1}{3^{\sum_{t=2}^{t=i} n_t}} \cdot \sum_{s=2}^{s=i} 3^{\sum_{t=s+1}^{t=i} n_t} \cdot 2^{\sum_{t=2}^{t=s-1} \alpha_t} \cdot (3^{n_s} - 2^{n_s}) - 2^{3^{\sum_{t=1}^{t=i} n_t - 1} \cdot j_{o.i}^{\delta} - \alpha_1} \cdot (3^{n_1} - 2^{n_1})$$

Then, we extract the element "s=2" of the summation

$$a_F^1 = 3^{n_1} \cdot 2^{\sum_{t=2}^{t=i} \alpha_t} \cdot Q_{O.i} - \frac{2^{3^{\sum_{t=1}^{t=i} n_t - 1} \cdot j_{o.i}^{\delta}} + 1}{3^{\sum_{t=2}^{t=i} n_t}} \cdot \sum_{s=3}^{s=i} 3^{\sum_{t=s+1}^{t=i} n_t} \cdot 2^{\sum_{t=2}^{t=s-1} \alpha_t} \cdot (3^{n_s} - 2^{n_s}) - \frac{2^{3^{\sum_{t=1}^{t=i} n_t - 1} \cdot j_{o.i}^{\delta}} + 1}{3^{\sum_{t=2}^{t=i} n_t}} \cdot 3^{\sum_{t=3}^{t=i} n_t} \cdot (3^{n_2} - 2^{n_2}) - 2^{3^{\sum_{t=1}^{t=i} n_t - 1} \cdot j_{o.i}^{\delta} - \alpha_1} \cdot (3^{n_1} - 2^{n_1})$$

As before, we simplify by eliminating the factor $3^{\sum_{t=3}^{t=i} n_t}$ that it is also in the denominator

$$(A6.5) \qquad a_F^1 = 3^{n_1} \cdot 2^{\sum_{t=2}^{t=i} \alpha_t} \cdot Q_{O.i} - \frac{2^{3^{\sum_{t=1}^{t=i} n_t - 1} \cdot j_{o.i}^{\delta}} + 1}{3^{\sum_{t=2}^{t=i} n_t}} \cdot \sum_{s=3}^{s=i} 3^{\sum_{t=s+1}^{t=i} n_t} \cdot 2^{\sum_{t=2}^{t=s-1} \alpha_t} \cdot (3^{n_s} - 2^{n_s})$$

$$-\frac{2^{3^{\sum_{t=1}^{t=i} n_t - 1} \cdot j_{o.i}^{\delta}} + 1}{3^{n_2}} \cdot (3^{n_2} - 2^{n_2})$$
$$-2^{3^{\sum_{t=1}^{t=i} n_t - 1} \cdot j_{o.i}^{\delta} - \alpha_1} \cdot (3^{n_1} - 2^{n_1})$$

Since $j_{o.i}^{\delta}$ complies with Inequality (A6.2)

$$(A6.2) \qquad 3^{\sum_{t=1}^{t=i} n_t - 1} \cdot j_{o.i}^{\delta} > \sum_{t=1}^{t=i} \alpha_t$$

it is also true

$$3^{\sum_{t=1}^{t=i} n_t - 1} \cdot j_{o.i}^{\delta} > \alpha_1 + \alpha_2$$

Therefore, we can extract the factor $2^{\alpha_2}$ from the last element of the Identity (A6.5) above.

As a result, we can extract the common factor $2^{\alpha_2}$ from all the remaining elements of the summation, also, the first and the last elements of Identity (A6.5) above. We obtain:

$$a_F^1 = 2^{\alpha_2} \cdot \Bigg[ 3^{n_1} \cdot 2^{\sum_{t=3}^{t=i} \alpha_t} \cdot Q_{O.i} - \frac{2^{3^{\sum_{t=1}^{t=i} n_t - 1} \cdot j_{o.i}^{\delta}} + 1}{3^{\sum_{t=2}^{t=i} n_t}} \cdot \sum_{s=3}^{s=i} 3^{\sum_{t=s+1}^{t=i} n_t} \cdot 2^{\sum_{t=3}^{t=s-1} \alpha_t} \cdot (3^{n_s} - 2^{n_s})$$
$$- 2^{3^{\sum_{t=1}^{t=i} n_t - 1} \cdot j_{o.i}^{\delta} - \alpha_1 - \alpha_2} \cdot (3^{n_1} - 2^{n_1}) \Bigg]$$
$$-\frac{2^{3^{\sum_{t=1}^{t=i} n_t - 1} \cdot j_{o.i}^{\delta}} + 1}{3^{n_2}} \cdot (3^{n_2} - 2^{n_2})$$

In a similar way as before, we extract $3^{n_2 - 1}$ from the power in the numerator of the last fraction

$$(A6.6) \quad a_F^1 = 2^{\alpha_2} \cdot \Bigg[ 3^{n_1} \cdot 2^{\sum_{t=3}^{t=i} \alpha_t} \cdot Q_{O.i} - \frac{2^{3^{\sum_{t=1}^{t=i} n_t - 1} \cdot j_{o.i}^{\delta}} + 1}{3^{\sum_{t=2}^{t=i} n_t}} \cdot \sum_{s=3}^{s=i} 3^{\sum_{t=s+1}^{t=i} n_t} \cdot 2^{\sum_{t=3}^{t=s-1} \alpha_t} \cdot (3^{n_s} - 2^{n_s})$$
$$- 2^{3^{\sum_{t=1}^{t=i} n_t - 1} \cdot j_{o.i}^{\delta} - \alpha_1 - \alpha_2} \cdot (3^{n_1} - 2^{n_1}) \Bigg]$$
$$-\frac{2^{3^{n_2 - 1} \cdot \left(3^{n_1 + \sum_{t=3}^{t=i} n_t} \cdot j_{o.i}^{\delta}\right)} + 1}{3^{n_2}} \cdot (3^{n_2} - 2^{n_2})$$

As in the case before, the Identity (A6.6) above follows the same format as Identity (A1.1) and represents the first odd number of the second cycle.

$$a_F^1 = a_O^2$$

$$a_O^2 = 2^{\alpha_2} \cdot \Bigg[ 3^{n_1} \cdot 2^{\sum_{t=3}^{t=i} \alpha_t} \cdot Q_{O.i} - \frac{2^{3^{\sum_{t=1}^{t=i} n_t - 1} \cdot j_{o.i}^{\delta}} + 1}{3^{\sum_{t=2}^{t=i} n_t}} \cdot \sum_{s=3}^{s=i} 3^{\sum_{t=s+1}^{t=i} n_t} \cdot 2^{\sum_{t=3}^{t=s-1} \alpha_t} \cdot (3^{n_s} - 2^{n_s})$$

$$- 2^{3^{\sum_{t=1}^{t=i} n_t - 1} \cdot j_{o.i}^{\delta} - \alpha_1 - \alpha_2} \cdot (3^{n_1} - 2^{n_1}) \Bigg]$$

$$- \frac{2^{3^{n_2 - 1} \cdot \left(3^{n_1 + \sum_{t=3}^{t=i} n_t} \cdot j_{o.i}^{\delta}\right)} + 1}{3^{n_2}} \cdot (3^{n_2} - 2^{n_2})$$

As a result, the final element of this cycle is Identity (A3.2). We apply Identity (A3.2) and obtain the final element of the second cycle

$$a_F^2 = 3^{n_2} \cdot \Bigg[ 3^{n_1} \cdot 2^{\sum_{t=3}^{t=i} \alpha_t} \cdot Q_{O.i} - \frac{2^{3^{\sum_{t=1}^{t=i} n_t - 1} \cdot j_{o.i}^{\delta}} + 1}{3^{\sum_{t=2}^{t=i} n_t}} \cdot \sum_{s=3}^{s=i} 3^{\sum_{t=s+1}^{t=i} n_t} \cdot 2^{\sum_{t=3}^{t=s-1} \alpha_t} \cdot (3^{n_s} - 2^{n_s})$$

$$- 2^{3^{\sum_{t=1}^{t=i} n_t - 1} \cdot j_{o.i}^{\delta} - \alpha_1 - \alpha_2} \cdot (3^{n_1} - 2^{n_1}) \Bigg]$$

$$- 2^{3^{\sum_{t=1}^{t=i} n_t - 1} \cdot j_{o.i}^{\delta} - \alpha_2} \cdot (3^{n_2} - 2^{n_2})$$

We follow the same procedure as before. We first insert the factor $3^{n_2}$ into the brackets, we extract the element "s=3" of the summation, and then we extract the common factor $2^{\alpha_3}$. We obtain:

$$a_F^2 = 2^{\alpha_3} \cdot \Bigg[ 3^{n_1 + n_2} \cdot 2^{\sum_{t=4}^{t=i} \alpha_t} \cdot Q_{O.i} - \frac{2^{3^{\sum_{t=1}^{t=i} n_t - 1} \cdot j_{o.i}^{\delta}} + 1}{3^{\sum_{t=3}^{t=i} n_t}} \cdot \sum_{s=4}^{s=i} 3^{\sum_{t=s+1}^{t=i} n_t} \cdot 2^{\sum_{t=4}^{t=s-1} \alpha_t} \cdot (3^{n_s} - 2^{n_s})$$

$$- 2^{3^{\sum_{t=1}^{t=i} n_t - 1} \cdot j_{o.i}^{\delta} - \alpha_1 - \alpha_2 - \alpha_3} \cdot 3^{n_2} \cdot (3^{n_1} - 2^{n_1}) - 2^{3^{\sum_{t=1}^{t=i} n_t - 1} \cdot j_{o.i}^{\delta} - \alpha_2 - \alpha_3} \cdot (3^{n_2} - 2^{n_2}) \Bigg]$$

$$- \frac{2^{3^{\sum_{t=1}^{t=i} n_t - 1} \cdot j_{o.i}^{\delta}} + 1}{3^{n_3}} \cdot (3^{n_3} - 2^{n_3})$$

The identity above, as in the previous cases, follows the same format as Identity (A1.1) and represents the initial element of the third cycle

$$a_F^2 = a_O^3$$

As a result, as in the previous cases, the final element of this cycle is Identity (A3.2). We apply Identity (A3.2) and obtain the final element of the third cycle

$$a_F^3 = 3^{n_3} \cdot \Bigg[ 3^{n_1 + n_2} \cdot 2^{\sum_{t=4}^{t=i} \alpha_t} \cdot Q_{O.i} - \frac{2^{3^{\sum_{t=1}^{t=i} n_t - 1} \cdot j_{o.i}^{\delta}} + 1}{3^{\sum_{t=3}^{t=i} n_t}} \cdot \sum_{s=4}^{s=i} 3^{\sum_{t=s+1}^{t=i} n_t} \cdot 2^{\sum_{t=4}^{t=s-1} \alpha_t} \cdot (3^{n_s} - 2^{n_s})$$

$$- 2^{3^{\sum_{t=1}^{t=i} n_t - 1} \cdot j_{o.i}^{\delta} - \alpha_1 - \alpha_2 - \alpha_3} \cdot 3^{n_2} \cdot (3^{n_1} - 2^{n_1}) - 2^{3^{\sum_{t=1}^{t=i} n_t - 1} \cdot j_{o.i}^{\delta} - \alpha_2 - \alpha_3} \cdot (3^{n_2} - 2^{n_2}) \Bigg]$$

$$- 2^{3^{\sum_{t=1}^{t=i} n_t - 1} \cdot j_{o.i}^{\delta} - \alpha_3} \cdot (3^{n_3} - 2^{n_3})$$

If we repeat the same process “i-1” times we obtain

$$a_F^{i-1} = 3^{n_{i-1}} \cdot \left[ 3^{\sum_{t=1}^{t=i-2} n_t} \cdot 2^{\alpha_i} \cdot Q_{O.i} - \frac{2^{3^{\sum_{t=1}^{t=i} n_{t-1}} \cdot j_{o.i}^{\delta}} + 1}{3^{n_i + n_{i-1}}} \cdot \left(3^{n_i} - 2^{n_i}\right) \right.$$

$$\left. - \sum_{s=1}^{s=i-2} 2^{3^{\sum_{t=1}^{t=i} n_{t-1}} \cdot j_{o.i}^{\delta} - \sum_{t=s}^{t=i-1} \alpha_t} \cdot 3^{\sum_{t=s+1}^{t=i-2} n_t} \cdot \left(3^{n_s} - 2^{n_s}\right) \right]$$

$$-2^{3^{\sum_{t=1}^{t=i} n_{t-1}} \cdot j_{o.i}^{\delta} - \alpha_{i-1}} \cdot \left(3^{n_{i-1}} - 2^{n_{i-1}}\right)$$

We repeat the same process one last time. We first insert the factor $3^{n_{i-1}}$ into the brackets, we extract the element “s=i” of the summation, and then we extract the common factor $2^{\alpha_i}$. We obtain:

$$a_F^{i-1} = 2^{\alpha_i} \cdot \left[ 3^{\sum_{t=1}^{t=i-1} n_t} \cdot Q_{O.i} - \sum_{s=1}^{s=i-1} 2^{3^{\sum_{t=1}^{t=i} n_{t-1}} \cdot j_{o.i}^{\delta} - \sum_{t=s}^{t=i} \alpha_t} \cdot 3^{\sum_{t=s+1}^{t=i-1} n_t} \cdot \left(3^{n_s} - 2^{n_s}\right) \right]$$

$$- \frac{2^{3^{\sum_{t=1}^{t=i} n_{t-1}} \cdot j_{o.i}^{\delta}} + 1}{3^{n_i}} \cdot \left(3^{n_i} - 2^{n_i}\right)$$

The identity above, as in the previous cases, follows the same format as Identity (A1.1) and represents the initial element of Cycle “i”

$$a_F^{i-1} = a_O^i$$

$$a_O^i = 2^{\alpha_i} \cdot \left[ 3^{\sum_{t=1}^{t=i-1} n_t} \cdot Q_{O.i} - \sum_{s=1}^{s=i-1} 2^{3^{\sum_{t=1}^{t=i} n_{t-1}} \cdot j_{o.i}^{\delta} - \sum_{t=s}^{t=i} \alpha_t} \cdot 3^{\sum_{t=s+1}^{t=i-1} n_t} \cdot \left(3^{n_s} - 2^{n_s}\right) \right]$$

$$- \frac{2^{3^{\sum_{t=1}^{t=i} n_{t-1}} \cdot j_{o.i}^{\delta}} + 1}{3^{n_i}} \cdot \left(3^{n_i} - 2^{n_i}\right)$$

As a result, as in the previous cases, the final element of this cycle is Identity (A3.2). We apply Identity (A3.2) and obtain the final element of Cycle “i”

$$a_F^{i} = 3^{n_i} \cdot \left[ 3^{\sum_{t=1}^{t=i-1} n_t} \cdot Q_{O.i} - \sum_{s=1}^{s=i-1} 2^{3^{\sum_{t=1}^{t=i} n_{t-1}} \cdot j_{o.i}^{\delta} - \sum_{t=s}^{t=i} \alpha_t} \cdot 3^{\sum_{t=s+1}^{t=i-1} n_t} \cdot \left(3^{n_s} - 2^{n_s}\right) \right]$$

$$-2^{3^{\sum_{t=1}^{t=i} n_{t-1}} \cdot j_{o.i}^{\delta} - \alpha_i} \cdot \left(3^{n_i} - 2^{n_i}\right)$$

We rearrange the above. We insert the factor $3^{n_i}$ within the brackets

$$a_F^{i} = 3^{\sum_{t=1}^{t=i} n_t} \cdot Q_{O.i} - \sum_{s=1}^{s=i-1} 2^{3^{\sum_{t=1}^{t=i} n_{t-1}} \cdot j_{o.i}^{\delta} - \sum_{t=s}^{t=i} \alpha_t} \cdot 3^{\sum_{t=s+1}^{t=i} n_t} \cdot \left(3^{n_s} - 2^{n_s}\right)$$

$$-2^{3^{\sum_{t=1}^{t=i} n_{t-1}} \cdot j_{o.i}^{\delta} - \alpha_i} \cdot \left(3^{n_i} - 2^{n_i}\right)$$

We insert the last element of the identity above in the summation

$$a_F^i = 3^{\sum_{t=1}^{t=i} n_t} \cdot Q_{O.i} - \sum_{s=1}^{s=i} 2^{3^{\sum_{t=1}^{t=i} n_t - 1} \cdot j_{o.i}^{\delta} - \sum_{t=s}^{t=i} \alpha_t} \cdot 3^{\sum_{t=s+1}^{t=i} n_t} \cdot (3^{n_s} - 2^{n_s})$$

Since

$$2^{3^{\sum_{t=1}^{t=i} n_t - 1} \cdot j_{o.i}^{\delta} - \sum_{t=s}^{t=i} \alpha_t} = 2^{3^{\sum_{t=1}^{t=i} n_t - 1} \cdot j_{o.i}^{\delta} - \sum_{t=1}^{t=i} \alpha_t} \cdot 2^{\sum_{t=1}^{t=s-1} \alpha_t}$$

We can extract the common factor

$$2^{3^{\sum_{t=1}^{t=i} n_t - 1} \cdot j_{o.i}^{\delta} - \sum_{t=1}^{t=i} \alpha_t}$$

from each addend of the summation and leave inside

$$2^{\sum_{t=1}^{t=s-1} \alpha_t}$$

and obtain

$$a_F^i = 3^{\sum_{t=1}^{t=i} n_t} \cdot Q_{O.i} - 2^{3^{\sum_{t=1}^{t=i} n_t - 1} \cdot j_{o.i}^{\delta} - \sum_{t=1}^{t=i} \alpha_t} \cdot \sum_{s=1}^{s=i} 2^{\sum_{t=1}^{t=s-1} \alpha_t} \cdot 3^{\sum_{t=s+1}^{t=i} n_t} \cdot (3^{n_s} - 2^{n_s})$$

We flip around the order of the power of 2 and the power of 3 in the summation above and obtain Identity (A6.3)

$$(A6.3) \qquad a_F^i = 3^{\sum_{t=1}^{t=i} n_t} \cdot Q_{O.i} - 2^{3^{\sum_{t=1}^{t=i} n_t - 1} \cdot j_{o.i}^{\delta} - \sum_{t=1}^{t=i} \alpha_t} \cdot \sum_{s=1}^{s=i} 3^{\sum_{t=s+1}^{t=i} n_t} \cdot 2^{\sum_{t=1}^{t=s-1} \alpha_t} \cdot (3^{n_s} - 2^{n_s})$$

$$3^{\sum_{t=1}^{t=i} n_t} \cdot j_{o.i}^{\delta} > \sum_{t=1}^{t=i} \alpha_t$$

Which is what we wanted to prove.

## Example

Let's say we want to generate a sequence made of three cycles with the following parameters:

| Cycle s | $n_s$ | $\alpha_s$ |
|---|---|---|
| 1 | 1 | 4 |
| 2 | 1 | 4 |
| 3 | 2 | 3 |

Note that $j_{o.i}^{\delta}$ can be any odd number but must be such that

$$3^{\sum_{t=1}^{t=i} n_t - 1} \cdot j_{o.i}^{\delta} > \sum_{t=1}^{t=i} \alpha_t$$

Since $\sum_{t=1}^{t=i} \alpha_t = 11$ and $3^{\sum_{t=1}^{t=i} n_t - 1} = 27$ , then the smallest odd-number $j_O^{\delta}$ that makes the above possible is for

$$j_{o.i}^{\delta} = 1$$

The initial element of this sequence is

$$(A6.1) \qquad a_O^1 = 2^{\sum_{t=1}^{t=i} \alpha_t} \cdot Q_{O.i} - \frac{2^{3^{\sum_{t=1}^{t=i} n_t - 1} \cdot j_{o.i}^{\delta}} + 1}{3^{\sum_{t=1}^{t=i} n_t}} \cdot \sum_{s=1}^{i} 3^{\sum_{t=s+1}^{t=i} n_t} \cdot 2^{\sum_{t=1}^{t=s-1} \alpha_t} \cdot (3^{n_s} - 2^{n_s})$$

$$a_O^1 = 2^{4+4+3} \cdot Q_{O.3} - \frac{2^{3^{2+1+1-1}} + 1}{3^{2+1+1}} \cdot \left(3^{2+1} \cdot (3^1 - 2^1) + 3^2 \cdot 2^4 \cdot (3^1 - 2^1) + 2^{4+4} \cdot (3^2 - 2^2)\right)$$

$$a_O^1 = 2^{11} \cdot Q_{O.3} - 1\,657\,009 \times (27 + 144 + 1280)$$

$$a_O^1 = 2^{11} \cdot Q_{O.3} - 1\,657\,009 \times 1451$$

The smallest $Q_O$ that complies with the above is

$$Q_{O.3} = 1\,173\,985$$

Therefore:

$$\boxed{a_O^1 = 1\,221}$$

The final element of this sequence after "$i = 3$" cycles is

vpadilla@aertecsolutions.com

$$(A6.3)\quad a_F^i = 3^{\sum_{t=1}^{t=i} n_t} \cdot Q_{O.i} - 2^{3^{\sum_{t=1}^{t=i} n_t - 1} \cdot j_{o.i}^{\delta} - \sum_{t=1}^{t=i} \alpha_t} \cdot \sum_{s=1}^{i} 3^{\sum_{t=s+1}^{t=i} n_t} \cdot 2^{\sum_{t=1}^{t=s-1} \alpha_t} \cdot (3^{n_s} - 2^{n_s})$$

$$a_F^3 = 3^{1+1+2} \cdot 1\,173\,985 - 2^{3^{1+1+2-1} - (4+4+3)}[3^{1+2} \cdot (3^{n_1} - 2^{n_1}) + 3^2 \cdot 2^4 \cdot (3^1 - 2^1) + 2^{4+4} \cdot (3^{n_2} - 2^{n_2})]$$

$$a_F^3 = 81 \times 1\,173\,985 - 65\,536 \times 1451$$

$$\boxed{a_F^3 = 49}$$

Then, the sequence generated by

$$a_O^1 = 1\,221$$

reaches

$$a_F^3 = 49$$

after three cycles

| | Steps up $n_s$ | Steps down $\alpha_s$ |
|---|---|---|
| Cycle 1 | 1 | 4 |
| Cycle 2 | 1 | 4 |
| Cycle 3 | 2 | 3 |

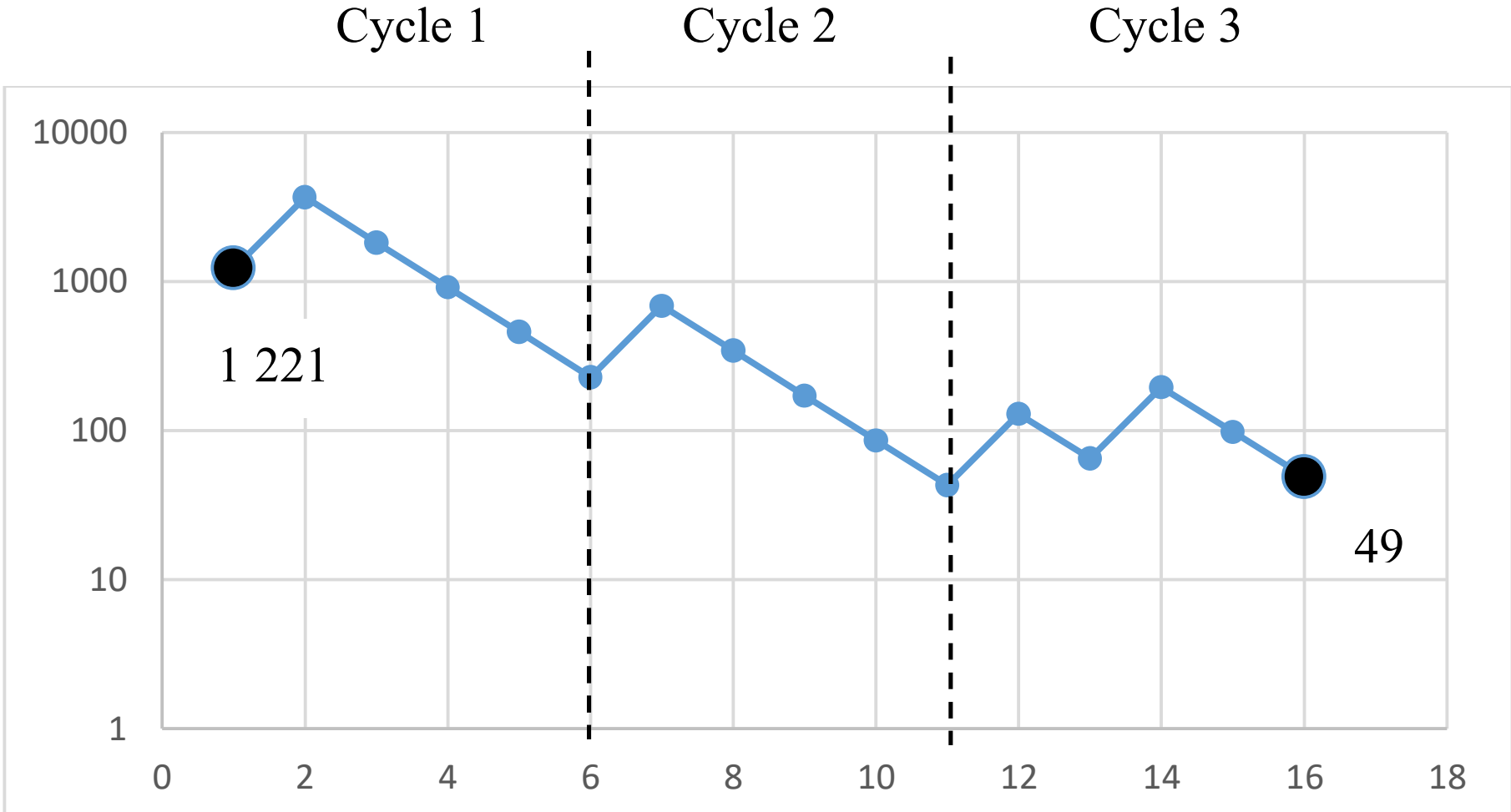

Note: The graph is in logarithm scale to facilitate the visualization of upward and downward steps.

Please, note that there are infinite numbers that satisfy the combination of $n_s, \alpha_s$ above and, as a result, follow the same pattern of upward and downward steps. This is the case since we picked the “smallest” $Q_{O.3}$; but the following numbers also satisfy the above.

$$a_O^1 = 2^{11} \times Q_{O.3} - 1\,657\,009 \times 1451$$

$$a_F^3 = 81 \times Q_{O.3} - 65\,536 \times 1451$$

## Lemma A6.2

Let $a_O^1$ - the initial element of a Sequence A made of "$i$" cycles – and $a_F^i$ the final element of this sequence after "i" cycles.

The following identity is true:

$$(A6.4) \qquad 2^{\sum_{t=1}^{t=i}\alpha_t}\cdot a_F^i - 3^{\sum_{t=1}^{t=i} n_t}\cdot a_O^1 = \sum_{s=1}^{i} 3^{\sum_{t=s+1}^{t=i} n_t}\cdot 2^{\sum_{t=1}^{t=s-1}\alpha_t}\cdot\left(3^{n_s}-2^{n_s}\right)$$

## Proof

According to Lemma A6.1, the initial element - $a_O^1$ - and the final element - $a_F^i$ - of a Sequence A made of "$i$" cycles are:

$$(A6.1) \quad a_O^1 = 2^{\sum_{t=1}^{t=i}\alpha_t}\cdot Q_{O.i} - \frac{2^{3^{\sum_{t=1}^{t=i} n_t-1}\cdot j_{o.i}^{\delta}}+1}{3^{\sum_{t=1}^{t=i} n_t}}\cdot\sum_{s=1}^{i} 3^{\sum_{t=s+1}^{t=i} n_t}\cdot 2^{\sum_{t=1}^{t=s-1}\alpha_t}\cdot\left(3^{n_s}-2^{n_s}\right)$$

$$(A6.3) \qquad a_F^i = 3^{\sum_{t=1}^{t=i} n_t}\cdot Q_{O.i} - 2^{3^{\sum_{t=1}^{t=i} n_t-1}\cdot j_{o.i}^{\delta}-\sum_{t=1}^{t=i}\alpha_t}\cdot\sum_{s=1}^{s=i} 3^{\sum_{t=s+1}^{t=i} n_t}\cdot 2^{\sum_{t=1}^{t=s-1}\alpha_t}\cdot\left(3^{n_s}-2^{n_s}\right)$$

$$Q_{O.i};\, j_{o.i}^{\delta}\in\{odd\}$$

$$(A6.2) \qquad 3^{\sum_{t=1}^{t=i} n_t-1}\cdot j_{o.i}^{\delta} > \sum_{t=1}^{t=i}\alpha_t$$

We multiply Identity (A6.3) by $2^{\sum_{t=1}^{t=i}\alpha_t}$

$$2^{\sum_{t=1}^{t=i}\alpha_t}\cdot a_F^i = 2^{\sum_{t=1}^{t=i}\alpha_t}\cdot 3^{\sum_{t=1}^{t=i} n_t}\cdot Q_{O.i} - 2^{3^{\sum_{t=1}^{t=i} n_t-1}\cdot j_{o.i}^{\delta}}\cdot\sum_{s=1}^{s=i} 3^{\sum_{t=s+1}^{t=i} n_t}\cdot 2^{\sum_{t=1}^{t=s-1}\alpha_t}\cdot\left(3^{n_s}-2^{n_s}\right)$$

and multiply Identity (A6.1) by $3^{\sum_{t=1}^{t=i} n_t}$

$$3^{\sum_{t=1}^{t=i} n_t}\cdot a_O^1 = 3^{\sum_{t=1}^{t=i} n_t}\cdot 2^{\sum_{t=1}^{t=i}\alpha_t}\cdot Q_{O.i} - \left(2^{3^{\sum_{t=1}^{t=i} n_t-1}\cdot j_{o.i}^{\delta}}+1\right)\cdot\sum_{s=1}^{i} 3^{\sum_{t=s+1}^{t=i} n_t}\cdot 2^{\sum_{t=1}^{t=s-1}\alpha_t}\cdot\left(3^{n_s}-2^{n_s}\right)$$

we subtract the latter from the former and obtain:

$$(A6.4) \qquad 2^{\sum_{t=1}^{t=i}\alpha_t}\cdot a_F^i - 3^{\sum_{t=1}^{t=i} n_t}\cdot a_O^1 = \sum_{s=1}^{s=i} 3^{\sum_{t=s+1}^{t=i} n_t}\cdot 2^{\sum_{t=1}^{t=s-1}\alpha_t}\cdot\left(3^{n_s}-2^{n_s}\right)$$

Which is what we wanted to prove.

## Lemma A6.3

Let $a_O^1$ - the initial element of a Sequence A made of "$i$" cycles – and $a_F^i$ the final element of this sequence after “i” cycles.

The following identity is true:

$$(A6.5) \qquad Q_{O.i} = \frac{2^{3^{\sum_{t=1}^{t=i} n_t - 1}\cdot j_{o.i}^{\delta}} + 1}{3^{\sum_{t=1}^{t=i} n_t}} \cdot a_F^i - 2^{3^{\sum_{t=1}^{t=i} n_t - 1}\cdot j_{o.i}^{\delta} - \sum_{t=1}^{t=i} \alpha_t} \cdot a_O^1$$

## Proof

According to Lemma A6.1, the initial element - $a_O^1$ - and the final element - $a_F^i$ - of a Sequence A made of "$i$" cycles are:

$$(A6.1) \quad a_O^1 = 2^{\sum_{t=1}^{t=i} \alpha_t} \cdot Q_{O.i} - \frac{2^{3^{\sum_{t=1}^{t=i} n_t - 1}\cdot j_{o.i}^{\delta}} + 1}{3^{\sum_{t=1}^{t=i} n_t}} \cdot \sum_{s=1}^{i} 3^{\sum_{t=s+1}^{t=i} n_t} \cdot 2^{\sum_{t=1}^{t=s-1} \alpha_t} \cdot (3^{n_s} - 2^{n_s})$$

$$(A6.3) \qquad a_F^i = 3^{\sum_{t=1}^{t=i} n_t} \cdot Q_{O.i} - 2^{3^{\sum_{t=1}^{t=i} n_t - 1}\cdot j_{o.i}^{\delta} - \sum_{t=1}^{t=i} \alpha_t} \cdot \sum_{s=1}^{s=i} 3^{\sum_{t=s+1}^{t=i} n_t} \cdot 2^{\sum_{t=1}^{t=s-1} \alpha_t} \cdot (3^{n_s} - 2^{n_s})$$

$$Q_{O.i}; j_{o.i}^{\delta} \in \{odd\}$$

$$(A6.2) \qquad 3^{\sum_{t=1}^{t=i} n_t - 1} \cdot j_{o.i}^{\delta} > \sum_{t=1}^{t=i} \alpha_t$$

We multiply Identity (A6.3) by $\frac{2^{3^{\sum_{t=1}^{t=i} n_t - 1}\cdot j_{o.i}^{\delta}+1}}{3^{\sum_{t=1}^{t=i} n_t}}$ and Identity (A6.1) by $2^{3^{\sum_{t=1}^{t=i} n_t - 1}\cdot j_{o.i}^{\delta} - \sum_{t=1}^{t=i} \alpha_t}$ and subtract the former from the latter:

$$(A6.5) \qquad Q_{O.i} = \frac{2^{3^{\sum_{t=1}^{t=i} n_t - 1}\cdot j_{o.i}^{\delta}} + 1}{3^{\sum_{t=1}^{t=i} n_t}} \cdot a_F^i - 2^{3^{\sum_{t=1}^{t=i} n_t - 1}\cdot j_{o.i}^{\delta} - \sum_{t=1}^{t=i} \alpha_t} \cdot a_O^1$$

Which is what we wanted to prove.

# Appendix 7 – Closing Cycle is always possible.

## Lemma A7.1

For any combination of "k-1" cycles made of different set of parameters:

$$\{n_t; \alpha_t\}$$
$$t \in \{1,2,\dots,k-1\}$$

$$n_t; \alpha_t \in \mathbb{N}$$
$$\alpha_t > n_t$$

and an additional parameter "$n_k$" that represents the consecutive upward steps of an additional Cycle-k.

We can always generate a sequence that converges to 1 in this additional and final cycle.

## Proof

From previous Appendix 6, the initial element of the sequence that combines "k-1" cycles together is Identity (A6.1)

$$(A6.1)\ \ a_O^1 = 2^{\sum_{t=1}^{t=k-1}\alpha_t} \cdot Q_{O.k-1} - \frac{2^{3^{\sum_{t=1}^{t=k-1} n_t - 1} \cdot j_{o.k-1}^{\delta}} + 1}{3^{\sum_{t=1}^{t=k-1} n_t}} \cdot \sum_{s=1}^{k-1} 3^{\sum_{t=s+1}^{t=k-1} n_t} \cdot 2^{\sum_{t=1}^{t=s-1}\alpha_t} \cdot (3^{n_s} - 2^{n_s})$$

And the inequality (A6.2) applies:

$$(A6.2)\qquad 3^{\sum_{t=1}^{t=k-1} n_t} \cdot j_{o.k-1}^{\delta} > \sum_{t=1}^{t=k-1} \alpha_t$$

and the final element of the sequence after $k-1$ cycles is identity (A6.3)

$$(A6.3)\ \ a_F^{k-1} = 3^{\sum_{t=1}^{t=k-1} n_t} \cdot Q_{O.k-1} - 2^{3^{\sum_{t=1}^{t=k-1} n_t - 1} \cdot j_{o.k-1}^{\delta} - \sum_{t=1}^{t=k-1}\alpha_t} \cdot \sum_{s=1}^{k-1} 3^{\sum_{t=s+1}^{t=k-1} n_t} \cdot 2^{\sum_{t=1}^{t=s-1}\alpha_t} \cdot (3^{n_s} - 2^{n_s})$$

Please be aware that the initial element of the sequence $a_O^1$ is not defined yet since the parameter $Q_{O.k-1}$ is not defined yet either. In other words, at this point we are still "looking for" the sequence that converges to 1 with the given set of parameters $n_t$ and $\alpha_t$.

In Appendix 3, Lemma A3.2 proves that any cycle whose initial element is:

$$(A3.7) \quad a_O^c = 2^{n_c} \cdot \frac{2^{3^{n_c-1} \cdot q_O} + 1}{3^{n_c}} - 1$$
$$q_O \in \{ODD\}$$

is a cycle whose final element is 1

$$(A3.8) \quad a_F = 1$$

Please, be aware, as we already said, that "$n_k$" - the number of upward steps – is already set. On the other hand, "$q_o$" - the number of downward steps after the upper bound of the closing cycle – is not established yet.

Since we want the sequence to converge to 1, the final element of the sequence after $k-1$ cycles must be the same as the initial element of the closing Cycle-k. Therefore:

$$a_F^{k-1} = a_O^c = a_O^k$$

And the number of upward steps in the closing Cycle-k

$$n_c = n_k$$

We equate Identity (A6.3) and Identity (A3.7) and obtain:

$$3^{\sum_{t=1}^{t=k-1} n_t} \cdot Q_{O.k-1} - 2^{3^{\sum_{t=1}^{t=k-1} n_t - 1} \cdot j_{o.k-1}^{\delta} - \sum_{t=1}^{t=k-1} \alpha_t} \cdot \sum_{s=1}^{k-1} 3^{\sum_{t=s+1}^{t=k-1} n_t} \cdot 2^{\sum_{t=1}^{t=s-1} \alpha_t} \cdot (3^{n_s} - 2^{n_s}) =$$
$$= 2^{n_k} \cdot \frac{2^{3^{n_k-1} \cdot q_O} + 1}{3^{n_k}} - 1$$

We reorganize the above. We pass the summation to the right side of the identity and the fraction to the left side of said identity.

$$(A7.1) \quad 3^{\sum_{t=1}^{t=k-1} n_t} \cdot Q_{O.k-1} - 2^{n_k} \cdot \frac{2^{3^{n_k-1} \cdot q_O} + 1}{3^{n_k}} =$$
$$= 2^{3^{\sum_{t=1}^{t=k-1} n_t - 1} \cdot j_{o.k-1}^{\delta} - \sum_{t=1}^{t=k-1} \alpha_t} \cdot \sum_{s=1}^{k-1} 3^{\sum_{t=s+1}^{t=k-1} n_t} \cdot 2^{\sum_{t=1}^{t=s-1} \alpha_t} \cdot (3^{n_s} - 2^{n_s}) - 1$$

For any different set of parameters

$$\{n_t; \alpha_t\}$$
$$t \in \{1,2,\dots,k\}$$

we obtain an odd number on the right side of the identity above.

$$(A7.2) \quad T_o = 2^{3^{\sum_{t=1}^{t=k-1} n_t - 1} \cdot j_{o.k-1}^{\delta} - \sum_{t=1}^{t=k-1} \alpha_t} \cdot \sum_{s=1}^{k-1} 3^{\sum_{t=s+1}^{t=k-1} n_t} \cdot 2^{\sum_{t=1}^{t=s-1} \alpha_t} \cdot (3^{n_s} - 2^{n_s}) - 1$$

In other words, the identity above can be rewritten as:

$$3^{\sum_{t=1}^{t=k-1} n_t} \cdot Q_{O.k-1} - 2^{n_k} \cdot \frac{2^{3^{n_k - 1} \cdot q_O} + 1}{3^{n_k}} = T_o = ODD$$

According to Lemma A5.6 from Appendix 5, the set of all odd numbers can be divided in $3^{\mathrm{m}}$ subsets as follows:

$$S_j \in \left\{ 3^m \cdot K_O - 2^n \cdot \frac{2^{3^{n-1} \cdot q_{o.j}} + 1}{3^n} \right\}$$

$$q_{o.j} \in \{1; 3; 5; 7; 11; \ldots.; 2 \cdot 3^m - 1\}$$

$$K_O \in \{Odd\}$$

$$j \in \{1,2,3,4,5, \ldots, 3^m\}$$

$$\forall\, n\,; m \in \mathbb{N}$$

In other words, any odd number $T_o$ such as the one defined in Identity (A7.2):

$$(A7.2) \quad T_o = 2^{3^{\sum_{t=1}^{t=k-1} n_t - 1} \cdot j_{o.k-1}^{\delta} - \sum_{t=1}^{t=k-1} \alpha_t} \cdot \sum_{s=1}^{k-1} 3^{\sum_{t=s+1}^{t=k-1} n_t} \cdot 2^{\sum_{t=1}^{t=s-1} \alpha_t} \cdot (3^{n_s} - 2^{n_s}) - 1$$

can be expressed as:

$$T_o = 3^m \cdot K_O - 2^n \cdot \frac{2^{3^{n-1} \cdot q_o} + 1}{3^n}$$

Since, according to Lemma A5.6, this is the case for any parameters “m” and “n”, we choose:

$$m = 3^{\sum_{t=1}^{t=k-1} n_t}$$
$$n = n_k$$

and therefore, there must always be a

$q_{o.j}$ and a $K_O$

that make Lemma A5.6 possible. Therefore, if we apply

$$Q_{O.k-1} = K_O$$
$$q_o = q_{o.j}$$

to Identity (A7.1), said identity is true for any set of parameters of an “k-1”-cycle sequence

$$\{n_t; \alpha_t\}$$
$$t \in \{1,2, \ldots, k-1\}$$

and the parameter of the closing Cycle-k

$$n_k$$

Which is what we wanted to prove.

We apply Identity (A3.8) from Appendix 3 and obtain the total number of downward steps in the closing cycle:

$$(A7.3) \quad \alpha_k = 3^{n_k-1} \cdot q_{o.j} + n_k$$

And the number of upward steps:

$$n_k$$

We apply the above to Identity (A6.1) and obtain the odd number $d_O^1$. In this document, whenever we use the letter D or d, such as Sequence D initiated by an odd number $d_O^1$, we use to refer to a sequence that converges to 1.

$$(A7.4)\ d_O^1 = 2^{\sum_{t=1}^{t=k-1}\alpha_t} \cdot K_O - \frac{2^{3^{\sum_{t=1}^{t=k-1} n_t-1} \cdot j_{o.k-1}^{\delta}} + 1}{3^{\sum_{t=1}^{t=k-1} n_t}} \cdot \sum_{s=1}^{k-1} 3^{\sum_{t=s+1}^{t=k-1} n_t} \cdot 2^{\sum_{t=1}^{t=s-1}\alpha_t} \cdot (3^{n_s} - 2^{n_s})$$

And the final element of the sequence after “k-1” cycles is:

$$(A7.5) \quad d_F^{k-1} = 3^{\sum_{t=1}^{t=k-1} n_t} \cdot K_O$$
$$- 2^{3^{\sum_{t=1}^{t=k-1} n_t-1} \cdot j_{o.k-1}^{\delta} - \sum_{t=1}^{t=k-1}\alpha_t} \sum_{s=1}^{k-1} 3^{\sum_{t=s+1}^{t=k-1} n_t} \cdot 2^{\sum_{t=1}^{t=s-1}\alpha_t} \cdot (3^{n_s} - 2^{n_s}) =$$
$$= 2^{n_k} \cdot \frac{2^{3^{n_k-1} \cdot q_{O.j}} + 1}{3^{n_k}} - 1$$

Therefore, the Sequence D generated by this number converges to 1 in the final Cycle-k.

If we make the following change:

$$K_O = 2^{3^{n_k-1} \cdot q_o + n_k} \cdot Q_{O.d} - \frac{2^{3^{\sum_{t=1}^{t=k} n_t-1} \cdot j_{o.k}^{\delta}} + 1}{3^{\sum_{t=1}^{t=k} n_t}} \cdot (3^{n_k} - 2^{n_k})$$

We can rewrite Identity (A7.5) for an odd number $d_O^1$ that generates a Sequence D and whose final element after “k” cycles is 1, we obtain:

$$(A7.6) \quad d_O^1 = 2^{\sum_{t=1}^{t=k-1}\alpha_t + 3^{n_k-1}\cdot q_o + n_k} \cdot Q_{O.d} - \frac{2^{3^{\sum_{t=1}^{t=k} n_t - 1}\cdot j_{o.k}^{\delta}} + 1}{3^{\sum_{t=1}^{t=k} n_t}} \cdot \sum_{s=1}^{k} 3^{\sum_{t=s+1}^{t=k} n_t} \cdot 2^{\sum_{t=1}^{t=s-1}\alpha_t} \cdot \left(3^{n_s} - 2^{n_s}\right)$$

$$d_F^k = 1 = 3^{\sum_{t=1}^{t=k} n_t} \cdot Q_{O.d} - 2^{3^{\sum_{t=1}^{t=k} n_t - 1}\cdot j_{o.k}^{\delta} - \left(\sum_{t=1}^{t=k-1}\alpha_t + 3^{n_k-1}\cdot q_o + n_k\right)} \cdot \left[\sum_{s=1}^{k} 3^{\sum_{t=s+1}^{t=k} n_t} \cdot 2^{\sum_{t=1}^{t=s-1}\alpha_t} \cdot \left(3^{n_s} - 2^{n_s}\right)\right]$$

Please be aware that the parameters $q_o$ and $Q_{O.d}$ of Cycle-k are related and unique for each set of parameters $\alpha_t$ and $n_t$ of the previous "k-1" cycles.

Also, since

$$q_{o.j} \in \{1; 3; 5; 7; 11; \ldots.; 2 \cdot 3^m - 1\}$$

And

$$m = 3^{\sum_{t=1}^{t=k-1} n_t}$$

Then

$$(A7.7) \qquad q_o = q_{o.j} < 2 \cdot 3^{\sum_{t=1}^{t=k-1} n_t}$$

In other words, $q_o$ in the final Cycle-k depends on the previous "k-1" cycles but cannot be larger than the power of 2 as expressed in Inequality (A7.6).

# Appendix 8 – Parallel paths

If two sequences follow the same pattern of upward and downward steps during "α" steps, we say that both sequences follow a parallel path. Two sequences can follow a parallel path during a finite number of steps, but any two sequences, at some point, follow different patterns.

## Lemma A8.1

Let $a_O^1$ - the initial element of a Sequence A made of "$i$" cycles – be:

$$(A6.1)\quad a_O^1 = 2^{\sum_{t=1}^{t=i}\alpha_t}\cdot Q_{O.i} - \frac{2^{3^{\sum_{t=1}^{t=i} n_t-1}\cdot j_{o.i}^{\delta}}+1}{3^{\sum_{t=1}^{t=i} n_t}}\cdot\sum_{s=1}^{i} 3^{\sum_{t=s+1}^{t=i} n_t}\cdot 2^{\sum_{t=1}^{t=s-1}\alpha_t}\cdot(3^{n_s}-2^{n_s})$$

$$Q_{O.i};j_{o.i}^{\delta}\in\{odd\}$$

$$(A6.2)\qquad 3^{\sum_{t=1}^{t=i} n_t-1}\cdot j_{o.i}^{\delta} > \sum_{t=1}^{t=i}\alpha_t$$

Let $b_O^1$- the initial element of a Sequence B made of "$i$" cycles – be:

$$(A8.1)\qquad b_O^1 = 2^{\sum_{t=1}^{t=i}\alpha_t}\cdot K_O + a_O^1$$

$$(A8.2)\quad b_O^1 = 2^{\sum_{t=1}^{t=i}\alpha_t}\cdot(Q_{O.i}+K_O) - \frac{2^{3^{\sum_{t=1}^{t=i} n_t-1}\cdot j_{o.i}^{\delta}}+1}{3^{\sum_{t=1}^{t=i} n_t}}\cdot\sum_{s=1}^{i} 3^{\sum_{t=s+1}^{t=i} n_t}\cdot 2^{\sum_{t=1}^{t=s-1}\alpha_t}\cdot(3^{n_s}-2^{n_s})$$

$$K_O\in\{odd\}$$

Sequence A and Sequence B follow parallel paths during "i" cycles – or $\sum_{t=1}^{t=i}\alpha_t$ downward steps. In other words, they both follow the same pattern of upward and downward steps until the end of Cycle "i".

Please note that $a_O^1$ and $b_O^1$ share the same parameters:

$$n_t,\alpha_t$$

for every cycle and the parameter $j_{o.i}^{\delta}$. The only difference is the addend:

$$2^{\sum_{t=1}^{t=i}\alpha_t}\cdot K_O$$

## Proof

We start with Identity (A8.1)

$$(A8.1) \qquad b_O^1 = a_O^1 + 2^{\sum_{t=1}^{t=i} \alpha_t} \cdot K_O$$

Since the addend is an even number, it won't affect the parity of $b_O^1$ as long as this addend remains even. As a result, Sequence A and Sequence B follow the same upward and downward pattern "i" cycles or $\sum_{t=1}^{t=i} \alpha_t$ downward steps.

We start with Cycle 1. The final element of the Cycle 1 of Sequence B after $n_1$ upward steps and $\alpha_1$ downward steps is:

$$b_F^1 = a_F^1 + 3^{n_1} \cdot 2^{\sum_{t=2}^{t=i} \alpha_t} \cdot K_O$$

This final element becomes the initial element of Cycle 2. This is the case for Sequence B

$$b_F^1 = b_O^2$$

and for Sequence A

$$a_F^1 = a_O^2$$

As a result, the initial element of Cycle 2 of Sequence B is:

$$b_O^2 = a_O^2 + 2^{\sum_{t=2}^{t=i} \alpha_t} \cdot 3^{n_1} \cdot K_O$$

The addend remains an even number, therefore it does not affect the parity of $b_O^2$. As a result, the final element of Cycle 2 of Sequence B is:

$$b_F^2 = a_F^2 + 2^{\sum_{t=3}^{t=i} \alpha_t} \cdot 3^{n_1+n_2} \cdot K_O$$

We repeat the same process during "i" cycles and obtain

$$b_F^i = a_F^i + 3^{\sum_1^i n_s} \cdot K_O$$

At this stage the addend is no longer an even number, therefore it affects the parity of $b_F^i$, and as a result, Sequence A and Sequence B no longer follow the same path.

Which is what we wanted to prove.

**Example:**

Let

$$a_o = 57$$

$$\alpha = 14$$

$$K_o = 1$$

then

$$b_o = 57 + 2^{14} = 16441$$

Sequence A and Sequence B follow a parallel path during 14 downward steps, until sequence A reaches the number 26 and Sequence B reaches the number 6587. At this point, paths follow different patterns.

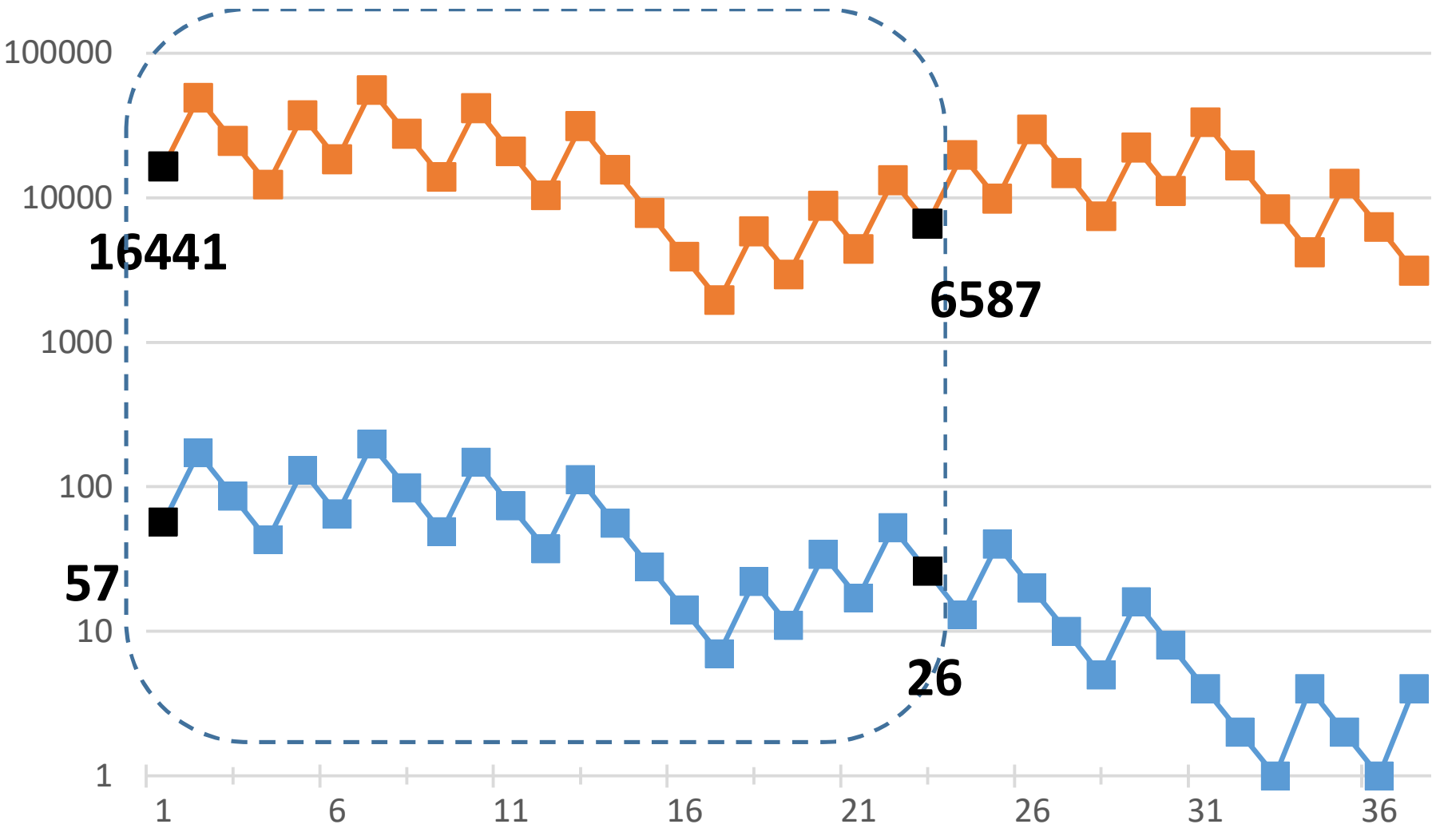

Note: The graph is in logarithm scale to facilitate the visualization of upward and downward steps.

## Lemma A8.2

Let

$$a_O^1 < 2^\alpha$$

This initial element generates a Sequence A that follows a pattern of upward and downward steps during the first $\alpha$ downward steps.

$a_O^1$ is the only number smaller than $2^\alpha$ that generates a sequence that follows such pattern during $\alpha$ downward steps.

## Proof

We start with Equation (A8.1) above

$$(A8.1) \quad b_O^1 = a_O^1 + 2^\alpha \cdot K_o$$

$$a_O^1 < 2^\alpha$$

and make

$$K_o = 1$$

then equation (A8.1) can be written as

$$b_O^1 = a_O^1 + 2^\alpha$$

Following Lemma A8.1, since $K_o = 1$, then $b_O$ is the smallest number that follows a parallel path to $a_O$ during $\alpha$ downward steps

Nevertheless, from the equation above we know that

$$b_O^1 > 2^\alpha$$

# Appendix 9 – Miscellaneous

In this appendix we will prove several lemmas that are necessary for the final demonstration.

**Lemma A9.1**

Let $a_O^1$ and $b_O^1$ two odd numbers defined as shown below. Please note that both numbers share all parameters except for

$$Q_{O.i}^{A} \neq Q_{O.i}^{B}$$

$$w_A \neq w_B$$

$$(A9.1)\quad a_O^1 = 2^{\sum_{t=1}^{t=i}\alpha_t} \cdot Q_{O.i}^{A} - \frac{2^{w_A}+1}{3^{\sum_{t=1}^{t=i} n_t}} \cdot \sum_{s=1}^{i} 3^{\sum_{t=s+1}^{t=i} n_t} \cdot 2^{\sum_{t=1}^{t=s-1}\alpha_t} \cdot (3^{n_s} - 2^{n_s})$$

$$(A9.2)\quad b_O^1 = 2^{\sum_{t=1}^{t=i}\alpha_t} \cdot Q_{O.i}^{B} - \frac{2^{w_B}+1}{3^{\sum_{t=1}^{t=i} n_t}} \cdot \sum_{s=1}^{i} 3^{\sum_{t=s+1}^{t=i} n_t} \cdot 2^{\sum_{t=1}^{t=s-1}\alpha_t} \cdot (3^{n_s} - 2^{n_s})$$

$$w_A = 3^{\sum_{t=1}^{t=i} n_t - 1} \cdot j_{o.i}^{\delta} > \sum_{t=1}^{t=i} \alpha_t$$

$$w_B = 3^{\sum_{t=1}^{t=k} n_t - 1} \cdot j_{o.i}^{\delta} > \sum_{t=1}^{t=i} \alpha_t$$

$$w_B > w_A$$

$$(A9.3)\quad Q_{O.i}^{B} = Q_{O.i}^{A} + 2^{w_A - \sum_{t=1}^{t=i}\alpha_t} \cdot \frac{2^{w_B - w_A} - 1}{3^{\sum_{t=1}^{t=i} n_t}} \cdot \sum_{s=1}^{i} 3^{\sum_{t=s+1}^{t=i} n_t} \cdot 2^{\sum_{t=1}^{t=s-1}\alpha_t} \cdot (3^{n_s} - 2^{n_s})$$

$$all\ numbers\ \in \mathbb{N}$$

Then $a_O^1$ and $b_O^1$ represent the same odd number:

$$a_O^1 = b_O^1$$

**Proof**

We start with Identity (A9.1)

$$(A9.1)\quad a_O^1 = 2^{\sum_{t=1}^{t=i}\alpha_t}\cdot Q_{O.i}^A - \frac{2^{w_A}+1}{3^{\sum_{t=1}^{t=i}n_t}}\cdot\sum_{s=1}^{i}3^{\sum_{t=s+1}^{t=i}n_t}\cdot 2^{\sum_{t=1}^{t=s-1}\alpha_t}\cdot(3^{n_s}-2^{n_s})$$

In order to simplify the identity above we make the following change:

$$(A9.4)\quad \sum I = \sum_{s=k+1}^{i}3^{\sum_{t=s+1}^{t=k}n_t}\cdot 2^{\sum_{t=1}^{t=s-1}\alpha_t}\cdot(3^{n_s}-2^{n_s})$$

And Identity (A9.1) becomes

$$a_O^1 = 2^{\sum_{t=1}^{t=i}\alpha_t}\cdot Q_{O.i}^A - \frac{2^{w_A}+1}{3^{\sum_{t=1}^{t=i}n_t}}\cdot\sum I$$

We add and subtract the addend

$$\frac{2^{w_B}+1}{3^{\sum_{t=1}^{t=i}n_t}}\cdot\sum I$$

and obtain

$$a_O^1 = 2^{\sum_{t=1}^{t=i}\alpha_t}\cdot Q_{O.i}^A - \frac{2^{w_A}+1}{3^{\sum_{t=1}^{t=i}n_t}}\cdot\sum I + \frac{2^{w_B}+1}{3^{\sum_{t=1}^{t=i}n_t}}\cdot\sum I - \frac{2^{w_B}+1}{3^{\sum_{t=1}^{t=i}n_t}}\cdot\sum I$$

We combine the second and the third addend and eliminate the common addend

$$\frac{1}{3^{\sum_{t=1}^{t=i}n_t}}\cdot\sum I$$

and obtain

$$a_O^1 = 2^{\sum_{t=1}^{t=i}\alpha_t}\cdot Q_{O.i}^A - \frac{2^{w_A}-2^{w_B}}{3^{\sum_{t=1}^{t=i}n_t}}\cdot\sum I - \frac{2^{w_B}+1}{3^{\sum_{t=1}^{t=i}n_t}}\cdot\sum I$$

We extract the common factor $2^{w_A}$ from the second addend

$$a_O^1 = 2^{\sum_{t=1}^{t=i}\alpha_t}\cdot Q_{O.i}^A + 2^{w_A}\cdot\frac{2^{w_B-w_A}-1}{3^{\sum_{t=1}^{t=i}n_t}}\cdot\sum I - \frac{2^{w_B}+1}{3^{\sum_{t=1}^{t=i}n_t}}\cdot\sum I$$

We extract the factor $2^{\sum_{t=1}^{t=i}\alpha_t}$ from the second addend and combine first and second addends

$$a_O^1 = 2^{\sum_{t=1}^{t=i}\alpha_t}\cdot\left(Q_{O.i}^A + 2^{w_A-\sum_{t=1}^{t=i}\alpha_t}\cdot\frac{2^{w_B-w_A}-1}{3^{\sum_{t=1}^{t=i}n_t}}\cdot\sum I\right) - \frac{2^{w_B}+1}{3^{\sum_{t=1}^{t=i}n_t}}\cdot\sum I$$

We use Identity (A9.3) above and make the following change:

$$(A9.3)\quad Q_{O.i}^{B} = Q_{O.i}^{A} + 2^{w_A - \sum_{t=1}^{t=i} \alpha_t} \cdot \frac{2^{w_B - w_A} - 1}{3^{\sum_{t=1}^{t=i} n_t}} \cdot \sum I$$

and obtain

$$a_O^1 = 2^{\sum_{t=1}^{t=i} \alpha_t} \cdot Q_{O.i}^{B} - \frac{2^{w_B} + 1}{3^{\sum_{t=1}^{t=i} n_t}} \cdot \sum I$$

We undo the change made with Identity (A9.4)

$$a_O^1 = 2^{\sum_{t=1}^{t=i} \alpha_t} \cdot Q_{O.i}^{B} - \frac{2^{w_B} + 1}{3^{\sum_{t=1}^{t=i} n_t}} \cdot \sum_{s=k+1}^{i} 3^{\sum_{t=s+1}^{t=k} n_t} \cdot 2^{\sum_{t=1}^{t=s-1} \alpha_t} \cdot (3^{n_s} - 2^{n_s})$$

The right side of the identity above is the same as in Identity (A9.2). Therefore

$$a_O^1 = b_O^1 = 2^{\sum_{t=1}^{t=i} \alpha_t} \cdot Q_{O.i}^{B} - \frac{2^{w_B} + 1}{3^{\sum_{t=1}^{t=i} n_t}} \cdot \sum_{s=k+1}^{i} 3^{\sum_{t=s+1}^{t=k} n_t} \cdot 2^{\sum_{t=1}^{t=s-1} \alpha_t} \cdot (3^{n_s} - 2^{n_s})$$

which is what we wanted to prove.

## Lemma A9.2

Let $a_O^1$ - the initial element of a Sequence A made of "$i$" cycles – and $a_F^i$ the final element of this sequence after "i" cycles.

If

$$a_F^i \geq a_O^1$$
$$\forall i$$

then

$$(A9.5) \qquad 2^{\sum_1^i \alpha_t} < \left(1+\frac{1}{a_O^1}\right)^{\sum_1^i n_t} \cdot 3^{\sum_1^i n_t}$$
$$\forall i$$
$$i \in \mathbb{N}$$

## Proof

Let's assume the opposite, in other words, let's assume it exists an "i" for which:

$$(A9.6) \qquad 2^{\sum_1^i \alpha_t} > \left(1+\frac{1}{a_O^1}\right)^{\sum_1^i n_t} \cdot 3^{\sum_1^i n_t}$$

but the inequality holds for any "s" smaller than "i":

$$(A9.7) \qquad 2^{\sum_1^s \alpha_t} < \left(1+\frac{1}{a_O^1}\right)^{\sum_1^i n_t} \cdot 3^{\sum_1^i n_t}$$
$$\forall s < i \ \in \mathbb{N}$$

We start with the Identity (A6.4) from Appendix 6:

$$(A6.4) \qquad 2^{\sum_{t=1}^{t=i} \alpha_t} \cdot a_F^i - 3^{\sum_{t=1}^{t=i} n_t} \cdot a_O^1 = \sum_{s=1}^{i} 3^{\sum_{t=s+1}^{t=i} n_t} \cdot 2^{\sum_{t=1}^{t=s-1} \alpha_t} \cdot (3^{n_s} - 2^{n_s})$$

We divide all elements by $a_O^1 \cdot 2^{\sum_{t=1}^{t=i} \alpha_t}$

$$\frac{a_F^i}{a_O^1} = \frac{3^{\sum_{t=1}^{t=i} n_t}}{2^{\sum_{t=1}^{t=i} \alpha_t}} + \frac{\sum_{s=1}^{s=i} 3^{\sum_{t=s+1}^{t=i} n_t} \cdot 2^{\sum_{t=1}^{t=s-1} \alpha_t} \cdot (3^{n_s} - 2^{n_s})}{a_O^1 \cdot 2^{\sum_{t=1}^{t=i} \alpha_t}}$$

We substitute the denominator of the first fraction for something smaller (see Inequality (A9.6)) and obtain a larger number:

$$\frac{a_F^i}{a_O^1} < \frac{3^{\sum_{t=1}^{t=i} n_t}}{\left(1+\frac{1}{a_O^1}\right)^{\sum_1^i n_t} \cdot 3^{\sum_1^i n_t}} + \frac{\sum_{s=1}^{s=i} 3^{\sum_{t=s+1}^{t=i} n_t} \cdot 2^{\sum_{t=1}^{t=s-1} \alpha_t} \cdot (3^{n_s} - 2^{n_s})}{a_O^1 \cdot 2^{\sum_{t=1}^{t=i} \alpha_t}}$$

We simplify the first fraction by eliminating the common factor $3^{\sum_{t=1}^{t=i} n_t}$ from the numerator and the denominator:

$$\frac{a_F^i}{a_O^1} < \frac{1}{\left(1+\frac{1}{a_O^1}\right)^{\sum_1^i n_t}} + \frac{\sum_{s=1}^{s=i} 3^{\sum_{t=s+1}^{t=i} n_t} \cdot 2^{\sum_{t=1}^{t=s-1} \alpha_t} \cdot (3^{n_s} - 2^{n_s})}{a_O^1 \cdot 2^{\sum_{t=1}^{t=i} \alpha_t}}$$

We use Inequality (A9.7) and substitute $2^{\sum_{t=1}^{t=s-1} \alpha_t}$ in the numerator for something larger $\left(1+\frac{1}{a_O^1}\right)^{\sum_1^{s-1} n_t} \cdot 3^{\sum_1^{s-1} n_t}$ and the inequality holds.

$$\frac{a_F^i}{a_O^1} < \frac{1}{\left(1+\frac{1}{a_O^1}\right)^{\sum_1^i n_t}} + \frac{\sum_{s=1}^{s=i} 3^{\sum_{t=s+1}^{t=i} n_t} \cdot \left(1+\frac{1}{a_O^1}\right)^{\sum_1^{s-1} n_t} \cdot 3^{\sum_1^{s-1} n_t} \cdot (3^{n_s} - 2^{n_s})}{a_O^1 \cdot 2^{\sum_{t=1}^{t=i} \alpha_t}}$$

We use Inequality (A9.6) and substitute $2^{\sum_{t=1}^{t=i} \alpha_t}$ in the denominator for something smaller $\left(1+\frac{1}{a_O^1}\right)^{\sum_1^i n_t} \cdot 3^{\sum_1^i n_t}$ and the inequality holds.

$$\frac{a_F^i}{a_O^1} < \frac{1}{\left(1+\frac{1}{a_O^1}\right)^{\sum_1^i n_t}} + \frac{\sum_{s=1}^{s=i} 3^{\sum_{t=s+1}^{t=i} n_t} \cdot \left(1+\frac{1}{a_O^1}\right)^{\sum_1^{s-1} n_t} \cdot 3^{\sum_1^{s-1} n_t} \cdot (3^{n_s} - 2^{n_s})}{a_O^1 \cdot \left(1+\frac{1}{a_O^1}\right)^{\sum_1^i n_t} \cdot 3^{\sum_1^i n_t}}$$

We eliminate the common factors $3^{\sum_{t=1}^{t=i} n_t}$ and $\left[1+\frac{1}{a_O^1}\right]^{\sum_1^{s-1} n_t}$ from the numerator and the denominator above:

$$\frac{a_F^i}{a_O^1} < \frac{1}{\left[1+\frac{1}{a_O^1}\right]^{\sum_1^i n_t}} + \frac{1}{a_O^1} \cdot \sum_{s=1}^{s=i} \frac{1}{\left[1+\frac{1}{a_O^1}\right]^{\sum_s^i n_t}} \cdot \frac{3^{n_s} - 2^{n_s}}{3^{n_s}}$$

We substitute $\frac{3^{n_s}-2^{n_s}}{3^{n_s}}$ for its maximum value which is $1$ for $n_s = \infty$ the inequality holds.

$$\frac{a_F^i}{a_O^1} < \frac{1}{\left[1+\frac{1}{a_O^1}\right]^{\sum_1^i n_t}} + \frac{1}{a_O^1} \cdot \sum_{s=1}^{s=i} \frac{1}{\left[1+\frac{1}{a_O^1}\right]^{\sum_s^i n_t}}$$

The summation above is a geometric progression with missing elements. If we add the missing elements the inequality holds:

$$\sum_{s=1}^{s=i} \frac{1}{\left[1+\frac{1}{a_O^1}\right]^{\sum_s^i n_t}} < \sum_{m=1}^{m=\sum_1^i n_t} \frac{1}{\left[1+\frac{1}{a_O^1}\right]^m}$$

Therefore, inequality above becomes

$$\frac{a_F^i}{a_O^1} < \frac{1}{\left[1+\frac{1}{a_O^1}\right]^{\sum_1^i n_t}} + \frac{1}{a_O^1} \cdot \sum_{m=1}^{m=\sum_1^i n_t} \left(\frac{1}{1+\frac{1}{a_O^1}}\right)^m$$

Since

$$\frac{1}{1+\frac{1}{a_O^1}} = \frac{a_O^1}{a_O^1+1}$$

We simplify the inequality above and obtain:

$$\frac{a_F^i}{a_O^1} < \left(\frac{a_O^1}{a_O^1+1}\right)^{\sum_1^i n_t} + \frac{1}{a_O^1} \cdot \sum_{m=1}^{m=\sum_1^i n_t} \left(\frac{a_O^1}{a_O^1+1}\right)^m$$

The summation above is a geometric progression that is solved as follows:

$$\sum_1^m r^m = \frac{r^{m+1}-r}{r-1}$$

We solve the geometric progression above and obtain:

$$\frac{a_F^i}{a_O^1} < \left(\frac{a_O^1}{a_O^1+1}\right)^{\sum_1^i n_t} + \frac{1}{a_O^1} \cdot \frac{\left(\frac{a_O^1}{a_O^1+1}\right)^{\sum_1^i n_t+1} - \left(\frac{a_O^1}{a_O^1+1}\right)}{\left(\frac{a_O^1}{a_O^1+1}\right)-1}$$

We extract the common factor

$$\frac{a_O^1}{a_O^1+1}$$

from the fraction of the geometric progression and obtain

$$\frac{a_F^i}{a_O^1} < \left(\frac{a_O^1}{a_O^1+1}\right)^{\Sigma_1^i n_t} + \frac{1}{a_O^1}\cdot\left(\frac{a_O^1}{a_O^1+1}\right)\cdot\frac{\left(\frac{a_O^1}{a_O^1+1}\right)^{\Sigma_1^i n_t}-1}{\left(\frac{a_O^1}{a_O^1+1}\right)-1}$$

We simplify the above and obtain:

$$\frac{a_F^i}{a_O^1} < \left(\frac{a_O^1}{a_O^1+1}\right)^{\Sigma_1^i n_t} + 1 - \left(\frac{a_O^1}{a_O^1+1}\right)^{\Sigma_1^i n_t}$$

We simplify and obtain:

$$\frac{a_F^i}{a_O^1} < 1$$

Which it cannot be since the initial premise was that $a_F^i \geq a_O^1$ . Therefore, there is a contradiction. As a result, Inequality (A9.5) is not possible and

$$2^{\Sigma_1^i \alpha_t} < \left(1+\frac{1}{a_O^1}\right)^{\Sigma_1^i n_t}\cdot 3^{\Sigma_1^i n_t}$$

$$\forall i$$

Which is what we wanted to prove.

## Lemma A9.3

Let $a_O^1$ - the initial element of a Sequence A made of "$i$" cycles – and $a_F^i$ the final element of this sequence after "i" cycles.

Sequence A does not converge

$$a_F^i \geq a_O^1$$

Then, the maximum value of the summation

$$(A9.8) \qquad \sum I = \sum_{s=1}^{s=i} 3^{\sum_{t=s+1}^{t=i} n_t} \cdot 2^{\sum_{t=1}^{t=s-1} \alpha_t} \cdot (3^{n_s} - 2^{n_s})$$

is given by inequality

$$(A9.9) \qquad \left(\sum I\right)_{MAX} < 3^{\sum_{t=1}^{t=i} n_t} \cdot (a_O^1 + 1) \cdot \left[\left(1 + \frac{1}{a_O^1}\right)\right]^{\sum_1^i n_t}$$

## Proof

We start with the summation

$$(A9.8) \qquad \sum I = \sum_{s=1}^{s=i} 3^{\sum_{t=s+1}^{t=i} n_t} \cdot 2^{\sum_{t=1}^{t=s-1} \alpha_t} \cdot (3^{n_s} - 2^{n_s})$$

Since

$$a_F^i \geq a_O^1$$

then lemma A9.2 is true

$$(A9.5) \qquad 2^{\sum_1^i \alpha_t} < \left(1 + \frac{1}{a_O^1}\right)^{\sum_1^i n_t} \cdot 3^{\sum_1^i n_t}$$

$$\forall i$$

We apply Inequality (A9.5) to Identity (A9.8) and obtain a larger number.

$$\sum I < \sum_{s=1}^{s=i} 3^{\sum_{t=s+1}^{t=i} n_t} \cdot \left(1 + \frac{1}{a_O^1}\right)^{\sum_1^{s-1} n_t} \cdot 3^{\sum_1^{s-1} n_t} \cdot (3^{n_s} - 2^{n_s})$$

We extract $3^{\sum_{t=1}^{t=i} n_t}$ from all addends of the summation

$$\sum I < 3^{\sum_{t=1}^{t=i} n_t} \cdot \sum_{s=1}^{s=i} \left(1 + \frac{1}{a_O^1}\right)^{\sum_1^{s-1} n_t} \cdot \left(1 - \frac{2^{n_s}}{3^{n_s}}\right)$$

If we substitute $\left(1 - \frac{2^{n_s}}{3^{n_s}}\right)$ for its maximum value which is $1$ for $n_s = \infty$ the inequality holds

$$(A9.10) \qquad \sum I < 3^{\sum_{t=1}^{t=i} n_t} \cdot \sum_{s=1}^{s=i} \left(1 + \frac{1}{a_O^1}\right)^{\sum_1^{s-1} n_t}$$

The summation above is a geometric progression with missing elements. If we add the missing elements the inequality holds:

$$\sum_{s=1}^{s=i} \left(1 + \frac{1}{a_O^1}\right)^{\sum_1^{s-1} n_t} < \sum_{m=1}^{m=\sum_1^i n_t} \left(1 + \frac{1}{a_O^1}\right)^m$$

Therefore, Inequality (A9.10) becomes

$$(A9.11) \quad \sum I < 3^{\sum_{t=1}^{t=i} n_t} \cdot \sum_{m=1}^{m=\sum_1^i n_t} \left(1 + \frac{1}{a_O^1}\right)^m$$

The summation above is a geometric progression that is solved as follows:

$$\sum_1^m r^m = \frac{r^{m+1} - r}{r - 1}$$

We apply the above to Inequality (A9.11)

$$\sum I < 3^{\sum_{t=1}^{t=i} n_t} \cdot \frac{\left(1 + \frac{1}{a_O^1}\right)^{\sum_{t=1}^{t=i} n_t + 1} - \left(1 + \frac{1}{a_O^1}\right)}{\left(1 + \frac{1}{a_O^1}\right) - 1}$$

We simplify the fraction above

$$\sum I < 3^{\sum_{t=1}^{t=i} n_t} \cdot (a_O^1 + 1) \cdot \left[\left(1 + \frac{1}{a_O^1}\right)^{\sum_{t=1}^{t=i} n_t} - 1\right]$$

Therefore

$$\left(\sum I\right)_{MAX} < 3^{\sum_{t=1}^{t=i} n_t} \cdot (a_O^1 + 1) \cdot \left(1 + \frac{1}{a_O^1}\right)^{\sum_{t=1}^{t=i} n_t}$$

Which is what we wanted to prove

**Lemma A9.4**

The Identity below is true

$$(A9.12) \quad \sum_{s=k+1}^{m} 3^{\sum_{t=s+1}^{t=m} 1} \cdot 2^{\sum_{t=k+1}^{t=s-1} 2} \cdot (3^1 - 2^1) = 4^{(m-k)} - 3^{(m-k)}$$

## Proof

We start with the summation above

$$\sum_{s=k+1}^{m} 3^{\sum_{t=s+1}^{t=m} 1} \cdot 2^{\sum_{t=k+1}^{t=s-1} 2} \cdot (3^1 - 2^1) = 3^{m-(k+1)} + 3^{m-(k+2)} \cdot 4 + 3^{m-(k+3)} \cdot 4^2 + \cdots + 4^{m-(k+1)}$$

If we multiply the summation by 1, the value of the summation doesn't change. Therefore, we multiply by

$$1 = (4-3)$$

And obtain:

$$(4-3) \cdot \sum \quad = (4-3) \cdot \left(3^{m-(k+1)} + 3^{m-(k+2)} \cdot 4 + 3^{m-(k+3)} \cdot 4^2 + \cdots + 4^{m-(k+1)}\right)$$

Therefore

$$\sum \quad = 3^{m-(k+1)} \cdot 4 + 3^{m-(k+2)} \cdot 4^2 + 3^{m-(k+3)} \cdot 4^3 + \cdots + 4^{m-k}$$

$$-\left(3^{m-k} + 3^{m-(k+1)} \cdot 4 + 3^{m-(k+2)} \cdot 4^2 + \cdots + 3 \cdot 4^{m-(k+1)}\right)$$

We eliminate the common elements and obtain:

$$\sum \quad = 4^{m-k} - 3^{m-k}$$

Which is what we wanted to prove.